\documentstyle[12pt, twoside]{article}
\input amssym.def
\input amssym.tex 

\newenvironment{demo}[1]%
{\vskip-\lastskip\medskip
  \noindent
  {\em #1.}\enspace
  }%
{\qed\par\medskip
  }

\newcommand{\qed}{
  \strut\hfill
  \mbox{$\Box$}
  }

{\vskip-\lastskip
  \begin{trivlist}
    \item[]
      {\bf Axiom #1 }
      }%
{\end{trivlist}
  }

\newtheorem{theorem}{Theorem}[section]

\newtheorem{corollary}{Corollary}[section]

\newtheorem{lemma}{Lemma}[section]

\newtheorem{convention}{Convention}[section]

\newtheorem{example}{Example}[section]

\newtheorem{remark}{Remark}[section]

\newtheorem{definition}{Definition}[section]

\newtheorem{proposition}{Proposition}[section]

\newcommand{\Binf}{ {\overline{b}}^-_{\infty} }

\newcommand{\Binftwo}{ {\overline{b}}^+_{\infty} }

\newcommand{\Binfm}{ {\overline{b}}_{\infty}^{- [m]} }

\newcommand{\Binfmtwo}{ {\overline{b}}_{\infty}^{+ [m]} }

\newcommand{\Cinf}{ {\overline{c}}_{\infty} }

\newcommand{\Cinfm}{ {\overline{c}}_{\infty}^{[m]} }

\newcommand{\Dinf}{ {\overline{d}}_{\infty} }

\newcommand{\Dinfm}{ {\overline{d}}_{\infty}^{[m]} }

\newcommand{\binf}{ b_{\infty} }

\newcommand{\binftwo}{ \tilde{b}_{\infty} }

\newcommand{\binfm}{ b_{\infty}^{[m]} }

\newcommand{\binfmtwo}{ \tilde{b}_{\infty}^{[m]} }

\newcommand{\binfmpm}{  {\binf}^{[m]}_{\pm} }

\newcommand{\cinf}{ c_{\infty} }

\newcommand{\cinfm}{ c_{\infty}^{[m]} }

\newcommand{\cinfmpm}{  {\cinf}^{[m]}_{\pm} }

\newcommand{\dinf}{ d_{\infty} }

\newcommand{\dinfm}{ d_{\infty}^{[m]} }

\newcommand{\dinfmpm}{  {\dinf}^{[m]}_{\pm} }

\newcommand{\B}{ {\cal D}^{-} }

\newcommand{\Bpm}{ {\cal D}^{\pm} }

\newcommand{\Bp}{ {\cal D}^{+} }

\newcommand{\BO}{ {\cal D}^{ { \cal O}, -} }

\newcommand{\BOpm}{ {\cal D}^{ { \cal O}, \pm} }

\newcommand{\BOp}{ {\cal D}^{ { \cal O}, +} }

\newcommand{\C}{ {\Bbb C} }

\newcommand{\co}{ \cal O}

\newcommand{\D}{ {\cal D} }

\newcommand{\DO}{ {\cal D}^{ { \cal O}} }

\newcommand{\ch}{ \mbox{ch} }

\newcommand{\cwe}{\Bbb C[w]^{(0)}}

\newcommand{\cwo}{\Bbb C[w]^{(1)}}

\newcommand{\End}{ \mbox{End } }

\newcommand{\Fhalf}{ {\cal F}^{\bigotimes \hf} }

\newcommand{\Fhalfminusl}{ {\cal F}^{\bigotimes -l + \hf} }

\newcommand{\Fl}{ {\cal F}^{\bigotimes l} }

\newcommand{\Flhalf}{ {\cal F}^{\bigotimes l + \hf} }

\newcommand{\Flminushalf}{ {\cal F}^{\bigotimes l - \hf} }

\newcommand{\Fminushalf}{ {\cal F}^{\bigotimes  - \hf} }

\newcommand{\Fminusl}{ {\cal F}^{\bigotimes  - l} }

\newcommand{\Fminuslhalf}{ {\cal F}^{\bigotimes  - l - \hf} }

\newcommand{\Fpairbb}{ ( Spin(2l+1), \binf ) }

\newcommand{\Fpairbbtwo}{ ( Spin(2l+1), \binftwo ) }

\newcommand{\Fpairbd}{ ( O (2l+1), \dinf ) }

\newcommand{\Fpaircc}{ (Sp (2l), \cinf ) }

\newcommand{\Fpairdb}{ ( Pin(2l), \binf ) }

\newcommand{\Fpairdbtwo}{ ( Pin(2l), \binftwo ) }

\newcommand{\Fpairospc}{ ({\frak {osp}} (1, 2l), \cinf ) }

\newcommand{\g}{ \frak g }

\newcommand{\glm}{ {\frak{gl}}^{[m]} }

\newcommand{\gl}{ {\frak{gl}} }

\newcommand{\glpm}{ \widehat{gl}_{\pm} }

\newcommand{\gm}{ {\frak{g}}^{[m]} }

\newcommand{\gvm}{ {\frak{g}}^{[\vec{m}]} }

\newcommand{\hB}{ {\widehat{\cal D}^{-} } }

\newcommand{\hBp}{ {\widehat{\cal D}^{+} } }

\newcommand{\hBpm}{ {\widehat{\cal D}^{\pm} } }

\newcommand{\hBO}{ {\widehat{\cal D}^{{ \cal O}, -} } }

\newcommand{\hBOpm}{ {\widehat{\cal D}^{{ \cal O}, \pm} } } 

\newcommand{\hD}{ {\widehat{\cal D}} }

\newcommand{\hf}{ \frac12 }

\newcommand{\hgl}{ \widehat{ \frak{gl} } }

\newcommand{\hglm}{ \widehat{ \frak{gl} }^{[m]} } 

\newcommand{\hL}{ \widehat{\Lambda} }

\newcommand{\hn}{ \hf + {\Bbb Z}_{+} }

\newcommand{\hphi}{ \widehat{\phi}_{1/2} }

\newcommand{\hphic}{ \widehat{\phi}_{ -1/2}}

\newcommand{\hphicm}{ \widehat{\phi}_{ -1/2}^{[m]}}

\newcommand{\hphio}{ \widehat{\phi}_0}

\newcommand{\hphim}{ \widehat{\phi}_{1/2}^{[m]}}

\newcommand{\hphiom}{ \widehat{\phi}_0^{[m]}}

\newcommand{\hphism}{ \widehat{\phi}_s^{[m]}}

\newcommand{\hphivsm}{ \widehat{\phi}_{ \vec{s} }^{ [ \vec{m}] } }

\newcommand{\hz}{ \hf+ \Bbb Z }

\newcommand{\Mpairbc}{ ( O(2l +1), \cinf ) }

\newcommand{\Mpaircd}{ ( Sp(2l), \dinf ) }

\newcommand{\Mpairdc}{ ( O(2l ), \cinf ) }

\newcommand{\Mpairospd}{ ({\frak {osp} } (1,2l), \dinf) }

\newcommand{\NN}{ z^{-n-1} }

\newcommand{\Oe}{ {\cal O}^{(0)}}

\newcommand{\Oo}{ {\cal O}^{(1)}}

\newcommand{\Pa}{ \cal P}

\newcommand{\phim}{\phi_{1/2}^{[m]}}

\newcommand{\phicm}{\phi_{ -1/2}^{[m]}}

\newcommand{\phiom}{\phi_0^{[m]}}

\newcommand{\phism}{\phi_s^{[m]}}

\newcommand{\phivsm}{ \phi_{ \vec{s} }^{ [ \vec{m}] } }

\newcommand{\Tr}{ \mbox{Tr} }

\newcommand{\Ts}{ \Theta_s }

\newcommand{\vac}{ |0 \rangle }

\newcommand{\W}{ {\cal W}_{1+\infty} }

\newcommand{\Wpairbb}{ ( Spin(2l+1), \hB ) }

\newcommand{\Wpairbbtwo}{ ( Spin(2l+1), \hBp ) }

\newcommand{\Wpairbc}{ ( O(2l+1), \hB ) }

\newcommand{\Wpairbd}{ ( O (2l+1), \hBp ) }

\newcommand{\Wpaircc}{ (Sp (2l), \hB ) }

\newcommand{\Wpaircd}{ ( Sp(2l), \hBp ) }

\newcommand{\Wpairdb}{ ( Pin(2l), \hB ) }

\newcommand{\Wpairdbtwo}{ ( Pin(2l), \hBp ) }

\newcommand{\Wpairdc}{ ( O(2l), \hB ) }

\newcommand{\Wpairdd}{ ( O(2l), \hBp ) }

\newcommand{\Wpairospc}{ ({\frak {osp}}(1, 2l), \hB ) }

\newcommand{\Wpairospd}{ ({\frak {osp} } (1,2l), \hBp) }

\newcommand{\Z}{ {\Bbb Z} }

\begin{document}
\title{ 
Quasifinite representations of classical Lie subalgebras of $\W$
  }
\author{
  Victor G. Kac${}^1$, $\;$
  Weiqiang Wang${}^2$\thanks{On leave from Department of Mathematics,
 Yale University, USA} $\;$
 and  $\;$ Catherine H. Yan${}^3$
\\
\\{\small ${}^1$Department of Mathematics}
\\{\small MIT, Cambridge, MA 02139, USA}
\\{\small E-mail: kac@math.mit.edu}
\\
\\{\small ${}^2$Max-Planck Institut f\"ur Mathematik}
\\{\small Gottfried-Claren Str. 26, 53225 Bonn, Germany}
\\{\small E-mail: wqwang@mpim-bonn.mpg.de}
\\
\\{\small ${}^3$Courant Institute of Mathematical Sciences}
\\{\small 251 Mercer Street, New York, NY 10012, USA}
\\{\small E-mail: yanhuaf@math1.cims.nyu.edu}
}

\date{}
\maketitle
\begin{abstract}
   We show that there are precisely two, up to conjugation,
 anti-involutions  $\sigma_{\pm}$
 of the algebra of differential operators on the circle
 preserving the principal gradation.
 We classify the irreducible quasifinite highest weight
 representations
 of the central extension $\hBpm$ of the Lie subalgebra of
 this algebra fixed by $ -\sigma_{\pm}$, and find the unitary ones.
 We realize them in terms of highest weight representations of
 the central extension of the Lie algebra of infinite matrices 
 with finitely many non-zero diagonals
 over the algebra ${\Bbb C} [u] / (u^{m+1} )$ and its classical
 Lie subalgebras of $B$, $ C$ and $D$ types.
 Character formulas for {\em positive primitive}
 representations of $\hBpm$ (including all the unitary ones) 
 are obtained. We also realize a class of primitive representations of
 $\hBpm$ in terms of free fields and establish a number of duality results
 between these primitive representations
 and finite-dimensional irreducible representations of
 finite-dimensional Lie groups and supergroups. 
 We show that the vacuum module $V_c$
 of $\hBp$ carries a vertex algebra structure and establish 
 a relationship between  $V_c$ for $c \in \hf \Z$ and $\cal W$-algebras.
\end{abstract}
\tableofcontents
\setcounter{section}{-1}
\section{Introduction}
    The systematic study of quasifinite highest weight modules 
of the universal central extension $\hD$ of the Lie algebra of
differential operators over the circle (described first in \cite{KP})
was initiated by
Kac and Radul in \cite{KR1} and further studied
in \cite{Ma, FKRW, AFMO, KR2, W1} and many others. 
The Lie algebra $\hD$ is also known in the literature as $\W$ \cite{PRS} 
and as one of the universal $\cal W$-algebras
(cf.  \cite{BS} \cite{BEH$^3$} and references therein),
and it has various connections with
conformal field theory, the quantum Hall effect \cite{CTZ} etc. 

The difficulty in understanding representation theory 
of a Lie algebra of this sort 
is that although $\hD$ admits a natural
``principal'' $\Z$-gradation and
thus the associated triangular decomposition, each of the graded subspaces
of $\hD$ is still infinite-dimensional in contrast to the more familiar
cases such as the Virasoro algebra and Kac-Moody algebras.
The study of the highest weight modules of $\hD$ with 
the finiteness requirement on the dimensions
of their graded subspaces (which we will refer to
as quasifiniteness condition throughout our paper) thus becomes a 
highly non-trivial problem. 
By analyzing for which
parabolic subalgebras of $\hD$ the corresponding generalized Verma
modules are quasifinite, Kac and Radul \cite{KR1} were able to
give an elegant characterization of quasifinite highest weight
$\hD$-modules in terms of a certain generating function of
highest weights. They constructed all these $\hD$-modules in terms
of representations of the Lie algebra $\hglm$ which is
the central extension of the Lie algebra $\glm$
of infinite matrices with finitely many non-zero diagonals
taking values in the truncated polynomial
algebra $R_m = {\Bbb C} [u] / (u^{m+1} )$.
They also classified all such $\hD$-modules which are unitary.

It is well known that the Lie algebra $\hgl$ \cite{KP, DJKM}, 
a special case of $\hglm$ with $ m =0$, has Lie subalgebras of B, C, D types
as in the finite dimensional case \cite{K}.
So a natural question which arises here is: what are the
Lie subalgebras of $\hD$ which correspond to the classical Lie
subalgebras of $\hgl$ of B, C, D types. 
It is the goal of this paper to give a complete answer
to this question and to present the representation theory of
these Lie subalgebras of $\hD$. 

We show that there are two, up to conjugation,
anti-involutions of $\hD$, denoted by $\sigma_{\pm}$,
which preserve the principal $\Z$-gradation.
We thus obtain two different Lie subalgebras,
denoted by $\widehat{\cal D}^{\pm}$, 
fixed by $ - \sigma_{\pm}$. The graded subspaces of
$\widehat{\cal D}^{\pm}$ are still infinite-dimensional.
We regard $\hBp$ to be more fundamental than $\hB$ as we shall
show that the vacuum module
of $\hBp$ carries a canonical vertex algebra structure.

We first give a description of parabolic subalgebras of $\hBpm$ in terms
of the so-called characteristic polynomials 
$\{b_k (w)\}_{ k \in \Bbb N}$. It turns out that 
the quasifiniteness of a generalized Verma module 
is completely determined by the non-vanishing
of the first characteristic polynomial $b_1 (w)$. The corresponding
condition in defining the generalized Verma module 
induced from a parabolic subalgebra $\Pa$ leads to a characterization
of highest weights for which the corresponding $\hBpm$-module
is quasifinite in terms of a generating function $\Delta (x)$. 
This is the content of our Theorem \ref{th_delta}.

Next for each $s \in \C$
we construct a Lie algebra homomorphism
$\hphism$ from $\hBpm$ to $\hglm$. While $\hBpm$ is the central
extension of a Lie algebra
of differential operators with polynomial coefficients,
it is important to consider a Lie algebra $\hBOpm$ which is 
an analytic completion of $\hBpm$ and to extend the  homomorphisms
$\hphism$ to $\hBOpm$. 
It turns out that for $s \notin \Z/2$ the homomorphism $\hphism$ 
from $\hBOpm$ to $\hglm$
is surjective, but for $s \in \Z/2$, this is no longer true.
The image for $s \in \Z/2$ turns out to be various classical Lie
subalgebras of $\hglm$.

More generally we define a family of Lie algebra homomorphisms
$\hphivsm$ from $\hBpm$ to $\gvm$ for a vector 
$\vec{m} = (m_1, \ldots, m_N) \in \Z_{+}$ and a vector 
$\vec{s} = (s_1, \ldots , s_N) \in \C^N$ 
satisfying a certain condition ($\star \pm$), where $\gvm$ is a direct
sum of Lie algebras $\binf^{[m_i]}$, $\binftwo^{[m_i]}$, $\cinf^{[m_i]}$,
$\dinf^{[m_i]}$, and $\hgl^{[m_i]}$ satisfying a certain consistency 
condition. These homomorphisms $\hphivsm$ again extend to $\hBOpm$ and
then become surjective.
The principal $\Z$-gradations on $\binf^{[m_i]}$, $\cinf^{[m_i]}$
and $\hgl^{[m_i]}$, and the specialized $\Z$-gradation on
$\dinf^{[m_i]}$ of type $(2, 1, 1, \cdots)$
induce one on $\gvm$. The homomorphism $\hphivsm$
matches the induced $\Z$-gradation on $\gvm$ with the 
principal $\Z$-gradation on $\hBpm$.

Quasifinite highest weight modules (QHWM's) of $\hglm$ 
and its classical subalgebras
are better understood (particularly when $m =0$).
Irreducible QHWM's over $\gvm$
can be regarded as modules over $\hBpm$ via the homomorphisms $\hphivsm$.
Our Theorem~\ref{th_pullback} asserts that they are
irreducible over $\hBpm$ as well.
Our Theorem~\ref{th_realization} asserts that all 
irreducible quasifinite highest weight modules over $\hBpm$
can be realized in this way.

There is a natural anti-involution $\omega$ on $\hBpm$
which is compatible with the standard Cartan involution
on $\hgl$ via the Lie algebra homomorphism $\widehat{\phi}_s$.
We show that unitary irreducible quasifinite highest weight modules 
over $\hBpm$ with respect to $\omega$ are the pullbacks
of those over $\hgl$ via the homomorphisms $\widehat{\phi}_s$ for real
vectors $\vec{s}$. We give 
a simple characterization for these unitary $\hBpm$-modules
in terms of the generating function
$\Delta (x)$ for the highest weights with respect to $\hBpm$
and the central charge $c$. This is our Theorem~\ref{th_unit}.
The proofs of Theorems~\ref{th_delta}, \ref{th_pullback}, \ref{th_realization}
and \ref{th_unit} are similar to those in \cite{KR1} for $\hD$.

The $q$-character formulas in accordance with
the $\Z$-gradation for unitary quasifinite highest weight
modules of $ \binf$, $\cinf$ and $\dinf$ are worked out explicitly. 
This leads to the $q$-character formulas
for the so-called {\em positive primitive} quasifinite $\hBpm$-modules, 
which include all unitary quasifinite $\hBpm$-modules.

In the second half of this paper we study 
primitive quasifinite representations
of $\hBpm$ in terms of free bosonic and fermionic fields. 
Recall that representations of $\binf$ or $\cinf$
(resp. $\binftwo$ and $\dinf$)
were realized in various Fock spaces \cite{DJKM1, W2}. 
A number of duality results in various Fock spaces between
finite dimensional irreducible representations of a
finite dimensional Lie group or superalgebra and 
quasifinite representations of $\binf$ or $\cinf$
(resp. $\binftwo$ or $\dinf$) were established in \cite{W2}.
Lie groups and Lie superalgebra appearing in these duality
results are $ Sp (2l), Pin(2l),$ $ Spin(2l +1),$ 
$O(2l)$, $O(2l +1)$ and $ {\frak {osp}}(1, 2l)$ respectively.
Combining with our homomorphisms $\hphio$ or $\hphi$
(resp. $\hphic$),
we obtain free field realizations of $\hB$ (resp. $\hBp$)
and corresponding dual pairs involving $\hB$ (resp. $\hBp$).
We present explicit descriptions of the Fock space decompositions into
isotypic subspaces with respect to the joint action of
the corresponding dual pairs and 
calculate the corresponding {\em exponents} and their {\em multiplicities}
(which characterize the highest weights for $\hBpm$) in each isotypic subspace.
The explicit formulas for the highest weight vectors in isotypic subspaces
were already given in \cite{W2}.
Free field realizations of $\hD$ were studied earlier
in \cite{FKRW, KR2} and a duality between the general linear
Lie group and $\hD$ was established in these papers. 

It is a remarkable new feature of
dual pairs in our infinite dimensional setting that $\binf$ and $\cinf$ 
(resp. $\binftwo$ and $\dinf$)
involved in different dual pairs are uniformly replaced now by a single 
Lie algebra $\hB$ (resp. $\hBp$).
The ranks of these finite dimensional Lie groups, 
while varying for different Fock spaces,
correspond neatly to the central charges of the primitive $\hBpm$-modules
appearing in these Fock spaces. 
Note that \cite{W2} the Lie algebras of these Lie groups
turn out to be the horizontal subalgebras of the twisted or
untwisted affine algebras acting on the corresponding Fock spaces 
with level $\pm 1$ \cite{F1} \cite{KP} \cite{FF}.

We further show that the vacuum module $M_c$ of $\hBp$
with central charge $c$ and its irreducible
quotient $V_c$ admit a natural structure of a vertex algebra
\cite{B} \cite{FLM} \cite{DL} \cite{K2}. The module
$M_c$ is irreducible if and only if $c \notin \hf \Z$.
For the central charge $c = l, l +1/2, -l$ or $-l +1/2$, 
$V_c$ is shown to be isomorphic to the vertex algebra of
invariants of the various Fock spaces
$\Fl, \Flhalf, \Fminusl$ or $\Fhalfminusl$ (see the text for notations)
with respect to the action of $O(2l)$, $O(2l +1)$, $Sp (2l)$ or
${\frak {osp}} (1,2l)$ respectively. Those quasifinite representations
of $\hBp$ appearing from the Fock space decompositions
are representations of the vertex algebra $V_c$ for the corresponding
$c \in \hf \Z$.

It is well known that
to any complex simple Lie algebra $\frak g$ and $c \in \C$
one canonically associates a vertex algebra with central charge $c$,
called $\cal W$-algebra and denoted by ${\cal W}{\frak g}$, cf. \cite{BS, FeF}
and references therein. It is known \cite{BS, F2}
that the associated $\cal W$-algebra ${\cal WD}_l$
with central charge $l$ is isomorphic to the vertex algebra
of the $SO(2l)$-invariants in the basic representation of the affine
algebra $\widehat{ {\frak so}} (2l)$. We prove that ${\cal WD}_l$
with central charge $l$ is a sum of the vertex algebra $V_l$ and
an irreducible representation of $V_l$. This provides us a new way of
computing the $q$-character formula of ${\cal WD}_l$.

In a similar fashion, we identify the vertex algebra
of $SO(2l +1)$-invariants in the Fock space $\Flhalf$
as a $\cal W$-superalgebra ${\cal WB}(0, l)$ \cite{Ito} with 
$l$ bosonic generating fields of conformal weight $2, 4, \ldots, 2l$
and a fermionic generating field of conformal weight $l +1/2$.
The even part of ${\cal WB}(0, l)$ is isomorphic to $V_{l + 1/2}$
while the odd part of ${\cal WB}(0, l)$ is 
an irreducible representation of $V_{l + 1/2}$. 

The approach of this paper can be 
modified to study the representation theory of another
interesting Lie subalgebra of $\hD$ considered by Bloch \cite{Bl}.
This will be treated in a separate publication \cite{W3}.

The paper is organized as follows. In Section~\ref{sec_classical}
we describe the infinite rank Lie algebras
$\hglm, \binfm, \binfmtwo, \cinfm, \dinfm$ and present
the $q$-character formulas
of unitary highest weight modules over them in the case $m =0$
(see also Appendix).
In Section~\ref{sec_walgebra} we review the Lie algebra
$\hD$ and classify the anti-involutions of
$\hD$ preserving the principal $\Z$-gradation of $\hD$. 
In Section~\ref{sec_parab} we study the structure 
of parabolic subalgebras of $\hBpm$. 
In Section \ref{sec_finiteness} we give 
a characterization for a $\hBpm$-module to be quasifinite.
In Section \ref{sec_embed} we study the connection
between $\hBpm$ and the infinite rank Lie algebras. 
In Section~\ref{sec_realization} we realize quasifinite  $\hBpm$-modules 
in terms of modules over these infinite rank Lie algebras.
In Section~\ref{sec_unit} we classify 
all unitary quasifinite $\hBpm$-modules. 
In Sections~\ref{sec_boson}, \ref{sec_fermion} and \ref{sec_fermionhalf},
we realize a class of primitive $\hB$-modules with a half-integral
central charge in some Fock spaces and 
establish some duality between this class of primitive
$\hB$-modules with central charge $c= -l$ (resp. $-l -1/2$, $l$, $l + 1/2$)
and finite dimensional irreducible modules of $O (2l)$ (resp. $O(2l +1)$,
$Sp (2l)$ and $Pin (2l)$, $ {\frak {osp}}(1, 2l)$, $Spin(2l +1)$).
In Sections \ref{sec_n}, \ref{sec_nhalf} and \ref{sec_minusn}, 
we realize in a similar way a class of primitive $\hBp$-modules and establish 
various duality results for $\hBp$. 
In Section \ref{sec_voa}, we study the
representations of $\hBpm$ from the viewpoint of vertex algebras.

{\bf Acknowledgement}
This paper is based on two preprints \cite{KWY}.
W.~W. wishes to thank Max-Planck-Institut f\"ur Mathematik for its hospitality.
\section{Lie algebra $\hglm$ and its classical subalgebras}
   \label{sec_classical}
\subsection{Lie algebra $\hglm$}
   \label{subsec_glinf}
  Denote by $R_m$ the quotient algebra $\C [u] / (u^{m+1} )$
of the polynomial algebra $ \C [u]$ by the ideal generated by
$u^{m+1}$ $(m \in \Z_{+})$. Denote by ${\bf 1}$ the identity 
element of $R_m$. Denote by $\gl_f^{[m]}$ 
the complex Lie algebra of all infinite matrices
$(a_{ij})_{i,j \in \Z}$ with finitely many non-zero entries in $R_m$.
Denote by $\gl^{[m]}$ the Lie algebra of all matrices
$(a_{ij})_{i,j \in \Z}$ with only finitely many nonzero
diagonals with entries in $R_m$. 
Denote by $E_{ij}$ the infinite matrix with $1$ at $(i, j)$
place and $0$ elsewhere. 
Obviously $\gl_f^{[m]}$ is a Lie subalgebra of $\gl^{[m]}$.
There is a natural automorphism $\nu$ of $\hglm$ given by
\begin{equation}
   \nu ( E_{i,j}) = E_{i+1, j+1}.
 \label{eq_autom}
\end{equation}
Let the weight of $ E_{ij}$ be $ j - i $. This defines
the {\em principal} $\Z$-gradation 
$\gl^{[m]} = \bigoplus_{j \in \Z} \gl^{[m]}_j$.
Denote by $\hglm = \gl^{[m]} \bigoplus R_m$ the central extension
of $\hglm$ given by the following $2$--cocycle with values in $R_m$
\cite{KP, DJKM}:
\begin{eqnarray}
 C(A, B) = \Tr \,\left(
                     [J, A]B \right)
  \label{eq_cocy}
\end{eqnarray}
where $J = \sum_{ j \leq 0} E_{ii}$. The $\Z$--gradation
of Lie algebra $\glm$ extends to $\hglm$ by putting the weight $R_m$ to
be $ 0$. In particular, we have the {\em triangular decomposition}
$$\hglm = \widehat{\gl}^{[m]}_{+}  
          \bigoplus \widehat{\gl}^{[m]}_{0}
          \bigoplus \widehat{\gl}^{[m]}_{-} , $$
where 
$$ \glpm = \bigoplus_{ j \in \Bbb N}
\widehat{\gl}^{[m]}_{\pm j}, \quad 
\widehat{\gl}^{[m]}_{0} = \gl^{[m]}_0 \bigoplus R_m.$$
Given $\Lambda \in \left( {\hgl}^{[m]}_{0} \right)^{*}$, we let
\begin{eqnarray*}
 c_j & = & \Lambda ( u^j ),   \\
 {}^a \lambda_i^{(j)} & = & \Lambda ( u^j E_{ii}) ,  \\
 {}^a H^{(j)}_i     & = & u^j E_{ii} - u^j E_{i+1, i+1} + \delta_{i,0} u^j ,\\
  {}^a h^{(j)}_i & = & \Lambda( {}^a H^{(j)}_i ) = 
   {}^a \lambda_i^{(j)} - {}^a \lambda_{i+1}^{(j)} + \delta_{i,0} \; c_j ,
 \label{eq_11}
\end{eqnarray*}
where $i \in \Z, \;j = 0, \ldots, m.$
The superscript $a$ here denotes $\hglm$ which is of A type.
Denote by $L(\hglm; \Lambda)$ 
the highest weight $\hglm$--module with
highest weight $\Lambda$. The $c_j$ are called the 
{\em central charges} and 
${}^a \lambda_i^{(j)}$ are called the {\em labels} of $L(\hglm; \Lambda)$.

In particular putting $m =0$ we recover the well-known Lie algebras
$\gl = \gl^{[0]}$, $\hgl = \hgl^{[0]}$.
In this case, we drop the superscript $[0]$.
Define ${}^a \Lambda_j \in \gl_0^* $ $(j \in \Z)$ as follows:
\begin{equation}
{}^a \Lambda_j ( E_{ii} ) =    
   \left\{
      \everymath{\displaystyle}
      \begin{array}{lll}
        1, & \mbox{for}\quad 0 < i \leq j \\
        -1, & \mbox{for}\quad j < i \leq 0  \\
        0, & \mbox{otherwise.}
      \end{array}
    \right. \\
  \label{eq_122}
\end{equation}
Define ${}^a \hL_0 \in \hgl_0^{*}$ by 
$$ {}^a \hL_0 (C) = 1, \quad
  {}^a \hL_0 (E_{ii}) = 0 \mbox{ for all } i \in \Z $$
and extend $\Lambda_j\/$ from $ \gl_0^{*}\/$ to $\hgl_0^{*}$ by letting
$\Lambda_j (C) = 0.$ Then 
$${}^a \hL_j = {}^a \Lambda_j + {}^a \hL_0 \quad (j \in \Z ) $$ 
are the {\em fundamental weights},
i.e. ${}^a \hL_j ( {}^a H_i ) = \delta_{ij}.$

Recall that the q-character (i.e. principally specialized character) 
formula for an integrable 
highest weight module $L(\lambda)$ of a Kac-Moody algebra $\frak g$ 
(\cite{K}, chap. 10) with respect to the
principal gradation of $\frak g$ is given by
\begin{eqnarray}
 \ch_q L( \lambda) = \prod_{\alpha^{\vee} \in \Delta_+^\vee}
        \left( \frac{ 1 - q^{ \langle \lambda + \rho, \alpha^{\vee} \rangle} }
                    { 1 - q^{ \langle \rho, \alpha^{\vee} \rangle} }
        \right)^{mult \; \alpha^{\vee}}
  \label{eq_qchar}
\end{eqnarray}
where $\Delta_+^\vee$ is the positive coroots of the Lie algebra $\frak g$, 
and $\rho $ satisfies 
$\langle \rho, \alpha_i^{\vee} \rangle =1$ 
for all simple coroots $\alpha_i^{\vee}$.

Since $\hgl$ is a completed infinite rank affine algebra of Kac-Moody
type, the above formula (\ref{eq_qchar}) applies. This 
gives the following explicit formula 
(cf. \cite{FKRW})
\begin{eqnarray}
 \ch_q { L(\hgl, \Lambda) }
    = \frac{ \prod\limits_{ 1\leq i <j \leq c}
                                (1 - q^{ n_i -n_j +j -i} )}
                            {\varphi (q)^c}
  \label{eq_Aqchar}
\end{eqnarray}
where $\Lambda = \hL_{n_1}+\hL_{n_2}+ \ldots +\hL_{n_c},$
$ n_1 \geq n_2 \geq \ldots \geq n_c, c \in \Bbb N,$ 
and $\varphi (q)=\prod_{j=1}^{\infty} (1 - q^j)$ is the Euler product.
\subsection{Lie algebra $\binfm$}
     \label{subsec_binf}
   Now consider the vector space $R_m [t, t^{-1}]$ and take its basis
$v_i = t^i \; ( i \in \Z) $ over $R_m $. 
The Lie algebra $\glm$ acts on this vector
space via the usual formula
\begin{eqnarray}
 E_{ij} v_k = \delta_{j,k}v_i. 
  \label{eq_natur}
\end{eqnarray}
Let us consider the following $\C$-bilinear forms on this space:
\begin{eqnarray*}
 B^{\pm} (u^m v_i, u^n v_j ) & = & u^m (-u)^n (\pm 1)^i \delta_{i, -j} .
\end{eqnarray*}
Denote by $\Binfm$ (resp. $\Binfmtwo$)
the Lie subalgebra of $\glm$ which preserves the
bilinear form $B^-$ (resp. $B^+$), namely we have 
\begin{eqnarray*}
  \Binfm  & = & \{ g \in \glm \mid B( a(u), v) + B (u, a(v)) = 0 \}  \\
    & = & \{ (a_{ij} (u))_{i,j \in \Z} \in \glm \mid
           a_{ij} (u) = ( -1)^{i +j +1} a_{-j,-i} ( -u) \}, \\
  \Binfmtwo & = & \{ g \in \glm \mid B( a(u), v) + B (u, a(v)) = 0 \}  \\
    & = & \{ (a_{ij} (u))_{i,j \in \Z} \in \glm \mid
           a_{ij} (u) = - a_{-j,-i} ( -u) \}.
\end{eqnarray*}
Denote by $\binfm = \Binfm \bigoplus R_m$ 
(resp. $\binfmtwo = \Binfmtwo \bigoplus R_m$)
the central extension of $\Binfm$
(resp. $\Binfmtwo$) given by the $2$-cocycle
(\ref{eq_cocy}) restricted to $\Binfm$ (resp. $\Binfmtwo$). Then $\binfm$
(resp. $\binfmtwo$) inherits from $\hglm$ the principal $\Z$-gradation and 
the triangular decomposition:
\begin{eqnarray*}
 \binfm = \bigoplus_{ j \in \Z} b_{\infty, j}^{[m]}, &
 \binfm = b_{\infty, +}^{[m]} \bigoplus b_{\infty, 0}^{[m]}
 \bigoplus b_{\infty, -}^{[m]}, \\
 \binfmtwo  = \bigoplus_{ j \in \Z} \tilde{b}_{\infty, j}^{[m]}, &
 \binfmtwo  = 
    \tilde{b}_{\infty, +}^{[m]} \bigoplus
    \tilde{b}_{\infty, 0}^{[m]} \bigoplus \tilde{b}_{\infty, -}^{[m]} 
\end{eqnarray*}
where 
\begin{eqnarray*}
b_{\infty, j}^{[m]} = \binfm \cap \hglm_j,  & 
\binfmpm = \binfm \cap \glm_{\pm},          & 
b_{\infty, 0}^{[m]} = \binfm \cap \glm_0,    \\
\tilde{b}_{\infty, j}^{[m]} = \binfmtwo \cap \hglm_j,  & 
\tilde{b}_{\infty, \pm}^{[m]} = \binfmtwo \cap \glm_{\pm}  & 
\tilde{b}_{\infty 0}^{[m]}  = \binfmtwo \cap \glm_0.
\end{eqnarray*}

\begin{remark} \rm
 The Lie algebra $\binfmtwo$ is isomorphic to $\binfm$ by
 sending the elements $u^k E_{ij} - ( -u)^k E_{-j, -i}$ to
 $u^k E_{ij} + (-1)^{i +j +1} ( -u)^k E_{-j, -i}$, $i,j \in \Z, k \in \Z_{+}$.
 Their Cartan subalgebras coincide. 
\end{remark}

Given $\Lambda \in \left( {\binfm}_0 \right)^* $, we let
\begin{eqnarray}
   c_j & = &  \Lambda (u^j),  \nonumber   \\
 {}^b\lambda^{(j)}_0 & = & \Lambda \left( 2 u^j E_{00}  \right)
                       \quad (j \mbox{ odd}), \nonumber   \\
 {}^b\lambda^{(j)}_i & = & \Lambda \left( u^j E_{ii} - (-u)^j E_{-i, -i}
                                   \right) ,
                                              \nonumber   \\
  {}^b H^{(j)}_i & = &
   u^j E_{ii} - u^j E_{i+1,i+1} + (-u)^j E_{-i-1, -i-1} - (-u)^j E_{-i, -i},
                                              \nonumber  \\
   {}^b H^{(j)}_0 & = & 2 \left(u^j E_{-1,-1} -u^j E_{1,1} \right) + 2 u^j
                                             \quad (j \mbox{ even})  ,     
                                              \nonumber   \\
   {}^b H^{(j)}_0 & = & \left(2 u^j E_{0,0} - u^j E_{-1,-1} - u^j E_{1,1}
                        \right) +  u^j
                                             \quad (j \mbox{ odd}) ,    
                                              \nonumber     \\
  {}^b h^{(j)}_i & = & \Lambda( {}^b H^{(j)}_i )
     = {}^b \lambda^{(j)}_i - {}^b \lambda^{(j)}_{i+1} ,  \label{eq_bweight} \\
  {}^b h^{(j)}_0 & = & \Lambda( {}^b H^{(j)}_0 )
                        = - 2 \; {}^b \lambda^{(j)}_1 + 2 c_j
                                             \quad (j \mbox{ even}) ,     
                                              \nonumber   \\
  {}^b h^{(j)}_0 & = & \Lambda( {}^b H^{(j)}_0 )
           = {}^b \lambda^{(j)}_0 - {}^b \lambda^{(j)}_1 +  c_j
                   \quad (j \mbox{ odd}),
                                              \nonumber  
 \end{eqnarray}
where $ i \in \Bbb N$ and $j = 0, \ldots, m. $
The superscript $b$ here means B type.
Denote by $L(\binfm; \Lambda)$ (resp. $L(\binfmtwo; \Lambda)$)
the highest weight module over $\binfm$ (resp. $\binfmtwo$) with
highest weight $\Lambda$. The $c_j$ are called the 
{\em central charges} and 
${}^b \lambda_i^{(j)}$ are called the {\em labels} of $L(\binfm; \Lambda)$
or $L(\binfmtwo; \Lambda)$.
               
In particular when $m =0$ we have the usual Lie subalgebras
of $\gl$ and $ \hgl$, denoted by $\Binf$ and $\binf$ \cite{K}
(resp. $\Binftwo$ and $\binftwo$ \cite{W2}). 
Denote by ${}^b \hL_i $ the $i$-th {\em fundamental weight}
of $\binf$ (and $\binftwo$), namely
${}^b \hL_i ({}^b h_j ) = \delta_{ij}.$

The set of simple coroots of $\binf$ (and $\binftwo$), denoted by $ \Pi^\vee $,
can be described as follows:
\begin{eqnarray*}
      \begin{array}{rcl}
      \lefteqn{
          \Pi^\vee = \{ \alpha_0^\vee = 2( E_{-1, -1} - E_{1,1}) + 2C,  
              }             \\
       & & \alpha_i^\vee = E_{i, i} + E_{-i-1,-i-1}
              - E_{i+1,i+1} - E_{-i,-i}, i \in \Bbb N \}.
      \end{array}
\end{eqnarray*}
The set of roots is: 
$$\Delta = 
 \{\pm \epsilon_0, \pm \epsilon_i, \pm \epsilon_i \pm \epsilon_j, 
i \neq j, i,j \in \Bbb N \}.
$$
The set of positive coroots is: 
\begin{eqnarray*}
   \Delta_+^\vee & = & \{ \alpha_i^\vee + \alpha_{i+1}^\vee
          + \dots + \alpha_j^\vee, \quad 1 \leq i \leq j\}    \\
     & \quad & \cup \{ \alpha_0^\vee + 2\alpha_1^\vee 
      + \dots + 2\alpha_i^\vee +\alpha_{i+1}^\vee + \dots + \alpha_{j-1}^\vee, 
       \quad 0 \leq i < j \}.
\end{eqnarray*}
The set of simple roots is: 
$$\Pi = \{ \alpha_0 = - \epsilon_1, 
 \alpha_i = \epsilon_i - \epsilon_{i+1}, i \in \Bbb N \}.$$
Here $\epsilon_i $ are viewed restricted to the restricted dual of
the Cartan subalgebra of $\binf$, so that
$ \epsilon_i = - \epsilon_{-i}.$

Given $\Lambda = {}^b \hL_{n_1} + {}^b \hL_{n_2}
                   + \ldots + {}^b \hL_{n_k} + {}^b h\, {}^b \hL_0,
\quad  n_1 \geq n_2 \geq \ldots \geq n_k \geq 1, {}^b h \in \Z_{+},$
the $\binf$-module $L( \binf; \Lambda)$ has
central charge $ c = k + {}^b h/2$. Here and further we let
${}^b h = {}^b h^{(0)}$ and ${}^b \lambda_i = {}^b \lambda_i^{(0)}$,
see (\ref{eq_bweight}). Let
$$ \varphi_i (q) = \prod_{j=1}^i (1 - q^j).
$$
Further on we shall use the following notations:
$\overline{m} = 1$ for $m \in 2\Z$
and $\overline{m} = 0$ for $m \in 2 \Z +1$; for a real number $x$ denote by
$ [x]$ the integer no greater than and closest to $x$.
\begin{proposition}
  The $q$-character formula of $L( \binf; \Lambda)$ 
 (the same for $L( \binftwo; \Lambda)$) corresponding
 to the principal gradation of $\binf$ is
\begin{eqnarray}
   \ch_q L(\binf; \Lambda)
  & = & \frac{ \prod\limits_{1 \leq i < j \leq k} (1 - q^{n_i - n_j +j -i})}
             { \prod\limits_{1 \leq i \leq k} \varphi_{n_i + k -i}(q)} \cdot
 \frac{\varphi(q^2)^{\overline{2c +1}}
                         \prod\limits_{j >0} \varphi_{2c - 2j}(q)} 
              {\varphi(q)^{ [ \frac{2c + 1}2] } } \cdot     \nonumber   \\
   & \quad & \prod_{i=0}^{n_1-1}
               \frac{ \varphi_{2c + i + n_1}(q)}
                    { \varphi_{2c + n_1 + i - {}^b \lambda_{i +1}} (q) } \cdot
             \prod_{0 \leq i <j \leq n_1} 
               \frac{ 1 - q^{2c +j +i - {}^b \lambda_{i +1} - {}^b \lambda_j} }
                    {1-q^{2c +j +i}}.
  \label{eq_Bqchar}
\end{eqnarray}
  \label{prop_bform}
\end{proposition}
\begin{demo}{Proof}
   Applying formula (\ref{eq_qchar}) to $L( \binf; \Lambda)$, we get
  \begin{eqnarray}
   \ch_q L(\binf; \Lambda) = \prod_{1 \leq i < j}
      \frac{1 - q^{{}^b \lambda_i - {}^b \lambda_j + j -i} }
           {1 - q^{j -i}} \cdot
    \prod_{0\leq i <j}
         \frac{1 -q^{2c - {}^b \lambda_{i +1} - {}^b \lambda_j +j +i} }
              {1-q^{i+j}}.
   \label{eq_Bch}
  \end{eqnarray}
 Denote by the first (resp. second) product on the right hand 
 side of (\ref{eq_Bch}) to be $A$ (resp. $B$).
 Consider the Young diagram $Y_\Lambda $ corresponding 
 to the partition $(n_1, n_2, \ldots, n_k)$, it is easy to see that  
 $({}^b \lambda_1, {}^b \lambda_2, \ldots, {}^b \lambda_{n_1})$ 
 is its conjugate (cf. \cite{M} for terminology).
 Also ${}^b \lambda_i = 0$ for $ i > n_1$. Thus
 \begin{eqnarray}
   A = \prod_{1 \leq i <j}
                       \frac{1 - q^{{}^b \lambda_i - {}^b \lambda_j +j -i} }
                            {1 - q^{j-i} }
     = \prod_{j \in Y_\Lambda}(1 - q^{h_j})^{-1}
     = \frac{ \prod\limits_{1 \leq i < j \leq k} (1 - q^{n_i - n_j +j -i})}
           { \prod\limits_{1 \leq i \leq k} \varphi_{n_i + k -i}(q)}
   \label{eq_a}
 \end{eqnarray}
 where $h_j$ are the hook numbers of the Young diagram $Y_\Lambda$ 
 (cf. \cite{M}). 

 On the other hand, $B$ can be computed as follows.
 \begin{eqnarray*}
  B & = & \prod_{0 \leq i<j}
               \frac{1 -q^{2c + j + i - {}^b \lambda_{i+1} - {}^b \lambda_j} }
                    {1-q^{i+j}}        \\
    & = & \prod_{0 \leq i<j}
               \frac{1 -q^{2c + j + i} }
                    {1-q^{i+j}}      \cdot  
               \frac{1 -q^{2c + j + i - {}^b \lambda_{i+1} - {}^b \lambda_j} }
                    {1 -q^{2c + j + i} } \\
    & = & \prod_{0 \leq i<j}
               \frac{1 -q^{2c + j + i} }
                    {1-q^{i+j}}        
          \prod_{0 \leq i< n_1 < j}
               \frac{1 -q^{2c + j + i - {}^b \lambda_{i+1} - {}^b \lambda_j} }
                    {1 -q^{2c + j + i} }  \cdot  \\
    && \prod_{0 \leq i<j \leq n_1}
               \frac{1 -q^{2c + j + i - {}^b \lambda_{i+1} - {}^b \lambda_j} }
                    {1 -q^{2c + j + i} }.
  \end{eqnarray*}

A little further manipulation shows that the first, second and third
products of the right hand side of the equation above give rise
to the second, third and fourth products of the right hand side of 
(\ref{eq_Bqchar}) respectively.
\end{demo}
\subsection{Lie algebra $\cinfm$}
   \label{subsec_cinf}
      As before we take a basis $v_i = t^i \; (i \in \Z)$
of $R_m [t, t^{-1}]$ over $R_m $. 
Consider on this space the following $\C$-bilinear form
$$
 C (u^m v_i, u^n v_j ) = u^m (-u)^n ( -1)^i \delta_{i, 1-j} .
$$
Denote by $\Cinfm$ the Lie subalgebra of $\glm$ which preserves this
bilinear form:
\begin{eqnarray*}
 \Cinfm & = & \{ g \in \glm \mid C( a(u), v) + C (u, a(v)) = 0 \}    \\
        & = & \{ (a_{ij} (u))_{i,j \in \Z} \in \glm \mid
              a_{ij} (u) = ( -1)^{i +j +1} a_{1-j,1-i} ( -u) \}.
\end{eqnarray*}
Denote by $\cinfm = \Cinfm \bigoplus R_m$ 
the central extension of $\Cinfm$ given by the $2$-cocycle
(\ref{eq_cocy}) restricted to $\Cinfm$. Then $\cinfm$
inherits from $\hglm$ the natural $\Z$-gradation and
the triangular decomposition:
\begin{eqnarray*}
 \cinfm = \bigoplus_{ j \in \Z} {\cinfm}_j , \quad
\cinfm = {\cinfm}_{+} \bigoplus {\cinfm}_0 \bigoplus  {\cinfm}_{-} ,
\end{eqnarray*}
where ${\cinfm}_j = \cinfm \cap \hglm_j$,
$\cinfmpm = \cinfm \cap \glm_{\pm} $ 
and $ {\cinfm}_0 = \cinfm \cap \glm_0 .$

Given $\Lambda \in \left( {\cinf}^{[m]}_0 \right)^* $, we let
\begin{eqnarray}
   c_j & = &  \Lambda (u^j) ,
                                              \nonumber   \\
 {}^c \lambda^{(j)}_i & = & \Lambda \left( u^j E_{ii} - (-u)^j E_{1-i, 1-i}
                               \right)   , 
                                              \nonumber    \\
   {}^c H^{(j)}_i & = &
      u^j E_{ii} - u^j E_{i+1, i+1} + (-u)^j E_{-i, -i} - (-u)^j E_{1-i, 1-i},
                                                      \nonumber         \\
   {}^c H^{(j)}_0 & = & \left( u^j E_{0,0} - u^j E_{1,1} \right) + u^j 
                 \quad (j \mbox{ even}) ,
                                            \label{eq_cweight}     \\
  {}^c h^{(j)}_i & = & \Lambda( {}^c H^{(j)}_i )
     = {}^c \lambda^{(j)}_i - {}^c \lambda^{(j)}_{i+1}   ,     
                                              \nonumber \\
  {}^c h^{(j)}_0 & = & \Lambda( {}^c H^{(j)}_0 ) = - {}^c \lambda^{(j)}_1 + c_j
           \quad (j \mbox{ even})  ,  \nonumber 
\end{eqnarray}
where $ i \in \Bbb N$ and $j = 0, \ldots, m. $
For later purpose, it is convenient to put 
${}^c h^{(j)}_0 = c_j$ ($j$ odd), $j = 0, \ldots, m.$
The superscript $c$ here denotes $\cinf$ which is of C type.
Denote by $L(\cinfm; \Lambda)$ the highest weight $\cinfm$--module with
highest weight $\Lambda$. The $c_j$ are called the 
{\em central charges} and 
${}^c \lambda_i^{(j)}$ are called the {\em labels} of $L(\cinfm; \Lambda)$.
               
In particular, when $m =0$ we have the usual
$\Cinf = \Cinf^{[0]}$, $\cinf = \cinf^{[0]}$ \cite{K}. In this case,
we drop the superscript $[0]$.
Then we denote by ${}^c \hL_i $ the $i$-th fundamental weight
of $\cinf$, namely
${}^c \hL_i ({}^c h_j ) = \delta_{ij}.$

The set of simple coroots of $\cinf$, denoted by $ \Pi^\vee $,
can be described as follows:
\begin{eqnarray*}
      \begin{array}{rcl}
      \lefteqn{
          \Pi^\vee = \{ \alpha_0^\vee = (E_{0,0} - E_{1,1}) + C,  
              }             \\
       & & \alpha_i^\vee = E_{i, i} + E_{-i,-i}
              - E_{i+1,i+1} - E_{1-i,1-i}, i \in \Bbb N \}.
      \end{array}
\end{eqnarray*}
The set of roots is: 
$$\Delta = 
 \{ \pm \epsilon_i  \pm  \epsilon_j, \pm 2 \epsilon_i,
i \neq j, i,j \in \Bbb N \}.
$$
The set of positive coroots is: 
\begin{eqnarray*}
   \Delta_+^\vee & = & \{ \alpha_i^\vee + \alpha_{i+1}^\vee
          + \dots + \alpha_j^\vee, \quad 0 \leq i \leq j \}    \\
     & \quad& \cup \quad \{2 \alpha_0^\vee 
      + \dots + 2\alpha_i^\vee +\alpha_{i+1}^\vee + \dots + \alpha_j^\vee, 
       \quad 0 \leq i< j \}.
\end{eqnarray*}
The set of simple roots is: 
$$\Pi = \{ \alpha_0 = - 2 \epsilon_1, 
 \alpha_i = \epsilon_i - \epsilon_{i+1}, i = 1,2, \dots \}.$$
Here $\epsilon_i $ are viewed restricted to ${\frak h}'$, so that
$ \epsilon_i = - \epsilon_{-i +1}.$
               %

Given $\Lambda = {}^c \hL_{n_1} + {}^c \hL_{n_2}
         + \ldots + {}^c \hL_{n_k} + {}^c h\, {}^c \hL_0,
\quad  n_1 \geq n_2 \geq \ldots \geq n_k \geq 1, {}^c h \in \Z_{+},$
the $\cinf$-module $L( \cinf; \Lambda)$ has
central charge $ c = k + {}^c h$. 
\begin{proposition}
  The $q$-character formula of $L( \cinf; \Lambda)$ corresponding
 to the principal gradation of $\cinf$ is
\begin{eqnarray}
   \ch_q L(\cinf; \Lambda)
  & = & \frac{ \prod\limits_{1 \leq i < j \leq k} (1 - q^{n_i - n_j +j -i})}
             { \prod\limits_{1 \leq i \leq k} \varphi_{n_i + k -i}(q)} \cdot
   \frac{ \prod\limits_{1 \leq j \leq c -1}\varphi_{2j}(q)}
                     {\varphi(q)^c} \cdot   \nonumber \\
  &  &    \prod_{i=0}^{n_1}
           \frac{\varphi_{2c + n_1 +i} (q)}
                {\varphi_{2c + n_1 +i -{}^c \lambda_i} (q)} \cdot
          \prod_{0 \leq i <j \leq n_1} 
           \frac{1 - q^{2c + i + j - {}^c \lambda_{i+1} - {}^c \lambda_j} }
                {1 - q^{2c + i + j} }.
   \label{eq_Cqchar}
\end{eqnarray}
\end{proposition}
\begin{demo}{Proof}
   Applying formula (\ref{eq_qchar}) to $L( \cinf; \Lambda)$, we get
  \begin{eqnarray}
   \ch_q L(\cinf; \Lambda)
   & = & \prod_{1 \leq i < j}
      \frac{1 - q^{{}^c \lambda_i - {}^c \lambda_j + j -i} }
           {1 - q^{j -i}} \cdot  \nonumber  \\
   & &   \prod_{0\leq i <j}
        \frac{1 - q^{2c - {}^c \lambda_{i +1} - {}^c \lambda_{j +1} +j +i +2} }
             {1 - q^{i +j +2}} \cdot
    \prod_{j >0} \frac{1 - q^{c  - {}^c \lambda_j +j} }{ 1 - q^j} \nonumber \\
   & = & \prod_{1 \leq i < j}
      \frac{1 - q^{{}^c \lambda_i - {}^c \lambda_j + j -i} }
           {1 - q^{j -i}} \cdot
    \prod_{0\leq i <j}
        \frac{1 - q^{2c - {}^c \lambda_i - {}^c \lambda_j +j +i } }
             {1 - q^{i +j}}
   \label{eq_Cch}
  \end{eqnarray}
 with the convention ${}^c \lambda_0 = 0$.
 The first product on the right hand 
 side of (\ref{eq_Cch}) is already given by formula (\ref{eq_a})
 with ${}^b \lambda_i$ substituted by ${}^c \lambda_i$.
 Denote by $C$ the second product on the right hand 
 side of (\ref{eq_Cch}); it can be calculated 
 in a similar way as calculating $B$ in the proof of Proposition
 \ref{prop_bform}:
 \begin{eqnarray*}
  C & = & \prod_{0 \leq i < j}
             \frac{1 -q^{2c - {}^c \lambda_i - {}^c \lambda_j +i +j} }
                  {1-q^{i+j}}    \\
    & = & \frac{ \prod\limits_{1 \leq j \leq c -1}\varphi_{2j}(q)}
               {\varphi(q)^c} \cdot
          \prod_{i=0}^{n_1}
           \frac{\varphi_{2c + n_1 +i} (q)}
                {\varphi_{2c + n_1 +i -{}^c \lambda_i} (q)} \cdot
          \prod_{0 \leq i <j \leq n_1} 
           \frac{1 - q^{2c + i + j - {}^c \lambda_{i+1} - {}^c \lambda_j} }
                {1 - q^{2c + i + j} }.
 \end{eqnarray*}
\end{demo}
\subsection{Lie algebra $\dinfm$}
   \label{subsec_dinf}
      As before we take a basis $v_i = t^i \; (i \in \Z)$ 
of $R_m [t, t^{-1}]$ over $R_m $. 
Consider the following $\C$-bilinear form
$$ D (u^m v_i, u^n v_j ) = u^m (-u)^n \delta_{i, 1-j}, i, j \in \Z . $$
Denote by $\Dinfm$ the Lie subalgebra of $\glm$ which preserves this
bilinear form:
\begin{eqnarray*}
 \Dinfm & = & \{ g \in \glm \mid D( a(u), v) + D (u, a(v)) = 0 \}    \\
        & = & \{ (a_{ij} (u))_{i,j \in \Z} \in \glm \mid
              a_{ij} (u) = ( -1)^{i +j +1} a_{1-j,1-i} ( -u) \}.
\end{eqnarray*}
Denote by $\dinfm = \Dinfm \bigoplus R_m$ 
the central extension of $\Dinfm$ given by the $2$-cocycle
(\ref{eq_cocy}) restricted to $\Dinfm$. Then $\dinfm$
inherits from $\hglm$ a natural $\Z$-gradation and
the triangular decomposition:
\begin{eqnarray*}
 \dinfm & = & \bigoplus_{ j \in \Z} {\dinfm}_j \\
\dinfm & = & {\dinfm}_{+} \bigoplus {\dinfm}_0 \bigoplus  {\dinfm}_{-} 
\end{eqnarray*}
where ${\dinfm}_j = \dinfm \cap \hglm_j$,
$\dinfmpm = \dinfm \cap \glm_{\pm} $ 
and $ {\dinfm}_0 = \dinfm \cap \glm_0 .$

Given $\Lambda \in \left( {\dinf}^{[m]}_0 \right)^* $, we let
\begin{eqnarray}
   c_j & = &  \Lambda (u^j) ,
                                              \nonumber   \\
 {}^d \lambda^{(j)}_i & = & \Lambda \left( u^j E_{ii} - (-u)^j E_{1-i, 1-i}
                               \right)    ,
                                              \nonumber    \\
   {}^d H^{(j)}_i & = &
      u^j E_{ii} - u^j E_{i+1, i+1} + (-u)^j E_{-i, -i} - (-u)^j E_{1-i, 1-i},
                                                      \nonumber         \\
   {}^d H^{(j)}_0 & = & \left( (-u)^j E_{0,0}  + (-u)^j E_{-1,-1} 
     - u^j E_{2,2} - u^j E_{1,1} \right) + 2 u^j  ,  \label{eq_dweight}     \\
  {}^d h^{(j)}_i & = & \Lambda( {}^d H^{(j)}_i )
     = {}^d \lambda^{(j)}_i - {}^d \lambda^{(j)}_{i+1} ,       
                                              \nonumber \\
  {}^d h^{(j)}_0 & = & \Lambda( {}^d H^{(j)}_0 )
     = - {}^d \lambda^{(j)}_1 - {}^d \lambda^{(j)}_2 + 2 c_j ,  \nonumber 
\end{eqnarray}
where $ i \in \Bbb N$ and $j = 0, \ldots, m. $
The superscript $d$ here denotes $\dinf$ which is of D type.
Denote by $L(\dinfm; \Lambda)$ the highest weight $\dinfm$--module with
highest weight $\Lambda$. The $c_j$ are called the central charges and 
${}^d \lambda_i^{(j)}$ are called the labels of $L(\dinfm; \Lambda)$.
               
In particular, when $m =0$ we have the usual
$\Dinf = \Dinf^{[0]}$, $\dinf = \dinf^{[0]}$, cf. \cite{K1}. In this case,
we drop the superscript $[0]$.
Denote by $\hL^d_i $ the $i$-th fundamental weight
of $\dinf$, namely
${}^d \hL_i ({}^d h_j ) = \delta_{ij}.$

The set of simple coroots of $\dinf$, denoted by $ \Pi^\vee $,
can be described as follows:
\begin{eqnarray*}
      \begin{array}{rcl}
      \lefteqn{
          \Pi^\vee = \{ \alpha_0^\vee
             = (E_{0,0} + E_{-1,-1} - E_{2,2} - E_{1,1}) + 2C,  
              }             \\
       & & \alpha_i^\vee = E_{i, i} + E_{-i,-i}
              - E_{i+1,i+1} - E_{1-i,1-i}, i \in \Bbb N \}.
      \end{array}
\end{eqnarray*}
The set of roots is: 
$$\Delta = 
 \{ \pm \epsilon_i  \pm  \epsilon_j, i \neq j, i,j \in \Bbb N \}.
$$
The set of positive coroots is: 
\begin{eqnarray*}
   \Delta_+^\vee & = & \{ \alpha_i^\vee + \alpha_{i+1}^\vee
          + \cdots + \alpha_j^\vee \quad (0 \leq i \leq j),  \quad   
          \alpha_0^\vee + \alpha_2^\vee
          + \cdots + \alpha_j^\vee \quad (j \geq 2)  \}        \\
     & \quad& \cup \quad \{ \alpha_0^\vee +  \alpha_1^\vee + 2 \alpha_2^\vee
      + \cdots + 2\alpha_i^\vee +\alpha_{i+1}^\vee + \cdots + \alpha_j^\vee, 
       \quad 1 < i< j \}.
\end{eqnarray*}
The set of simple roots is: 
$$ \Pi = \{ \alpha_0 = - \epsilon_1 - \epsilon_2, 
 \alpha_i = \epsilon_i - \epsilon_{i+1}, i \in \Bbb N \}.$$
Here $\epsilon_i $ are viewed restricted to ${\frak h}'$, so that
$ \epsilon_i = - \epsilon_{-i +1}.$ The root vectors 
corresponding to $ \pm \alpha_0$ are
$e_0 = E_{0,2} - E_{-1, 1}$ and $f_0 = E_{2, 0} - E_{1, -1}.$
Introduce $\rho \in \left( {\dinf}^{[m]}_0 \right)^* $ such
that $(\rho, \alpha_i^{\vee}) =1, i \in \Z_{+}.$
\begin{remark} \rm
 The $q$-character formula (corresponding to the induced
 $\Z$-gradation on $\dinf$ from $\hgl$) of a $\dinf$-module 
 is given in the Appendix. Note that such a $\Z$-gradation
 on $\dinf$ is of type $(2,1, 1, \ldots)$
 rather than the principal $\Z$-gradation of type $(1, 1, \ldots)$
 as for $\hgl$, $\binf$ and $\cinf$ since 
 $f_0 = E_{2, 0} - E_{1, -1}$ has degree $-2$ instead of $-1$.
\end{remark}

\begin{example} \rm
   Take a parabolic subalgebra of $\Dinf$ 
 $$ {\cal P}_0 = \{ (a_{ij}) \in \Dinf
     \mid a_{ij} = 0 \mbox{ if } i > 0 \geq j \}, 
 $$
 and let $\widehat{\cal P}_0 = {\cal P}_0 \bigoplus \C C.$ The 2-cocycle 
 (\ref{eq_cocy}) when restricted to ${\cal P}_0$ is trivial.
 Consider the so-called vacuum module $M_c ( \dinf)$, i.e. the
 induced $\dinf$-module from a 1-dimensional
 ${\cal P}_0$-module on which ${\cal P}_0 $ acts as zero and
 $C$ as $ c \in \C$. The irreducible quotient of $M_c ( \dinf)$
 is $L(\dinf; 2c\; {}^d \hL_0 )$. 
 It follows by a standard argument (cf. \cite{K1}, Chapter 10) that
 \begin{eqnarray*}
  L(\dinf; 2c \; {}^d \hL_0 )
   = M_c ( \dinf) / \langle f_0^{2c +1} \rangle \mbox{   if }
  c \in  \Z_{+}/2.
 \end{eqnarray*}
  \label{ex_para}
\end{example}
\section{Anti-involutions of $\hD$ preserving its principal gradation}
   \label{sec_walgebra}
Let $\D_{as}$ be the associative algebra of regular differential operators on
the circle. The elements 
\begin{eqnarray*}
  J^l_k = - t^{l+k} ( \partial_t )^l \quad
    ( l \in \Z_{+}, k \in \Z )
\end{eqnarray*} 
form its basis, where $ \partial_t$ denotes $\frac{d}{dt}$. 
Another basis of $\D$ is
\begin{eqnarray*}
  L^l_k = - t^{k} D^l \quad
    (  l \in \Z_{+}, k \in \Z)   \nonumber
\end{eqnarray*}
where $D = t \partial_t$. It is easy to see that $J^l_k = -t^k [D]_l$.
Here and further we use the notation 
\begin{eqnarray}
 [x]_l = x(x-1)\ldots (x-l+1). 
   \label{eq_symb}
\end{eqnarray}
Let $\D$ denote the Lie algebra obtained from $\D_{as}$ by
taking the usual bracket $[a, b] =ab - ba.$
Denote by $\hD$ the central extension of $ \D $
by a one-dimensional center with a generator $C$: $\hD = \D + \C C$.
The Lie algebra $\hD$
has the following commutation relations \cite{KR1}
\begin{eqnarray}
  \left[
     t^r f(D), t^s g(D)
  \right]
    & = & t^{r+s} 
    \left(
      f(D + s) g(D) - f(D) g(D+r) 
    \right) \nonumber \\ 
   & + & \Psi 
      \left(
        t^r f(D), t^s g(D)
      \right)
      C
  \label{eq_12}
\end{eqnarray}
where 
\begin{eqnarray}
  \Psi 
      \left(
        t^r f(D), t^s g(D)
      \right)
   = 
   \left\{
      \everymath{\displaystyle}
      \begin{array}{ll}
        \sum_{-r \leq j \leq -1} f(j) g(j+r),& r= -s \geq 0  \\
        0, & r + s \neq 0. 
      \end{array}
    \right. 
  \label{cocy}
\end{eqnarray}
The 2-cocycle $\Psi$ can be equivalently given in terms of
another formula \cite{KP}:

\begin{eqnarray}
  \Psi \left( f( t)(\partial_t)^m, g( t)(\partial_t)^n
       \right)
  =\frac{m!n!}{(m+n+1)!} \mbox{Res}_{t = 0} f^{(n+1)}( t)g^{(m)}( t) dt.
  \label{cocy_kp}
\end{eqnarray}

Let the weight of $J^l_k$ be $ k$ and the weight of $ C$ be $0$.
This defines the principal $\Z$-gradations of $\D_{as}, \D$ and $\hD$:
\begin{eqnarray*}
  \D = \bigoplus_{j \in \Z} {\cal D}_j, \quad
  \hD = \bigoplus_{j \in \Z} \hat{\cal D}_j
\end{eqnarray*}
and so we have the triangular decomposition
\begin{eqnarray*}
  \hD = \hD_{+} \bigoplus \hD_{0}  \bigoplus \hD_{-},
\end{eqnarray*}
where 
\begin{eqnarray*}
  \hD_{\pm} = \bigoplus_{j \in \pm \Bbb N} \hat{\cal D}_j,
  \quad
  \hD_{0} = {\D}_{0}  \bigoplus {\Bbb C} C.
\end{eqnarray*}

 An {\em anti-involution} 
$\sigma$ of $\D_{as}$ is an involutive anti-automorphism of $\D$, i.e.
$ \sigma^2 = I, \sigma (a X + bY) =a\sigma (X) +b \sigma (Y)$ and
$ \sigma (XY)=\sigma (Y)\sigma (X),$ 
where $a, b \in \Bbb C, X, Y \in \D.$

\begin{proposition}
  Any anti-involution $\sigma$ of $\D_{as}$ which preserves the
 principal $\Z$-gradation is one of the following:
 \begin{enumerate}
  \item[1)] $\sigma_{-,b}( t) = -t, \quad \sigma_{-,b}(D) = -D +b;$ 
  \item[2)] $\sigma_{+,b}( t)= t, \quad \sigma_{+,b}(D) = -D +b, 
     \quad b \in \C. $
 \end{enumerate}
\end{proposition}

\begin{demo}{Proof}
 Since $\sigma$ preserves the principal $\Z$-gradation,
 we can assume that $\sigma(t )= tf(D),$ and $\sigma(D)=g(D)$, where
 $f$ and $g$  are polynomials. Since $\sigma$ is an anti-involution
 we have
 \begin{eqnarray}
    g(g(D)) & = & D,    \label{eq_lin}\\
   f(g(D))\cdot ( t f(D) ) & = & t.  \label{eq_f}
 \end{eqnarray}

 Equation (\ref{eq_lin}) means that the polynomial $g(w)$ 
 satisfies $ g(g( w)) = w$, which can only be possible when 
 $g (w)$ is linear in $w$. Hence
 $ g(w) = w$ or $g(w) = -w + b$ for some $b \in \C$.

 It follows from equation (\ref{eq_f}) that
 $$ 
   t = f(g(D) ) \cdot (t f(D)) = t f(g(D+1) ) f(D)
 $$
 which implies $ f(w)=\pm 1.$
 By applying the anti-involution $\sigma$ to the equation
 $ [t, D]= -t $, we have 
 $ [\sigma(D), \sigma(t)] = -\sigma(t), $ which excludes the 
 possibility $ g(w) = w. $
 
 An anti-involution $\sigma$ of $\D_{as}$ is completely determined by 
 $\sigma(t)$ and $\sigma(D)$. 
 On the other hand, it is straightforward to check the $\sigma_{\pm, b}$
 listed in the proposition are indeed anti-involutions of $\D_{as}$.
 This completes the proof of the proposition.
\end{demo}

It follows immediately that
 $
  \sigma_{\pm, b}  (\partial_t) = \mp \left(\partial_t - (b+1) t^{-1} \right).
 $ 
Given $ s \in \C$, denote by $\Ts$ the automorphism of $\D_{as}$ which sends
$t$ to $t$ and $D$ to $D + s$. Equivalently $\Ts$ is
given by sending $a \in \D$ to $t^{-s}a t^s $, 
the conjugate of $a$ by $t^s$.
Clearly $\Ts$
preserves the principal $\Z$-gradation of $\D_{as}$.
We have 
\begin{equation}
  \sigma_{\pm, b} \cdot \Ts = \sigma_{\pm, b + s}, \quad
  \Theta_{-s} \cdot \sigma_{\pm, b} = \sigma_{\pm, b + s} . 
 \label{eq_autoinvol}
\end{equation}
Denote by ${\cal D}^{\pm, b}$ the fixed Lie subalgebra
of $\D$ by $- \sigma_{\pm, b}$, namely 
$$ {\cal D}^{\pm, b} = \{ a \in \D \mid \sigma_{\pm, b} (a) = -a \}. $$
It inherits a $\Z$-gradation from $\D$
since $\sigma_{\pm, b}$ preserves the principal $\Z$-gradation of $\D$:
$ {\cal D}^{\pm, b} = \bigoplus_{j \in \Bbb Z} {\cal D}^{\pm, b}_j,$
where
 $$ {\cal D}^{\pm, b}_j = \{\  t^j f(D) \mid  f(w) \in \Bbb C[w]
    \mbox{ and } \sigma_{\pm, b} (t^j f(D)) = - t^j f(D) \}.
 $$

Let us denote by $\cwo$ the set of all odd polynomials 
in $ \C [w]$, and by $\cwe$ the set of all even polynomials 
in $ \C [w]$. As before we let $ \overline{k} = 0$
if $k$ is an odd integer and $ \overline{k} = 1$ if $k$ is even.
The following lemma gives a complete description of 
$ {\cal D}^{\pm, b}_j$.
\begin{lemma}
    We have
 \begin{eqnarray*}
   {\cal D}^{-, b}_j & = & \left\{ t^j g \left( D + (j -b)/2 \right)
      \mid  g(w) \in \C [w]^{(\bar{j})}, \; j \in \Z \right \} ,  \\
   {\cal D}^{+, b}_j & = & \left\{ t^j g \left(D + (j -b)/2 \right)
      \mid  g(w) \in \cwo, \; j \in \Z \right \}.
 \end{eqnarray*}
  \label{lem_graded}
\end{lemma}
\begin{demo}{Proof}
    Given $t^j f(D) \in {\cal D}^{-, b}_j $, we have
 $ \sigma_{-, b} (t^j f(D)) = -t^j f(D)$, which means
 $ (-1)^j t^j f(-D - j + b) = -t^j f(D). $ Equivalently we have
 $ (-1)^j f(-w - j + b) = -f( w). $
 Letting $g(w ) = f( w - \frac{j -b}2 )$, we have
 $ g( -w) = (-1)^{j +1}  g( w). $
\end{demo}

The relation among $ {\cal D}^{\pm, b}$ for different $b \in \C$ is
given by the following lemma which follows from Lemma~\ref{lem_graded}.
\begin{lemma}
  The Lie algebras $ {\cal D}^{+, b}$ (resp. $ {\cal D}^{-, b}$) 
 for different $b \in \C$ 
  are all isomorphic. More precisely, we have
  $ \Theta_s ( {\cal D}^{\pm, b}) = {\cal D}^{\pm, b -2s}.$
    \label{lem_iso}
\end{lemma}
 
Due to Lemma \ref{lem_iso} we may choose among $ {\cal D}^{-, b}$ a
Lie algebra, say $ \B = {\cal D}^{-, 0}$.
We see from Lemma \ref{lem_graded} that
$ D^n \in \B$ for $ n\in 2 \Z +1 $,
and $ D^n \notin \B$ for $ n\in 2 \Z $. Let (see notation (\ref{eq_symb}))
\begin{equation}
  T^{n, s}_k = - t^k \left(
                    [D -s]_n + (-1)^{k+1} [ -D -k -s]_n
                \right) \quad ( k \in \Z, n \in  \Z_{+}, s \in \C).
  \label{eq_basis}
\end{equation}
Clearly $T^{n, s}_k \in \B$. It is also clear that 
$ T^{n, s}_k$ $ ( k \in \Z, n \in 2 \Z_{+} +1)$
form a $\C$-basis of $\B$ for a
fixed $s \in \C$. We
denote $ T^{n, 0}_k$ by $ T^n_k$. A straightforward computation shows 
\begin{eqnarray}                
  T^{n, s}_k = - \left( 
              t^{k +n +s } \partial_t^{n} t^{-s}
               + (-1)^{k+n+1} t^{1 -s} \partial_t^n t^{k +n +s -1}
            \right).
   \label{eq_formal}
\end{eqnarray}
Although $t^{m \pm s}$ $( m \in \Z, s \in \C)$ does not lie in $\D$
in general, the above expression (\ref{eq_formal}) does lie in $\B$.

We denote again by $\Psi$ the restriction of the
2-cocycle $\Psi$ to $\B$, namely
\begin{eqnarray*}
  \begin{array}{rcl}
    \lefteqn{
  \Psi 
      \left(
        t^r f( D + r/2), t^s g( D + s/2)
      \right)} \\
   & = & 
   \left\{
      \everymath{\displaystyle}
      \begin{array}{ll}
        \sum_{-r \leq j \leq -1} f(j + r/2 ) g(j + r/2),
           & r= -s \geq 0  \\
        0, & r + s \neq 0 
      \end{array}
    \right. \\
  \end{array}
\end{eqnarray*}
where $ f(w) \in \C [w]^{( \bar{r})}, g(w) \in \C [w]^{( \bar{s})} .$
Denote by $\hB$ the central extension of $ \B $
by $ \C C$ corresponding to the 2-cocycle $\Psi$. 
$\hB$ is a Lie subalgebra of $\hD$ by definition.

The most convenient choice for the other family of
Lie algebras is $b = -1$ ($\sigma_{+, -1} ( \partial_t) = - \partial_t$).
Denote $ \Bp = {\cal D}^{+, -1}$.
We shall see that there is a canonical structure of a 
vertex algebra on
the vacuum module of the central extension of $\Bp$. Let
\begin{equation}
W^{n,s}_k = - \hf t^k \left(
                    [D + s]_n - [ -D -k -1 +s]_n
                \right) \quad ( k \in \Z, n \in \Bbb N, s \in \C).
  \label{eq_basisw}
\end{equation}
An important linear basis of $\Bp$ is
$ W^{n,s}_k$ $( k \in \Z, n \in 2 \Z_{+} +1)$ for a fixed $s \in \C$.
In particular we denote $W^{n, 0}_k$ by $W^n_k$.
A straightforward computation shows that
\begin{eqnarray}                
 W^n_k = - \hf \left( 
              t^{k +n} \partial_t^n + (-1)^{n +1}  \partial_t^n t^{k +n}
            \right).
  \label{eq_newbasis}
\end{eqnarray}
Note that the $W^1_k \; (k \in \Z)$ span a Virasoro algebra, namely we have
$$ \left[ W^1_m, W^1_n 
  \right]
  = (m -n) W^1_{m +n} + \delta_{m,-n} \frac{m^3 -m}{12} C.
$$

By abuse of notation we again denote by $\Psi$ the restriction of the
2-cocycle $\Psi$ to $\Bp$:
\begin{eqnarray*}
  \begin{array}{rcl}
    \lefteqn{
  \Psi 
      \left(
         t^r f ( D + (r +1)/2 ),
         t^s g ( D + (s +1)/2 )
      \right)    } \\
   & = &
   \left\{
      \everymath{\displaystyle}
      \begin{array}{ll}
        \sum_{-r \leq j \leq -1} f(j + (r +1)/2 ) g(j + (r +1)/2),
           & r= -s \geq 0  \\
        0, & r + s \neq 0,
      \end{array}
    \right. \\
  \end{array}
\end{eqnarray*}
where $ f(w), g(w) \in \cwo.$
Denote by $\hBp$ the central extension of $ \Bp $
by $ \C C$ corresponding to the 2-cocycle $\Psi$. The Lie algebra
$\hBp$ is a subalgebra of $\hD$ by definition.
\section{Structure of parabolic subalgebras of $\hBpm$}
   \label{sec_parab}
\subsection{$\hB$ case}
  We define a {\em parabolic subalgebra} ${\Pa}$ of $\hB$ 
as a subalgebra of the form
$ {\Pa} =\oplus_{j\in \Z} {\Pa}_j, $
where ${\Pa}_j = \hB_j$ if $j \geq 0$, and ${\Pa}_j \neq 0$ for some $j <0$. 

For each positive integer $k$ we have  
${\Pa}_{-k} = \{ t^{-k} h(D- k/2)\mid h(w) \in I_{-k} \}$, where
$ I_{-k}$ is some subspace of $\C [w]^{( \bar{k} )}$. 
Given $p(w) \in \C [w]^{( \bar{k} )}$ and $f (w) \in \C [w]^{(1)},$
we have $ f(D) , t^{-k} p( D - k/2) \in \hB$. One calculates that
\begin{eqnarray}
  \left[ f(D), t^{-k} p( D - k/2)
  \right]
    & = &  t^{-k} ( f(D -k) -f(D) ) p(D- k/2)  \nonumber   \\
    & = &  t^{-k} g(D - k/2) p(D- k/2)
  \label{eq_par}
\end{eqnarray}
where $ g (w) = f(w - k/2) - f( w + k/2).$
As $f$ ranges over all odd polynomials, $g(w)$ ranges over all 
even polynomials. Thus equation (\ref{eq_par}) implies that 
if $p(w) \in I_{-k} $ then $\cwe p(w) \subset I_{-k}$. 
This means that $I_{-k}$ is a submodule of $\C [w]^{( \bar{k} )}$
over $\C [w]^{( 0 )}$, where $\C [w]^{( \bar{k} )}$
is regarded as a module over $\C [w]^{( 0 )}$ by multiplication. 
Clearly every non-zero submodule of $\C [w]^{( \bar{k} )}$ 
over $\C [w]^{( 0 )}$
is a free rank 1 submodule generated by a monic polynomial. Denote by $b_k (w)$
such a generator for $I_{-k}$ if $I_{-k} \neq 0$,
and let $b_k (w) = 0$ if $I_{-k} = 0$. We call $b_k(w) \; ( k = 1, 2, \ldots)$
the {\em characteristic polynomials} of ${\Pa}$. 
\begin{lemma}
  Let $\{b_k,k \in \Bbb N \}$ be the sequence of characteristic 
polynomials of a parabolic subalgebra ${\Pa}$ of the Lie algebra $\hB$. 
Then
 \begin{enumerate}
  \item[1)] $b_k(w ) \in \C [w]^{(\bar{k})}$;
  \item[2)] $b_k(w)$ divides $ w b_{k+1} (w - 1/2)$ and 
   $w b_{k+1} (w + 1/2)$ for all $k \in \Bbb N$;
  \item[3)] $b_{k+l} (w)$ divides $w b_k (w -l/2) b_l ( w + k/2)$ 
   for all $l, k \in \Bbb N$;
  \item[4)] ${\Pa}_{-k} \neq 0$ for all $k \in \Bbb N$.
 \end{enumerate}
 In particular, all $b_k (w)$ are non-zero.
 \label{lem_charact}
\end{lemma}
\begin{demo}{Proof}
  Part 1) follows from the definition of the characteristic polynomials. From 
 the commutation relation
 \begin{eqnarray*}
   \begin{array}{rcl}
    \lefteqn{
     \left[ t, t^{-k-1} b_{k+1} \left( D - (k +1)/2
                                \right) 
     \right]
            }   \\
   & & = t^{-k} \left( b_{k+1} \left( D - (k+1)/2
                         \right)
        - b_{k+ 1} \left( D - (k-1)/2
                   \right)
               \right)
   \end{array}
 \end{eqnarray*}
 we see that $ b_k ( w - k/2 ) $ divides 
  \begin{equation}
      b_{k+1} \left( w -  (k +1)/2
                             \right)
    - b_{k+1} \left( w- (k -1)/2
              \right).
    \label{eq_divid1}
 \end{equation}
 In particular $ b_{k +1} \neq 0$ implies $ b_k \neq 0$.
 From  the commutation relation
 \begin{eqnarray*}
   \begin{array}{rcl}
    \lefteqn{
  \left[ t, \; t^{-k-1} \left( D- (k+1)/2
                    \right)^2 b_{k+1} \left( D- (k+1)/2
                                      \right)
  \right]    
            }           \\
   & & = t^{-k} \left\{ \left\{ D- (k -1)/2
                   \right)^2 \left( b_{k+1} \left( D - (k+1)/2
                                            \right)
                                    - b_{k+1}\left( D - (k -1)/2
                                            \right)
                        \right\}
            \right.         \\
   & & \quad           -2 \left( D- k/2
                      \right) 
                      b_{k+1} \left( D - (k+1)/2
                              \right)
            \Big\}
   \end{array}
\end{eqnarray*}
we see that $b_k \left( w - k/2 \right)$ divides 
\begin{eqnarray}
            \left( w - (k -1)/2
                                       \right)^2 
            \left[ b_{k+1} \left( w - (k+1)/2
                           \right) -b_{k+1} \left(w - (k -1)/2
                                            \right)
           \right]            \nonumber   \\
           - 2 \left( w - k/2
               \right) b_{k+1} \left( w -(k +1)/2
                               \right).
   \label{eq_divid2}
\end{eqnarray}
 From (\ref{eq_divid1}) and (\ref{eq_divid2}), we see that
 $b_k (w) $ divides $ w b_{k +1} ( w - 1/2)$.
 Noting  that
 $b_k (w) = \pm b_k ( -w)$ and $b_{k+1} (w) = \pm b_{k+1} ( -w)$, 
 $ b_k (w)$ also divides $ w b_{k +1} (w + 1/2 ). $ This proves part 2).

 Part 3) can be similarly proved by calculating
 the following two commutators
 $$ \left[ t^{-k} b_k (D - k/2), t^{-l} b_l (D - l/2)
    \right] 
 $$
 and
 $$ \left[ t^{-k} ( D- k/2)^2 b_k (D - k/2), t^{-l} b_l (D- l/2)
    \right]. 
 $$
 In particular, $b_k \neq 0$ and $ b_l \neq 0$ imply $b_{k +l} \neq 0$.
 Part 4) follows from 2) and 3).
\end{demo}
   
Given a monic even polynomial $b = b(w)$, denote
by $ \hB_0 (b) $ the subspace of $\hB_0$ spanned by
 $$  f(D- 1/2 ) b( D - 1/2 ) - f (D + 1/2) b ( D + 1/2)
     + f ( - 1/2) b ( - 1/2) C
$$
where $ f(w) \in \cwe .$
We have the following proposition.
\begin{proposition}
   Let ${\Pa}$ be a parabolic subalgebra of $\hB$ and let $b = b_1 (w)$ 
 be its first characteristic polynomial. Then
 $$ 
  [{\Pa}, {\Pa}] ={\hB}_0 (b) \bigoplus \bigoplus_{k\neq 0} {{\Pa}}_k.
 $$
  \label{prop_commutator}
\end{proposition}
\begin{demo}{Proof}
    It follows from $\left[ D, t^k f(D) \right] = k t^k f(D)$
 that $ {{\Pa}}_k \subset [{\Pa}, {\Pa} ]$ for $ k \neq 0.$ Next we claim
 $\left[ \B_1, \B_k \right] = \B_{k +1}$ for $k \in \Bbb N$. Indeed, first
 consider the case when $k$ is an even integer.
 \begin{equation}
  \begin{array}{rcl}
   \lefteqn{
      \left[ t ( D + 1/2)^{2l}, t^k ( D + k/2)
    \right]
           }      \\
       & = & t^{k+1} \left( (D +k + 1/2 )^{2l} ( D + k/2) 
                        - ( D + 1/2)^{2l} ( D + k/2 +1)
                 \right).
   \end{array}
    \label{eq_span}
 \end{equation}
 Note that the leading term of 
 $ (D +k + 1/2 )^{2l} ( D + k/2) - ( D + 1/2)^{2l} ( D + k/2 +1)$ is
 $ (2kl -1) D^{2l} $. It is clear that when $l$ ranges over $\Bbb Z_+$, 
 the right hand side of (\ref{eq_span}) form  a basis of $\hB_{k+1}$. 
 The case when $k$ is odd can be treated similarly.

 Now from the fact that
 ${{\Pa}}_k = \B_k $ for $k \geq 0$, we have
 \begin{eqnarray*}
  \begin{array}{rcl}
   \lefteqn{
   [{\Pa}_{k +1}, {\Pa}_{-k -1} ]
   \subset \left[ [{{\Pa}}_1, {{\Pa}}_k], {{\Pa}}_{-1-k} \right]  
           }      \\
   & \subset & \left[ [{\Pa}_k, {\Pa}_{-k -1}], {\Pa}_1 \right]
         + \left[ [{\Pa}_{-k -1}, {\Pa}_1], {\Pa}_k \right]   
   \subset [{\Pa}_{-1}, {\Pa}_1] + [{\Pa}_{-k}, {\Pa}_k].
   \end{array}
 \end{eqnarray*}
 Hence by induction, $[{\Pa}, {\Pa}]_0=[{\Pa}_1, {\Pa}_{-1}]$.
 A direct computation shows that $[{\Pa}_1, {\Pa}_{-1}]$ 
 is exactly $\hB_0 (b)$.
\end{demo}
\subsection{$\hBp$ case}
 Given a parabolic subalgebra $ {\Pa} =\bigoplus_{j\in \Z} {\Pa}_j $
of $\hBp$, we have:
$$
 {\Pa}_{-k} = \{ t^{-k} h(D + (-k +1)/2)\mid h(w) \in I_{-k} \} \quad
 (k \in \Bbb N ), 
$$
where $ I_{-k}$ is some subspace of $\cwo$. 
Given $p(w), f (w) \in \cwo$ we have:
$ t^{-k} p( D + (-k +1)/2)$, $ f(D +1/2) \in \hBp$. One calculates that
\begin{eqnarray}
 \begin{array}{rcl}
   \lefteqn{
  \left[ f\left( D + 1/2 \right), t^{-k} p \left( D + (1 -k)/2 \right)
  \right]} \\
    & = &
   t^{-k} g \left( D + (1 -k)/2 \right)
    p \left( D + (1 -k)/2 \right) ,
 \end{array}
  \label{eq_pard}
\end{eqnarray}
where $ g (w) = f(w - k/2) - f( w + k/2).$
As $f$ ranges over all odd polynomials, $g(w)$ ranges over all 
even polynomials. Thus (\ref{eq_pard}) implies that 
if $p(w) \in I_{-k} $ then $p(w)$ multiplied by any
even polynomial belongs to $ I_{-k}$. 
Let $b_k (w) \; (k \in \Bbb N)$ be the unique monic odd polynomial
in $ I_{-k} $ of minimal degree when $ I_{-k} \neq 0$ and 
let $b_k (w) = 0$ when $ I_{-k} =0 $.
We call $b_k(w) \; ( k = 1, 2, \ldots)$
the {\em characteristic polynomials} of ${\Pa}$. 
\begin{lemma}
  Let $\{b_k, k \in \Bbb N \}$ be the sequence of characteristic 
polynomials of a parabolic subalgebra ${\Pa}$ of the Lie algebra $\hBp$. 
Then
 \begin{enumerate}
  \item[1)] $b_k(w)$ divides $ w ( w + (k +1)/2 ) b_{k+1} (w + 1/2)$ 
   for all $k \in \Bbb N$;
  \item[2)] $b_{k+l} (w)$ divides $w b_k (w +l/2) b_l ( w - k/2)$ 
   for all $l, k \in \Bbb N$;
  \item[3)] ${\cal P}_{-k} \neq 0$ for all $ k \in \Bbb N$.
 \end{enumerate}
 \label{lem_chara}
\end{lemma}
\begin{demo}{Proof}
 It follows from  
 \begin{eqnarray*}
 \begin{array}{rcl}
   \lefteqn{
 [t (D +1),t^{-k -1} b_{k +1} ( D- k/2) ] } \\
  & = & t^{-k} \left(
              (D -k)b_{k +1} (D - k/2) - (D +1)b_{k +1} (D - k/2 +1)
           \right)
 \end{array}
 \end{eqnarray*}
that $b_k (w +(-k +1)/2)$ divides 
 \begin{equation}
(w -k)b_{k +1} (w - k/2) - (w +1)b_{k +1} (w - k/2 +1).   \label{eq_one}
 \end{equation}
 We see that
 $b_k (w +(-k +1)/2)$ divides 
 \begin{equation}
 (w -k) (w -k/2)^2 b_{k +1} ( w - k/2)
   - (w +1)(w - k/2 +1)^2 b_{k +1} (w - k/2 +1)
   \label{eq_two}
 \end{equation}
 by computing 
 $[t (D +1),t^{-k -1} ( D- k/2)^2 b_{k +1} ( D- k/2) ].$
 Thus $b_k (w +(-k +1)/2)$ divides (\ref{eq_two}) subtracted by
 (\ref{eq_one}) multiplied with $(w - k/2)^2$
 which is equal to $ (w +1) (w + (-k +1)/2 )b_{k +1} (w - k/2 +1).$
 This proves 1). The above computation shows that
 $b_{k+1}(w) \neq 0$ implies $b_k (w) \neq 0$.
 
 Part 2) can be similarly proved by computing the following
 two commutators 
 $$
 [ t^{ -k} b_k (D + (-k +1)/2), t^{ -l} b_l (D + (-l +1)/2 )] ,  $$
 $$
[ t^{ -k}(D + (-k +1)/2)^2  b_k (D + (-k +1)/2), t^{ -l} b_l (D + (-l +1)/2) ].  $$
 Similarly it follows
 that  $b_k(w), b_l (w) \neq 0$ implies $b_{k +l} (w) \neq 0$.
 Now part 3) follows from 1) and 2).
\end{demo}
\begin{lemma}
  $ \left[ \hBp_1, \hBp_k \right] = \hBp_{ k +1} \;\; (k >1), \;
  \hBp_2 = \left[ \hBp_1, \hBp_1 \right] \bigoplus \C \;t^2 ( D+3/2).$
     \label{lem_add}
\end{lemma}
\begin{demo}{Proof}
   First we have 
 \begin{equation}
  \begin{array}{rcl}
   \lefteqn{
      \left[ t ( D + 1)^l, t^k ( D + (k +1)/2)^m
    \right]
           }      \\
       & = & t^{k+1} \left( (D +k + 1 )^l ( D + (k +1)/2)^m 
                        - ( D + 1)^l ( D + (k +3)/2 )^m
                 \right).
   \end{array}
    \label{eq_spand}
 \end{equation}
 For odd positive integers $l$ and $ m$, 
 $t ( D + 1)^l, t^k ( D + (k +1)/2)^m \in \hBp$.
 For $m =1$, the leading term of the 
 right hand side is $( lk -1) D^l$.
 When $k >1$, $( lk -1) \neq 0$ for $l \in \Bbb N$.
 Thus in the case $k >1$ the right hand side of (\ref{eq_spand})
 span the whole $\hBp_{k +1}$ when $l$ ranges over all odd positive integers. 
 In the case $k =1$, the right hand side of (\ref{eq_spand})
 when $l$ ranges over all odd positive integers
 together with $ t^2 (D +1)$ span the whole $\hBp_2.$

 On the other hand, putting $k =1$ in equation (\ref{eq_spand})
 we see that the right hand side always
 contains a factor $(D + 1)(D +2)$. 
 So $[ \hBp_1, \hBp_1 ] $ does not contain $ t^2 (D +3/2)$.
\end{demo}

Let $\hBp_0 (b_1, b_2)$ denote the subspace of $\hBp_0$ spanned by
\begin{eqnarray*}
  \{   g(D- 1/2 ) b_2( D - 1/2 ) - g (D + 3/2) b_2 ( D + 3/2) 
          + 2 g ( 1/2) b_2 (  1/2) C  ,&    \\
f(D ) b_1( D  ) - f (D + 1) b_1 ( D + 1),\;   f , g \in \cwo \} .    &    
\end{eqnarray*}

\begin{proposition}
   Let ${\Pa}$ be a parabolic subalgebra of $\hBp$ and let $b_i = b_i(w)$ 
 $ (i = 1, 2)$ be its first and second characteristic polynomials. Then
 $$ 
  [{\Pa}, {\Pa}] ={\hBp}_0 (b_1, b_2) \bigoplus \bigoplus_{k\neq 0} {{\Pa}}_k.
 $$
  \label{prop_commutat}
\end{proposition}
\begin{demo}{Proof}
      It follows from $\left[ D, z^k f(D) \right] = k z^k f(D)$
 and $D \in \hBp$
 that $ {{\Pa}}_k = [{\Pa}, {\Pa} ]_k$ for $ k \neq 0.$ Since
 ${{\Pa}}_k = \B_k $ for $k \geq 0$, it follows from
 Lemma \ref{lem_add} that for $k >1$
 \begin{eqnarray*}
  \begin{array}{rcl}
   \lefteqn{
   [{\Pa}_{k +1}, {\Pa}_{-k -1} ]
   \subset \left[ [{{\Pa}}_1, {{\Pa}}_k], {{\Pa}}_{-1-k} \right]  
           }      \\
   & \subset & \left[ [{\Pa}_k, {\Pa}_{-k -1}], {\Pa}_1 \right]
         + \left[ [{\Pa}_{-k -1}, {\Pa}_1], {\Pa}_k \right]   
   \subset [{\Pa}_{-1}, {\Pa}_1] + [{\Pa}_{-k}, {\Pa}_k].
   \end{array}
 \end{eqnarray*}
 Hence it follows by induction that
\begin{equation}
  [{\Pa}, {\Pa}]_0 = [{\Pa}_1, {\Pa}_{ -1}] + [{\Pa}_2, {\Pa}_{ -2}].
   \label{eq_cartan}
\end{equation}
 
 A direct computation shows the right
 hand side of (\ref{eq_cartan}) is indeed  $\hBp_0 (b_1, b_2)$.
\end{demo}
\begin{example} \rm
   Let
 ${\cal P} = \{ W^n_k \mid n +k \geq 0, k \in \Z, n \in 2 \Z_{+} +1 \}$
 and let $\widehat{\cal P} = {\cal P} \bigoplus \C C.$
 By using (\ref{eq_newbasis}) it is easy to see that
 ${\cal P}$ is closed under the Lie bracket and
 thus is a parabolic subalgebra of $\Bp$. Geometrically,
 ${\cal P}$ consists of those differential operators in $\hBp$
 which extend to the interior of the circle.
 It is clear from formula (\ref{cocy_kp}) that
 the 2-cocycle $\Psi$ when restricted to ${\cal P}$ is trivial.
 Denote by $M_c$ the generalized Verma module
 $M (\hBp, {\cal P}, \xi_0)$ where $ \xi_0 $ has labels
 $\Delta_i =0$ for all $i$ and central charge $c$. Denote by
 $V_c$ the irreducible quotient of the $\hBp$-module $M_c$.
  \label{ex_para2}
\end{example}
\section{Characterization of quasifiniteness of HWM's of $\hB$}
  \label{sec_finiteness}
Let $ \frak g =\oplus_{j \in \Bbb Z} \frak g_j$ 
(possibly $\dim {\frak g}_j = \infty$)
be a $\Bbb Z$-graded Lie algebra over $\Bbb C$ satisfying
$ [\g_i, \g_j] \subset \g_{i+j},$ and let 
${\frak g}_{+} =\oplus_{j >0} \frak g_j$.
A $\g$-module $V$ is called {\it $\Z$-graded} if 
$ V=\oplus_{j \in \Z} V_j,$ and $\g_i V_j \subset V_{j -i}. $
A graded $\frak g$-module is called {\it quasifinite} if 
$ \dim V_j < \infty$ for all $j$.

Given $\Lambda \in \frak g_0^*$, a {\em highest weight module} (HWM) is 
a $\Bbb Z$-graded $\frak g$-module 
$V(\frak g, \Lambda)=\oplus_{j \in \Bbb Z_+} V_j$ 
generated by a highest weight vector $v_{\Lambda} \in V_0$ which 
satisfies 
$$ h v_{\Lambda} = \Lambda (h) v_{\Lambda} \quad(h \in \frak g_0), \quad \quad
\frak g_+ v_\Lambda =0.
$$
A non-zero vector $v \in V(\frak g, \Lambda)$ is called {\it singular} 
if $\g_+ v =0$. 

A Verma module is defined as
$$
 M(\frak g; \Lambda) = {\cal U}(\frak g)
  \bigotimes_{ {\cal U}({\frak g}_0 \oplus {\frak g}_+)} \Bbb C_\Lambda
$$
where $\Bbb C_\Lambda$ is the 1-dimensional $( \g_0 \oplus\g_+ )$-module 
given by $h \mapsto \Lambda (h)$ if $h \in \frak g_0$ and $\g_+ \mapsto 0$.
Here and further $\cal U(\frak s)$ stands for the universal enveloping 
algebra of the Lie algebra $\frak s$.
Any highest weight module $V(\frak g, \Lambda)$ is a quotient module of 
$M(\frak g; \Lambda)$. The irreducible module $L(\frak g; \Lambda)$ 
is the quotient of $M(\frak g; \Lambda)$ 
by the maximal proper graded submodule.

Let ${\Pa} = \oplus_j {\Pa}_j$ be a parabolic subalgebra 
of $\frak g$, and let $\Lambda \in \frak g_0^*$ be such that 
$\Lambda|_{ \g_0 \cap [{\Pa}, {\Pa}] } = 0$. Then the 
$( \g_0 \oplus \g_+)$-module $\Bbb C_{\Lambda}$ extends to 
a ${\Pa}$-module by letting
${\Pa}_j $ act as $ 0$ for $j<0$, and we may construct 
the highest weight module
$$
  M(\frak g, {\Pa}, \Lambda) = {\cal U} (\frak g)
                 \bigotimes_{ {\cal U}( {\Pa})}\Bbb C_\Lambda
$$
which is usually referred to as the {\it generalized Verma module}. 
Clearly all the highest weight modules are graded.

In the following we consider $\g = \hBpm$ and $\xi \in (\hBpm)_0^*$.
Let $b(w)$ be a monic even polynomial
(resp. $b_1(w), b_2 (w)$ be two monic odd polynomials).
Let $\xi \in \left( \hB_0 \right)^{*}$ 
be such that $\xi \mid_{\hB_0 (b) } = 0$, 
(resp. $\xi \in \left( \hBp_0 \right)^{*}$ be such that 
$\xi \mid_{\hBp_0 (b_1, b_2) } = 0$).
Consider a parabolic subalgebra $\Pa$ of $\hB$ (resp. $\hBp$)
whose first characteristic polynomial is $b(w)$ (resp.
whose first and second characteristic polynomials are $b_1(w), b_2 (w)$).
Denote by $M(\hB; \xi, b)$ (resp. $M(\hBp; \xi, b_1, b_2)$) the generalized 
Verma module $M( \hB, \Pa, \xi)$ (resp. $M( \hBp, \Pa, \xi)$). 
The polynomial $b$ (resp. $b_1, b_2$) does not determine 
$\Pa$ uniquely, but for our need, any corresponding parabolic $\Pa$ will do.
\begin{proposition}
 The following conditions on 
 $\xi \in \left( \hBpm_0 \right)^{*}$ are equivalent:
  \begin{enumerate}
   \item[1)] $M(\hB; \xi, b)$ (resp. $M(\hBp; \xi, b_1, b_2)$)
  contains a singular vector in its first graded subspace;
   \item[2)] $L(\hBpm; \xi) $ is quasifinite;
   \item[3)] $L(\hBpm; \xi)$ is a quotient of a generalized Verma module 
     $M(\hB; \xi, b)$ (resp. $M(\hBp; \xi, b_1, b_2)$)
   for some monic even polynomial $b$ (resp. some monic odd polynomials
   $b_1, b_2$).
  \end{enumerate}
 \label{prop_equiv}
\end{proposition}
\begin{demo}{Proof} We give the proof for the $-$ case. The proof for
 the $+$ case is essentially the same.
   1) $\Rightarrow$ 3). Denote by $( t^{-1} b( D- 1/2) ) v_\xi$
 the singular vector where 
 $ b(w)$ is a monic even polynomial of minimal degree (note that
 $ t^{-1} b( D- 1/2) \in \hB_{-1}$).
 Then it is easy to see that 3) holds for this particular 
 monic even polynomial $b$. 3) $\Rightarrow$ 2) and
 2) $\Rightarrow$ 1) follow by Lemma \ref{lem_charact}. 
\end{demo}
       
Let $L(\xi)$ be an irreducible QHWM 
over $\hBpm$. According to Proposition \ref{prop_equiv}, we have
$( t^{-1} b( D- 1/2) ) v_\xi = 0$ for some monic even 
polynomial $b(w)$ in the $-$ case, and
$( t^{-1} b( D ) ) v_\xi = 0$ for some monic odd
polynomial $b(w)$ in the $+$ case.
Such a monic polynomial of minimal degree 
is uniquely determined by $\xi$ and
is called the {\em characteristic polynomial}
of $L(\xi)$. 

We shall characterize a weight $\xi \in ( \hB_0 )^{*}$ 
(resp. $\xi \in ( \hBp_0 )^{*}$) by its labels $\Delta^-_n = -\xi( D^n)$ 
(resp. $\Delta^+_n = -\xi( (D +1/2)^n)$), where 
$n \in \Bbb N_{{odd}} = \{ 1, 3, 5, \ldots \}$,
and the central charge $c=\xi (C)$. Introduce the generating series
\begin{eqnarray*}
 \Delta^{\pm}_{\xi}(x) = \sum_{n \in {\Bbb N}_{{ odd} } } 
\frac{x^n}{n!} \Delta^{\pm}_n.
\end{eqnarray*}
Sometimes we simply write $ \Delta^{\pm}(x)$ instead of 
$ \Delta^{\pm}_{\xi}(x)$, or even drop $\pm$ when no confusion
may arise. Clearly we have
\begin{eqnarray}
 \Delta^- (x) & = & - \hf \xi \left( e^{xD} - e^{ -xD} \right), \\
 \Delta^+ (x) & = & - \hf \xi \left( e^{x(D + 1/2)} - e^{ -x(D + 1/2)} \right).
   \label{eq_triv}
\end{eqnarray}

A {\it quasipolynomial} is a finite
linear combination of functions of the form 
$p(x) e^{ \alpha x}$, where $p(x)$ is a polynomial 
and $\alpha \in \Bbb C$. Quasipolynomials have the following well-known 
simple characterization:
a formal power series is a quasipolynomial if and only if 
it satisfies a non-trivial linear differential equation 
with constant coefficients. 
We have the following characterization of quasi-finiteness of
an irreducible module $ L(\xi)$. 
\begin{theorem}
 A $\hBpm$-module $L(\hBpm; \xi)$ is quasifinite if and only if
 $$
   \Delta^{\pm} (x) = \frac{ F(x)}{2 \sinh (x/2 )}
 $$
 where $F(x)$ is an even quasipolynomial such that $F(0) = 0$.
   \label{th_delta}
\end{theorem}
\begin{demo}{Proof}
 We prove the $+$ case first.
   It follows from Propositions \ref{prop_commutat} and \ref{prop_equiv}
that $L(\hBp;  \xi)$ is quasifinite if and only if there 
 exist two monic odd polynomials $b_1(w)$ and $b_2(w)$ such that 
 for all $l \geq 1$ the following two equations hold:
 \begin{eqnarray}
  \xi ( D^{2l -1} b_1( D  ) - (D + 1)^{2l -1} b_1 ( D + 1) ) & = & 0 ,
  \label{eq_crit1}
 \end{eqnarray}
 \begin{eqnarray}
 \begin{array}{rcl}
  \lefteqn{
    \xi \left( (D- 1/2 )^{2l -1} b_2 ( D - 1/2 ) \right.  } \\
      &&   \left.     - (D + 3/2)^{2l -1} b_2 ( D + 3/2)  
       + 2 ( 1/2)^{2l -1} b_2 (  1/2) C 
      \right) 
       =  0.
  \label{eq_crit2}
 \end{array}
 \end{eqnarray}
 Let
 $$
  b_1 (w) = \sum_{n=0}^M a_n w^{2n+1},
  \quad b_2(w) = \sum_{n=0}^N c_n w^{2n+1}.
 $$
 Then we can rewrite equations (\ref{eq_crit1}) and (\ref{eq_crit2}) as
 \begin{eqnarray}
  \sum_{n=0}^M a_n \left( \sum_{i \in {\Bbb N}_{odd}} { {2n+2l} \choose i}
                      \left( \frac 12 
                      \right)^{2n+2l-i} \Delta^+_i 
                  \right) & = & 0   \label{eq_expd}   \\
 \sum_{n=0}^M c_n \left( \sum_{i \in {\Bbb N}_{odd} } { {2n+2l} \choose i}
          \Delta^+_i + c \left( \frac 12
                       \right)^{2n+2l} 
                  \right) & = & 0.
  \label{eq_explic}
\end{eqnarray}
Let $F(x) = \Delta^+ (x) \sinh{\frac x2}$ 
and $G(x) = \Delta^+ (x) \sinh{x} + c \cosh{\frac x2}$.
It is straightforward to check that equations (\ref{eq_expd})
and (\ref{eq_explic}) can be equivalently reformulated as follows:
\begin{eqnarray*}
  (\sum_{n=0}^M a_n (\frac d{dx})^{2n+1}) F (x) & = & 0  \\
  (\sum_{n=0}^N c_n (\frac d{dx})^{2n+1}) G(x) & = & 0.
\end{eqnarray*}
Since $G(x) = (2F(x) + c) \cosh{\frac x2}$, 
we see that $L(\hBp; \xi)$ is quasifinite if and only if
$F(x)$ is an even quasipolynomial. 

In the $-$ case, it follows from 
Propositions \ref{prop_commutator} and \ref{prop_equiv}
 that $L(\hB; \xi)$ is quasifinite if and only if there exists 
 a monic even polynomial $ b(w) = \sum_{n=0}^M {\beta}_n w^{2n}$ such that 
 $$
  \xi \left( (D- 1/2)^{2l} b (D - 1/2)
            -(D + 1/2)^{2l} b(D + 1/2) + ( - 1/2)^{2l} b( - 1/2)C
          \right) = 0
 $$
 for all $l \in \Z_{+}.$ As in the $+$ case, one can show that
 this condition can be reformulated as 
 $$
  \left(\sum_{n=0}^M {\beta}_n \left( {d^2}/{dx^2} \right)^n 
                    \right) H(x) = 0.
 $$
 where
 \begin{eqnarray} 
   H (x) = 2  \Delta^- (x) \sinh (x/2) + c \cosh (x/2), \quad H(0) = c.
    \label{eq_fg}
 \end{eqnarray}
 Thus $L(\xi)$ is quasifinite if and only if $ H(x)$ is an even 
 quasipolynomial such that $H (0) = 0$. By (\ref{eq_fg}), letting
 $ F(x)  = H (x) - c \cosh (x/2)$ completes the proof.
\end{demo}
        
From the proof of Theorem \ref{th_delta} we obtain the following corollary.
\begin{corollary}
    Let $L(\hBp; \xi)$ (resp. $L(\hB; \xi)$) be an irreducible 
 quasifinite highest weight module 
 over $\hBp$ (resp. $\hB$) with $b(w)$ (resp. $b_1 (w)$)
 as its first characteristic polynomial. Then
 $$ F(x) = \Delta^+ (x) \sinh{\frac x2},
  \quad  H (x) = 2 \Delta^- (x) \sinh (x/2) + c \cosh (x/2) $$
 are even quasipolynomials.
 Let $F^{(M)} + a_{M-1} F^{(M-1)} + \cdots + a_0 = 0$ be 
 the minimal order linear differential equation with constant 
 coefficients satisfied by $F(x)$ such that 
 $w^M + a_{M-1} w^{M-1} + \ldots + a_0 $ is an odd polynomial. 
 Let $H^{(N)} + {\beta}_{N-1} H^{(N-1)} + \cdots + {\beta}_0 = 0$ be 
 the minimal order linear differential equation with constant 
 coefficients satisfied by $H(x)$ such that 
 $w^N + {\beta}_{N -1} w^{N-1} + \ldots + {\beta}_0 $ is an even polynomial. 
 Then $b(w) =w^N + {\beta}_{N -1} w^{N-1} + \ldots + {\beta}_0,
  b_1(w) = w^M + a_{M-1} w^{M-1} + \ldots + a_0 $.
    \label{cor_diff}
\end{corollary}

Given a quasifinite irreducible highest weight $\hBpm$-module $V$ with
 central charge $c \in \C$ and with $\Delta (x)$ 
as in Theorem~\ref{th_delta}, write $F(x) + c$ in the $+$ case
or $F(x) + c \cosh (x/2)$ in the $-$ case
as a finite sum of the form
 \begin{eqnarray}
  \sum_i p_i (x) \cosh (e_i^+ x) + \sum_j q_j (x) \sinh (e_j^- x),
   \label{eq_ka}
 \end{eqnarray}
where $p_i (x)$ (resp. $q_j (x)$) are non-zero even (resp. odd) polynomials
and $e_i^+$ (resp. $ e_j^-$) are distinct complex
numbers. Clearly 
\begin{eqnarray}
 \sum_i p_i (0) = c. 
   \label{eq_charge}
\end{eqnarray}
 The expression (\ref{eq_ka})
is unique up to a sign of $e_i^+$ or a simultaneous change of signs of
$ e_j^-$ and $q_j (x)$.
We call $e_i^+$ (resp. $e_j^-$) the {\em even type} (resp. {\em odd type})
{\em exponents} of $V$ with {\em multiplicities}
$ p_i (x)$ (resp. $ q_j (x)$). We denote by $e^+$ 
the set of even type exponents $e_i^+$ 
with multiplicity $p_i (x)$ 
and by $e^-$ the set of odd type exponents $e_j^-$ 
with multiplicity $q_j (x)$.
Then the pair $ (e^+, e^-)$ determines $V$ uniquely. 
We shall therefore denote this module by  $L( \hBpm; e^+, e^-)$.
As we shall see, the following class of $\hBpm$-modules is
especially important. 
\begin{definition}
  A quasifinite irreducible highest weight $\hBpm$-module $V$ with
 central charge $c \in \C$ is called {\em primitive} if the multiplicities
 of its exponents $e_i^+$ are nonzero constants $n_i \in \C$
 and $e^- = \varnothing$. 
 A primitive $\hB$-module $V$ is called {\em positive} if
 $n_i \in \Bbb N$ when $e_i \neq \pm \hf$ and 
 $n_i \in \hf \Bbb N$ when $e_i = \pm \hf$.
 A primitive $\hBp$-module $V$ is called {\em positive} if
 $n_i \in \Bbb N$ when $e_i \neq  0$ and
 $- \hf n_{i_0} \leq n_i \in \hf \Z$ when $e_i =  0$, where
 ${i_0}$ is the index such that $e_{i_0} = 1$.
 (In both $\pm$ cases the central charge $c = \sum_i n_i \in \hf \Bbb N$).
\end{definition}

  It is convenient to make the following convention.
\begin{convention}
 For a primitive module $V$ of $\hBpm$,
 we let $e$ stand for the set of (even type) exponents
 with their multiplicities in $e^+$ which are not equal to $\pm \hf$ 
 (resp. 0) in the $-$ (resp. $+$) case.
 The pair $(e, c)$ determines uniquely the module $V$. We will
 denote this primitive module by $L( \hBpm; e, c)$.
   \label{conven_1}
\end{convention}
\section{Embedding of $\hB$ into infinite rank classical Lie algebras}
  \label{sec_embed}

Let $\co $ be the algebra of all holomorphic functions on $\C$ 
with topology of uniform convergence on compact sets. Denote 
\begin{eqnarray*}
  \Oo & = & \{ f \in {\co} \mid f(w) = - f( -w) \}  \\
  \Oe & = & \{ f \in {\co} \mid f(w) = f(-w) \}.
\end{eqnarray*}
We define a completion $\D^{\cal O}$ of $\D$ consisting  of 
all differential operators of the form $ t^j f(D)$ where
$f \in \co$ and $j \in \Z$. We similarly
define a completion $\BO$ (resp. $\BOp$) of $\B$ (resp. $\Bp$) consisting  of 
all differential operators of the form $ t^j f(D)$ 
(resp. $ t^j f( D + j/2)$) where
$f \in \co$ (resp. ${\co}^{ ( \bar{j} ) }$) and $j \in \Z$.

Note that $ t^j f(D)$ acts on $\C [t, t^{-1}]$ by
$$
 t^j f(D) t^k = f(k)t^{k +j}.
$$
Formula (\ref{cocy}) for the 2-cocycle $\Psi$ on $\D$ (resp. $\B$) 
extends to a 2-cocycle 
on $\D^{\cal O}$ (resp. $\BO$). We denote  
the corresponding central extension by 
$\hD^{\cal O} =\D^{\cal O} \bigoplus \C C$
and $\hBO =\BO \bigoplus \C C$. The commutation relations (\ref{eq_12})
extend as well.

The vector space $R_m [t, t^{-1}] t^s \; (s \in \C)$ has a basis
$v_i = t^{ -i+s}$ $(i \in \Z)$ over $R_m$. 
The Lie algebra $\glm$ acts on this vector
space by (\ref{eq_natur}). The Lie algebras
$\D$ and $\D^{\cal O}$ 
also act on $R_m [t, t^{-1}]$ naturally as differential operators.
In this way we obtain a family of embeddings
$\phism$ of the Lie algebra $\D$ (resp. $\DO$) to $\glm$ defined by
$$ \phism ( t^k f( D) ) = \sum_{j \in \Z}
        f (-j +s +u)  E_{j-k, j}
= \sum_{i =0}^m \sum_{j \in \Z}
        \frac{ f^{(i)} (-j +s ) }{i!} u^i E_{j-k, j}
$$
where $f^{(i)}$ denotes the $i$-th derivative.
When restricted to $\B$ and $\BO$, we have
\begin{eqnarray}
  \phism \left( t^k f( D + k/2) 
          \right) = \sum_{i=0}^m \sum_{j\in \Bbb Z}
      \frac { f^{(i)}(-j +k/2 +s) }{i!} u^i E_{j-k, j}
        \mbox{ for } f \in \C [w]^{( \bar{k})}.
  \label{eq_emb}
\end{eqnarray}
When restricted to $\Bp$ and $\BOp$, we have
\begin{eqnarray}
  \phism \left( t^k f( D + (k +1)/2) 
          \right) = \sum_{i=0}^m \sum_{j\in \Bbb Z}
      \frac { f^{(i)}(-j +(k +1)/2 +s) }{i!} u^i E_{j-k, j}
  \label{eq_embd}
\end{eqnarray}
for $ f \in \Oo.$

\begin{remark} \rm
  The principal $\Z$-gradations on $\Bpm$ and $\glm$ are compatible
 under the homomorphisms $\phism$.
\end{remark}
Let 
\begin{eqnarray*}
 I_{s,k}^{[m], -} & = & \{ f \in {\co}^{ ( \bar{k} )}
    \mid f^{(i)} (n + k/2 + s) =0  \mbox{ for all } n \in \Bbb Z,
          i=0,1,...,m \}  \\
I_{s,k}^{[m], +} & = & \{ f \in \Oo
    \mid f^{(i)} (n + (k +1)/2 + s) =0  \mbox{ for all } n \in \Bbb Z,
          i=0,1,...,m \}
\end{eqnarray*}
and let
\begin{eqnarray*}
 J_s^{[m], -} & = & \bigoplus_{k \in \Bbb Z}
                 \{ t^k f(D + k/2) \mid f \in I_{s,k}^{[m], -} \},   \\
 J_s^{[m], +} & = & \bigoplus_{k \in \Bbb Z}
                 \{ t^k f(D + (k +1)/2) \mid f \in I_{s,k}^{[m], +} \}.
\end{eqnarray*}
Now fix $\vec{s} = ( s_1, s_2, \ldots, s_N) \in \C^N$, 
$s_i - s_j \notin \Z$ if $i \neq j$, $s_i + s_j \notin \Z$ for all $i,j$. 
Also fix $\vec{m} = ( m_1, m_2, \ldots, m_N ) \in \Z_+^N$. Let
$$
 \gl^{[\vec{m}]} = \bigoplus_{i=1}^N \gl^{[m_i]}
$$
and consider the homomorphism
$$ \phivsm = \bigoplus_{i=1}^N \phi_{s_i}^{[m_i], \pm}: 
   \BOpm \longrightarrow \gl^{[\vec{m}]}.
$$
\begin{proposition}
     Given $\vec{s}$ and $\vec{m}$ as above,
  we have the following exact sequence of Lie algebras:

 $$ 0 \longrightarrow J_{\vec s}^{[\vec m], \pm} \longrightarrow \BOpm
     \stackrel{\phivsm}{\longrightarrow} \gl^{[\vec{m}]} \longrightarrow 0
 $$
 where
 $J_{\vec{s} }^{ [\vec{m} ], \pm } = \cap_{i=1}^N J_{s_i}^{ [m_i], \pm}.$
  \label{prop_surj}
\end{proposition}
\begin{demo}{Proof}
  We will prove the proposition 
 only in the $-$ case. The proof in the $+$ case is parallel.
  For the sake of simplicity of notations, we prove it
 in the case $N =1$: $\vec{m} =m \in \Z_{+}$ and $ \vec{s} = s \in \C$
 ($s \notin \Z/2$ by the assumption on $\vec{s}$). The general
 case is similar.

   It is clear from the definition of $J_s^{[m]}$ that
 $\ker \phism = J_s^{[m]} $. 
 To show the surjectivity of $\phism$, it suffices to find a preimage of 
 $$ g = \sum_{i=0}^m \sum_{j \in \Bbb Z} \frac{ p_{ij} }{i!}
      u^j E_{j-k, j} \in \glm \quad ( p_{ij} \in \C)
 $$
 for a fixed $k \in \Z$. 
 We need to quote the following well-known
 theorem: for every discrete sequence of points in $\C$ and a non-negative 
 integer $m$ there exists $p(w) \in \cal O$ having prescribed values
 of its first $m$ derivatives at these points.

 Since $s \notin \Bbb Z/2$, the two sequences
 $\{ -j + \frac k2 +s\}_{j \in \Bbb Z}$
 and $\{ j -\frac k2-s\}_{j \in \Bbb Z}$ are disjoint. Thus
 there exists $p(w) \in \cal O$ such that
 \begin{eqnarray*}
   p^{(i)} ( -j + k/2 + s) = p_{ij}, \quad
   p^{(i)} ( j - k/2 - s ) = (-1)^{k+i+1}p_{ij}.
 \end{eqnarray*}
 Now  let 
 $$ f(w) = \frac{ p(w) - (-1)^k p(-w) }2 \in {\cal O}^{(\bar{k})}.
 $$
 Then $t^k f( D + k/2)$ is the preimage of $g$ 
 via $\phism$.
\end{demo}

Now we want to extend the homomorphism $\phism$ to 
a homomorphism between the central extensions of the corresponding
Lie algebras. Introduce the following functions:
$$ \eta_j (x; \mu ) = \frac{1}{2 j!} x^j (e^{ \mu x} + (-1)^j e^{ - \mu x})
    \quad (j \in \Z_{+}, \mu \in \C ).
$$
The functions $\eta_j (x; \mu )$ satisfy:
\begin{equation}
   \eta_j (-x; \mu ) = \eta_j (x; \mu ) , \quad 
   \eta_j (x; -\mu ) = (-1)^j \eta_j (x; \mu ),
     \quad  \eta_0 (x; \mu ) = \cosh (\mu x).
 \label{eq_property}
\end{equation}
 Note that, being an even quasipolynomial, 
  $F(x)$ in Theorem~\ref{th_delta} is a finite
 linear combination of the functions $ \eta_j (x; \mu)$.
\begin{proposition}
 \begin{enumerate}
  \item[1)] 
 The $\C$-linear map $\hphism: \hBp \longrightarrow \hglm$ defined by
 \begin{eqnarray}
  \phism (C) & = & 1     \label{eq_cen} \\
  \hphism |_{\hBp_j} & = & \phism |_{\Bp_j} \quad \mbox{ if } j \neq 0 
                      \nonumber \\
   \hphism \left( e^{x(D + 1/2)} - e^{ -x(D + 1/2)} \right)
   & = & \phism \left( e^{ x(D + 1/2)} - e^{ -x(D + 1/2)} \right) \nonumber \\
   & &   - \frac{\cosh {sx}-1}{\sinh{\frac x2}} {\bf 1}
     - \sum_{j=1}^m \frac{u^j \eta_j(x; s)}{2\sinh{\frac x2}}
       \label{eq_cent}
 \end{eqnarray}
 is a homomorphism of Lie algebras over $\Bbb C$.
  \item[2)] The $\Bbb C$-linear map $\hphism: \hB \longrightarrow \hglm$ 
 defined by
 \begin{eqnarray}  
  \hphism( C)  =  {\bf 1} \in R_m , \quad
  \hphism \mid_{ {\hB}_j} = \phism \mid_{\B_j} \quad (j \neq 0 ),
    \label{eq_center}
 \end{eqnarray} 
 \begin{eqnarray}
   \hphism \left( e^{xD} - e^{-xD} \right)
      & = & \phism \left( e^{xD} - e^{-xD} \right)   
       - \frac{\cosh (s- 1/2)x - \cosh (x/2)}{\sinh (x/2)}{\bf 1} \nonumber \\
      & & - \sum_{j =1}^m \frac{u^j \eta_j (x; s -1/2 ) }{\sinh (x/2)}
    \label{eq_central}
 \end{eqnarray}
 is a homomorphism of Lie algebras over $\Bbb C$. 
 \end{enumerate}
   \label{prop_ext}
\end{proposition}
\begin{demo}{Proof}
 We will prove part 1). The proof of part 2) is similar.
Part 1) follows directly from the computation of the difference 
$$ \hphism \left( e^{x(D + 1/2)} - e^{ -x(D + 1/2)} \right) $$
by using the following lemma (also cf. Proposition 4.4, \cite{KR1}).
\end{demo}

\begin{lemma}
    The $\C$-linear map $\hphism: \hD \longrightarrow \hglm$ defined by
 \begin{eqnarray}
  \phism (C) & = & 1  \nonumber \\
  \hphism |_{\hD_j} & = & \phism |_{\D_j} \quad \mbox{ if } j \neq 0 
                 \nonumber       \\
   \hphism \left( e^{x(D + 1/2)} \right)
   & = & \phism \left( e^{ x(D + 1/2)}  \right) \nonumber \\
   & & - \frac{e^{(s + 1/2) x} - e^{x/2} }{e^x -1} {\bf 1}
     -  \sum^m_{i = 1} \frac{x^i e^{(s + 1/2) x}}{e^x - 1} u^i / i!.
   \label{eq_homom}
 \end{eqnarray}
 is a homomorphism of Lie algebras over $\Bbb C$.
  \label{lem_exten}
\end{lemma}
\begin{demo}{Proof}
 Introduce two formal variables $\alpha, \beta$
 and let $ x = \alpha + \beta$. It suffices to check
 that for a given $k \in \Bbb N$ 
 \begin{eqnarray}
   \left[ \hphism (t^k e^{ \alpha (D + 1/2)}), 
   \hphism (t^{ -k} e^{ \beta (D + 1/2)}) 
  \right]
   = \hphism  \left[ (t^k e^{ \alpha (D + 1/2)}), 
                 (t^{ -k} e^{ \beta (D + 1/2)}) 
             \right].
 \label{eq_equal}
\end{eqnarray}

 By a straightforward computation using (\ref{eq_embd}) we obtain
 \begin{eqnarray}
  \begin{array}{rcl}
   \lefteqn{
  \left[ \hphism (t^k e^{ \alpha (D + 1/2)}), 
   \hphism (t^{ -k} e^{ \beta (D + 1/2)}) 
  \right]  } \\
   & = &
  \left[ \phism (t^k e^{ \alpha (D + 1/2)}), 
   \phism (t^{ -k} e^{ \beta (D + 1/2)}) 
  \right]                              \\
  & = & \left( e^{ - \alpha k} - e^{  \beta k  }
    \right)
     \left( \sum_{j \in \Z}
      e^{ x (-j +s + 1/2)} E_{j,j}
    + e^{xu} e^{ x (s + 1/2)} \frac 1{1 - e^x} \right)     \\
  & = & \left( e^{ - \alpha k} - e^{  \beta k  }
    \right)
     \left( \phism \left( e^{ x (D + 1/2)}
                   \right)
          + \frac {e^{ x (s + 1/2)}}{1 - e^x} {\bf 1}
          + \sum_{j+1}^m \frac{x^i e^{ x (s + 1/2)}} {1 - e^x} u^i / j!
    \right). 
  \end{array}
  \label{eq_sutb}
\end{eqnarray}
 On the other hand, using (\ref{eq_12}) we have
\begin{eqnarray}
     \left[ t^k e^{ \alpha (D + 1/2)}, t^{ -k} e^{ \beta (D + 1/2)}
    \right]
  =  \left( e^{ - \alpha k} - e^{  \beta k  }
    \right)
     \left( e^{ x (D + 1/2)} + \frac {e^{ x/2}}{1 - e^x} C
    \right).
  \label{eq_calcu}
\end{eqnarray}
 By applying $\hphism$ to (\ref{eq_calcu}) and then comparing with
 (\ref{eq_sutb}) with the help of (\ref{eq_homom}),
 we obtain (\ref{eq_equal}). 
\end{demo}

The homomorphism $\phism: \BO \longrightarrow \glm$ is defined for 
any $s \in \C$. However for $s \in \Z/2$, it is no longer surjective.
The case $ s =0$ is described by the following proposition.
\begin{proposition}
    We have the following exact sequence of Lie algebras:
 \begin{eqnarray*}
   0 \longrightarrow J_0^{[m], -} \longrightarrow \BO
     \stackrel{\phiom}{\longrightarrow} \Binfm \longrightarrow 0 \\
  0 \longrightarrow J_0^{[m], +} \longrightarrow \BOp
     \stackrel{\phiom}{\longrightarrow} \Dinfm \longrightarrow 0.  \\
 \end{eqnarray*}
 \label{prop_phizero}
\end{proposition}
\begin{demo}{Proof}
  By the definition of $\phiom$, it is easy to see that
 the image of $\phiom$ lies in $\Binfm$ (resp. $\Dinfm $). The proof of the
 rest of the proposition is similar to that of Proposition~\ref{prop_surj}.
\end{demo}

Similarly for $s = \hf$ we have the following proposition.
\begin{proposition}
   We have the following exact sequence of Lie algebras:
 \begin{eqnarray*}
  0 \longrightarrow J_{1/2}^{[m], -} \longrightarrow \BO
     \stackrel{\phim}{\longrightarrow} \Cinfm \longrightarrow 0 \\
  0 \longrightarrow J_{ -1/2}^{[m], +} \longrightarrow \BOp
     \stackrel{\phicm}{\longrightarrow} \Binfmtwo \longrightarrow 0.
 \end{eqnarray*}
   \label{prop_phihalf}
\end{proposition}
\begin{remark} \rm
  For $ s \in \Z$, the image of $\hB$ (resp. $\hBp$) under 
 the homomorphism $\phism$ is $\nu^s (\Binfm)$ (resp. $\nu^s (\Dinfm)$)
(recall that $\nu$ is defined in (\ref{eq_autom}) ). 
 For $ s \in \hz$, the image of $\hB$ (resp. $\hBp$) under $\phism$ is 
 $\nu^{s - 1/2} (\Cinfm)$ (resp. $\nu^{s + 1/2} (\Binfmtwo)$). 
 Hence we will only need to consider $ s = 0, \hf$ in the $-$ case
 and  $ s = 0, -\hf$ in the $+$ case
 whenever $s \in \Z/2$ throughout the paper.
 Note that the principal $\Z$-gradation
 of $\hBp$ is compatible with the gradation of type $(2, 1, 1, \ldots)$
 on $\Dinf$ via the homomorphism $\phi_0$.
\end{remark}
\begin{definition}
 We say that the vector 
 $\vec{s} = ( s_1, s_2, \ldots, s_N) \in \C^N$ satisfies
 the ($\star -$) (resp. ($\star +$)) condition if 
 $s_i \in \Z$ implies $s_i = 0$, $s_i \in \hz$ implies $s_i = \hf$
 (resp. $- \hf$), and $s_i \neq \pm s_j \bmod{\Z} $ for $i \neq j$.
   \label{def_star}
\end{definition}

 Given $\vec m =(m_1, \ldots, m_N) \in \Z_+^N $ and
$\vec s = (s_1, s_2, \ldots, s_N) \in \C^N$ satisfying
the ($\star -$) (resp. ($\star +$)) condition, 
combining Propositions \ref{prop_surj}, 
\ref{prop_phizero} and \ref{prop_phihalf}, we obtain 
a homomorphism of Lie algebras over $\Bbb C$:
 \begin{equation}
   \hphivsm = \bigoplus_{i=1}^N \widehat{\phi}_{s_i}^{[m_i]}: 
   \hBpm \longrightarrow \g^{ [\vec m]} := \bigoplus_{i=1}^N \g^{[m_i]} ,
    \label{eq_map}
 \end{equation}
where {\bf the following
consistency condition is always assumed throughout the paper}:
in the $-$ case, 
\begin{eqnarray*}
  \g^{[m_i]} = 
   \left\{
      \everymath{\displaystyle}
      \begin{array}{ll}
        \binf^{[m_i]}, & \mbox{if } s_i = 0      \\
        \cinf^{[m_i]}, & \mbox{if } s_i = \hf    \\
        \hgl^{[m_i]} , & \mbox{if } s_i \notin \Bbb Z/2,
      \end{array}
    \right. \\
\end{eqnarray*}
while in the $+$ case,
\begin{eqnarray*}
  \g^{[m_i]} = 
   \left\{
      \everymath{\displaystyle}
      \begin{array}{ll}
        \dinf^{[m_i]}, & \mbox{if } s_i = 0      \\
        \binftwo^{[m_i]}, & \mbox{if } s_i = - \hf    \\
        \hgl^{[m_i]} , & \mbox{if } s_i \notin \Bbb Z/2.
      \end{array}
    \right. \\
\end{eqnarray*}
Furthermore, we can prove the following proposition
in the same way as Proposition~\ref{prop_surj}.
\begin{proposition}
 The homomorphism $\hphivsm$ extends to a surjective homomorphism of 
 Lie algebras over $\Bbb C$ which is denoted again by $\hphivsm$:
 \begin{eqnarray*}
   \hphivsm = \bigoplus_{i=1}^N \widehat{\phi}_{s_i}^{[m_i]}: 
   \hBOpm \longrightarrow \g^{ [\vec m]}.
 \end{eqnarray*}
   \label{prop_phivector}
\end{proposition}
\section{Realization of QHWM's of $\hBpm$}
  \label{sec_realization}
     Let $\gm$ stand for $\hglm$, or one of its classical Lie subalgebras.
The proof of the following simple proposition is standard.
\begin{proposition}
   The $\gm$-module $L(\gm; \Lambda)$ is quasifinite if and only if 
 all but finitely many of the ${}^{*}h_k^{(j)} $ are zero, where $*$ 
 represents $a, b, c$ or $ d$ depending on whether
 $\gm$ is $\hglm$, $\binfm$, $\binfmtwo$, $\cinfm$ or $\dinfm$.
   \label{prop_qfin}
\end{proposition}

Take a quasifinite $\Lambda(i) \in (\frak g_0^{[m_i]})^*$ 
for each $i = 1, \ldots, N$,
and let $L \left( \frak g^{[m_i]}; \Lambda(i) \right) $ 
be the corresponding irreducible $\g^{[m_i]}$-module. 
Then the outer tensor product
$$
 L \left( \gvm; \vec{\Lambda} \right) \equiv \bigotimes_{i=1}^N 
  L\left( \g^{[m_i]}; \Lambda(i) \right)
$$
is an irreducible $\gvm$-module.
The module $L \left(\gvm; \vec{\Lambda} \right)$ can be regarded
as a $\hBpm$-module via the homomorphism $\hphivsm$ 
given by (\ref{eq_map}), 
which we shall denote by $L_{\vec s}^{[\vec m], \pm}(\vec{\Lambda} )$.

We will need a technical lemma whose proof is analogous to that
of Proposition~4.3 in \cite{KR1}.

\begin{lemma}
  Let $V$ be a quasifinite $\hBpm$-module. Then the action of $\hBpm$ 
 on $V$ naturally extends to the action of $\hBOpm_k$ on $V$ 
 for any $k \neq 0$.
   \label{lem_tech}
\end{lemma}

\begin{theorem}
   Let $V$ be a quasifinite $\g^{[\vec m]}$-module, which
 can be regarded as a quasifinite 
 $\hBpm$-module via the homomorphism $\hphivsm$. 
 Then any $\hBpm$-submodule of $V$ is also a $\gvm$-submodule. 
 In particular, the $\hBpm$-modules
 $L_{\vec s}^{[\vec m], \pm}(\vec \Lambda)$ 
 are irreducible if $\vec s = (s_1, s_2, \ldots, s_N)$
 satisfies the $(\star \pm)$ condition.
   \label{th_pullback}
\end{theorem}
\begin{demo}{Proof}
   Take any $\hBpm$-submodule $W$ of $V$. By the Lemma \ref{lem_tech} 
 we can extend the action of $\hBpm$ to $\hBOpm_j \;(j \neq 0)$.
 Then by Proposition \ref{prop_phivector}, we see that the subspace
 $W$ is preserved by the action of the graded subspace
 $\gvm_j \; (j \neq 0)$ of $\gvm$. Since $\gvm$
 coincides with its derived algebra, $W$ is preserved by the action 
 of the whole $\gvm$.
\end{demo}

We will show that in fact all the quasifinite $\hBpm$-modules can be 
realized as some $L_{\vec s}^{[\vec m], \pm}(\vec{\Lambda} )$.
But first let us calculate the generating function
$\Delta^{\pm}_{m, s, \Lambda} (x)$ of highest weight
for $L_{\vec s}^{[\vec m], \pm}(\vec{\Lambda} )$ in some typical cases.
\begin{proposition}
    Consider the embedding
 $\hphism : \hB \longrightarrow \hglm$ with $ s \notin \Z/2$.
 The $\hglm$-module $L( \hglm; \Lambda)$ regarded as 
 a $\hB$-module is isomorphic to $L(\hB; e^+, e^-)$ where 
 $e^+$ and $e^-$ consist of the 
 exponents $s -i -1/2$ $(i \in \Z)$ with multiplicities 
 \begin{eqnarray}
  \sum_{ 0 \leq j \leq m,\; j\, {\small even}} {}^a h_i^{(j)} x^j / j! 
       \quad and \quad
    \sum_{ 0 \leq j \leq m,\; j\, {\small odd}} {}^a h_i^{(j)} x^j / j!
 \end{eqnarray}
 respectively (see Section \ref{subsec_glinf} for notations; the exponents
 with zero multiplicities are dropped).
   \label{prop_deltagl}
\end{proposition}
\begin{demo}{Proof}
    By Theorem \ref{th_pullback} and Proposition \ref{prop_qfin}
 the $\B$-module $L_{\vec s}^{[\vec m], -}(\vec{\lambda} )$ 
 is an irreducible quasifinite highest weight module. 
 By formula (\ref{eq_center}) we see the central charge $c = c_0$.
 By applying $\Lambda$ to (\ref{eq_central}) and using formulas 
 (\ref{eq_triv}) and (\ref{eq_emb}) we obtain:
 \begin{eqnarray*}
   2 \Delta^-_{m, s, \Lambda} (x) & = &
    -  \sum_{j=0}^m \sum_{i \in \Z}
          ({}^a \lambda_i^{(j)}/ j! )x^j 
            \left( e^{(s -i)x} - (-1)^j e^{(i -s)x}
            \right)        \\
    & \quad & + c_0 \frac{\cosh (s- 1/2)x - \cosh (x/2)}{\sinh (x/2)} 
    + \sum_{j =1}^m \frac{ c_j \eta_j (x; s -1/2 ) }{\sinh (x/2)}
                                     \\
   & =  & - \sum_{j=0}^m \sum_{i \in \Z}
          ({}^a \lambda_i^{(j)}/ j! )
            \frac{ (\eta_j (x; s -i + 1/2) - \eta_j (x; s-i - 1/2) ) }
                 {\sinh (x/2)}   \\
   & \quad & + c_0 \frac{\cosh (s- 1/2)x - \cosh (x/2)}{\sinh (x/2)} 
    + \sum_{j =1}^m \frac{c_j \eta_j (x; s -1/2 ) }{\sinh (x/2)}
                                 \\
   & \stackrel{(1)}{=} &    \frac{ 
                      \sum_{j=0}^m \sum_{i \in \Z}
                       {}^a h_i^{(j) } \eta_j (x; s -i -1/2)} 
                          {  \sinh (x/2)}     
                -  \frac{c_0 \cosh (x/2)}{ \sinh (x/2)} .
 \end{eqnarray*}
 The identity (1) above is obtained by shifting the index $i$ to $i+1$
 in the first half of the first summation of the left hand side of (1).
 Now the proposition
 follows from the definition of exponents and their multiplicities.
\end{demo}
\begin{proposition}
    Consider the embedding
 $\hphiom : \hB \longrightarrow \binfm$. 
 The $\binfm$-module $L(\binfm; \Lambda)$ regarded as 
 a $\hB$-module is isomorphic to $L(\hB; e^+, e^-)$ where 
 $e^+$ and $e^-$ consist of the
 exponents $- i -1/2$ $(i \in \Z_{+})$ with multiplicities 
 \begin{eqnarray}
  \sum_{ 0 \leq j \leq m,\;j \,{\small even}}
  {}^b \widetilde{h}_i^{(j)} x^j / j! 
       \quad and \quad
  \sum_{ 0 \leq j \leq m, \;j \, {\small odd}}
  {}^b \widetilde{h}_i^{(j)} x^j / j!
 \end{eqnarray}
 respectively,
 where ${}^b\tilde{h}_i^{(j)} = {}^bh_i^{(j)}$ $( i > 0)$ and
 ${}^b \tilde{h}_0^{(j)} = \hf \; {}^b h_0^{(j)}$
 (see Section \ref{subsec_binf} for  notations).
   \label{prop_deltab}
\end{proposition}
\begin{demo}{Proof}
 We will only need to calculate $\Delta^-_{m, s, \Lambda} (x)$. The rest
 of the statement is clear, cf. the proof of Proposition \ref{prop_deltagl}.
 We have:
 \begin{eqnarray*}
   2 \Delta^-_{m, s, \Lambda} (x) 
   & \stackrel{(1)}{=} & - \sum_{j=0}^m \sum_{i \in \Z}
          ( \Lambda( u^j E_{ii}) / j! )
            \frac{ (\eta_j (x; -i + 1/2) - \eta_j (x; -i - 1/2) ) }
                 {\sinh (x/2)}   \\
     &\quad & + \sum_{j =1}^m \frac{  c_j \eta_j (x;  -1/2 ) }{\sinh (x/2)}
                                 \\
   & \stackrel{(2)}{=} &    \frac{ 
                      \sum_{j=0}^m \sum_{i \in \Z}
                         \Lambda ( u^j E_{ii} - u^j E_{i+1, i+1})
                           \eta_j (x; -i -1/2) 
                  }     { \sinh (x/2)}     \\
   &\quad &  + \sum_{j =1}^m \frac{c_j \eta_j (x;  -1/2 ) }{ \sinh (x/2)}   \\
    & \stackrel{(3)}{=} & 
                     \frac{ 
                      \sum_{j=0}^m \sum_{i \in \Z_{+} } 
                        {}^b \tilde{h}_i^{(j)} \eta_j (x;  -i -1/2)
                           }   
            { \sinh (x/2)}  
   -  \frac{c_0 \cosh (x/2)}{ \sinh (x/2)} .
 \end{eqnarray*}
 The identity (1) is obtained from a similar identity in the proof of
 Proposition \ref{prop_deltagl} by putting $s = 0$.
 The identity (2) is obtained by shifting the index $i$ to $i+1$
 in the first half of the first summation of the left hand side of (2).
 The identity (3) is obtained by splitting the summation into two,
 $ \sum_{i \in \Z} = \sum_{i \in \Z_{+}} + \sum_{i \in - \Bbb N}$,
 changing the index $i$ to $-i-1$ in the second summation,
 and then using formulas (\ref{eq_bweight}) and (\ref{eq_property}).
\end{demo}
\begin{proposition}
    Consider the embedding
 $\hphim : \hB \longrightarrow \cinfm$. 
 The $\cinfm$-module $L(\cinfm; \Lambda)$ regarded as 
 a $\hB$-module is isomorphic to $L(\hB; e^+, e^-)$ where 
 $e^+$ and $e^-$ consist of the 
 exponents $s -i -1/2$ $(i \in \Z_{+})$ with multiplicities 
 \begin{eqnarray}
  \sum_{ 0 \leq j \leq m,\; j\, {\small even}} {}^c h_i^{(j)} x^j / j! 
       \quad and \quad
    \sum_{ 0 \leq j \leq m,\; j\, {\small odd}} {}^c h_i^{(j)} x^j / j!
 \end{eqnarray}
 respectively (see Section \ref{subsec_cinf} for more notations).
   \label{prop_deltac}
\end{proposition}
\begin{demo}{Proof}
 It suffices to calculate $\Delta^-_{m, s, \Lambda} (x)$:
 \begin{eqnarray*}
   2 \Delta^-_{m, s, \Lambda} (x) 
   & \stackrel{(1)}{=} & - \sum_{j=0}^m \sum_{i \in \Z}
          ( \Lambda( u^j E_{ii}) / j! )
            \frac{ (\eta_j (x; -i + 1) - \eta_j (x; -i ) ) }
                 {\sinh (x/2)}   \\
    & \quad & + c_0 \frac{1 - \cosh (x/2)}{\sinh (x/2)} 
    + \sum_{j =1}^m \frac{ c_j \eta_j (x; 0 ) }{\sinh (x/2)}
                                     \\
    & \stackrel{(2)}{=} &    \frac{ 
                      \sum_{j=0}^m \sum_{i \in \Z}
                         \Lambda ( u^j E_{ii} - u^j E_{i+1, i+1})
                           \eta_j (x; -i ) 
                  }     {  \sinh (x/2)}  \\
    & \quad & + c_0 \frac{1 - \cosh (x/2)}{\sinh (x/2)} 
    + \sum_{j =1}^m \frac{ c_j \eta_j (x; 0 ) }{\sinh (x/2)}
                                     \\
    & \stackrel{(3)}{=} & 
                     \frac{ 
                      \sum_{j=0}^m \sum_{i \in \Z_{+} } 
                       {}^c h_i^{(j)} \eta_j (x;  -i )
                           }   
            { \sinh (x/2)}  
   -  \frac{c_0 \cosh (x/2)}{\sinh (x/2)} .
 \end{eqnarray*}
 The identity (1) is obtained from a similar identity in the proof of
 Proposition \ref{prop_deltagl} by putting $s = \hf$.
 The identity (2) is obtained by shifting the index $i$ to $i+1$
 in the first half of the first summation of the left hand side of (2).
 The identity (3) is obtained by splitting the summation into two,
 $ \sum_{i \in \Z} = \sum_{i \in \Z_{+}} + \sum_{i \in - \Bbb N}$,
 changing the index $i$ to $-i$ in the second summation,
 and then using formulas(\ref{eq_cweight}) and (\ref{eq_property}).
\end{demo}
 
\begin{proposition}
    Consider the embedding
 $\hphism : \hBp \longrightarrow \hglm$ with $ s \notin \Z/2$.
 The $\hglm$-module $L( \hglm; \Lambda)$ regarded as 
 a $\hBp$-module is isomorphic to $L (\hBp; e^+, e^- )$
 where $e^+$ and $e^-$ consist of the 
 exponents $s -i$ $(i \in \Z)$ with multiplicities 
 \begin{eqnarray}
  \sum_{ 0 \leq j \leq m,\; j\, {\small even}} {}^a h_i^{(j)} x^j / j! 
       \quad and \quad
    \sum_{ 0 \leq j \leq m,\; j\, {\small odd}} {}^a h_i^{(j)} x^j / j!
 \end{eqnarray}
 respectively (see Section \ref{subsec_glinf} for notations).
   \label{prop_dtagl}
\end{proposition}
\begin{demo}{Proof}
    By Theorem \ref{th_pullback} and Proposition \ref{prop_qfin}
 $L_{\vec s}^{[\vec m], +}(\vec{\lambda} )$ 
 is an irreducible quasifinite highest weight $\hBp$-module. 
 By formula (\ref{eq_cen}) we see the central charge $c = c_0$.
 By applying $\Lambda$ to (\ref{eq_cent}) and using formulas 
 (\ref{eq_triv}) and (\ref{eq_embd}) we obtain
 \begin{eqnarray*}
   2 \Delta^+_{m, s, \Lambda} (x) & = &
    -  \sum_{j=0}^m \sum_{i \in \Z}
          ({}^a \lambda_i^{(j)}/ j! )x^j 
            \left( e^{(s -i + 1/2)x} - (-1)^j e^{(i -s - 1/2)x}
            \right)        \\
    & \quad & + c_0 \frac{\cosh (sx) - 1}{\sinh (x/2)} 
    + \sum_{j =1}^m \frac{ c_j \eta_j (x; s ) }{\sinh (x/2)}
                                     \\
   & =  & - \sum_{j=0}^m \sum_{i \in \Z}
          ({}^a \lambda_i^{(j)}/ j! )
            \frac{ (\eta_j (x; s -i + 1) - \eta_j (x; s -i) ) }
                 {\sinh (x/2)}   \\
    & \quad & + c_0 \frac{\cosh (sx) - 1}{\sinh (x/2)} 
    + \sum_{j =1}^m \frac{ c_j \eta_j (x; s ) }{\sinh (x/2)}
                                 \\
   & \stackrel{(1)}{=} &    \frac{ 
                      \sum_{j=0}^m \sum_{i \in \Z}
                         {}^a h_i^{(j)}\eta_j (x; s -i)  }
                          {  \sinh (x/2)}     
                -  \frac{c_0 }{ \sinh (x/2)} .
 \end{eqnarray*}
 The identity (1) above is obtained by shifting the index $i$ to $i+1$
 in the first half of the first summation of the left hand side of (1).
 This completes the proof of this proposition.
\end{demo}
\begin{proposition}
    Consider the embedding
 $\hphiom : \hBp \longrightarrow \dinfm$. 
 The $\dinfm$-module $L(\dinfm; \Lambda)$ regarded as 
 a $\hBp$-module is isomorphic to 
 $L(\hB; e^+, e^-)$ where 
 $e^+$ and $e^-$ consist of the 
 exponents $- i $ $(i \in \Z_{+})$ with multiplicities 
 \begin{eqnarray}
  \sum_{ 0 \leq j \leq m,\;j \, {\small even}}
     {}^d \widetilde{h}_i^{(j)} x^j / j! 
       \quad and \quad
  \sum_{ 0 \leq j \leq m, \;j \, {\small odd}} 
     {}^d \widetilde{h}_i^{(j)} x^j / j!
 \end{eqnarray}
 respectively,
 where ${}^d \tilde{h}_i^{(j)} = {}^d h_i^{(j)}$ $( i > 0)$ and
 ${}^d \tilde{h}_0^{(j)} = \hf \; ( {}^d h_0^{(j)} - {}^d h_1^{(j)})$
 (see Section~\ref{subsec_dinf} for  notations).
   \label{prop_deltad}
\end{proposition}
\begin{demo}{Proof}
 We will only need to calculate $\Delta^+_{m, s, \Lambda} (x)$. The rest
 of the statement is clear, cf. the proof of Proposition \ref{prop_deltagl}.
 \begin{eqnarray*}
   2 \Delta^+_{m, s, \Lambda} (x) 
   & \stackrel{(1)}{=} & - \sum_{j=0}^m \sum_{i \in \Z}
          ( \Lambda( u^j E_{ii}) / j! )
            \frac{ (\eta_j (x; -i + 1) - \eta_j (x; -i) ) }
                 {\sinh (x/2)}  
 + \sum_{j =1}^m \frac{ c_j \eta_j (x; 0 ) }{\sinh (x/2)}
                                 \\
   & \stackrel{(2)}{=} &    \frac{ 
                      \sum_{j=0}^m \sum_{i \in \Z}
                         \Lambda ( u^j E_{ii} - u^j E_{i+1, i+1})
                           \eta_j (x; -i) 
                  }     { \sinh (x/2)} 
  + \sum_{j =1}^m \frac{c_j \eta_j (x;  0) }{ \sinh (x/2)}   \\
    & \stackrel{(3)}{=} & 
                     \frac{ 
                      \sum_{j=0}^m 
                       \left( \sum_{i \in \Bbb N } 
                        {}^d h_i^{(j)} \eta_j (x;  -i)
                        + \hf( h_0^{(j)} - h_1^{(j)}) \eta_j (x; 0)
                       \right)
                           }   
            { \sinh (x/2)}  
   -  \frac{c_0 }{ \sinh (x/2)} .
 \end{eqnarray*}
 The identity (1) is obtained from a similar identity in the proof of
 Proposition \ref{prop_deltagl} by putting $s = 0$.
 The identity (2) is obtained by shifting the index $i$ to $i+1$
 in the first half of the first summation of the left hand side of (2).
 The identity (3) is obtained by splitting the summation into two,
 $ \sum_{i \in \Z} = \sum_{i \in \Z_{+}} + \sum_{i \in - \Bbb N}$,
 changing the index $i$ to $-i$ in the second summation,
 and then using formulas (\ref{eq_dweight}) and (\ref{eq_property}).
\end{demo}
\begin{proposition}
    Consider the embedding
 $\hphicm : \hBp \longrightarrow \binfmtwo$. 
 The $\binfmtwo$-module $L(\binfmtwo; \Lambda)$ regarded as 
 a $\hBp$-module is isomorphic to $L(\hB; e^+, e^-)$ where 
 $e^+$ and $e^-$ consist of the 
 exponents $- i -1/2$ $(i \in \Z_{+})$ with multiplicities 
 \begin{eqnarray}
  \sum_{ 0 \leq j \leq m,\;j \, {\small even}}
   {}^b \widetilde{h}_i^{(j)} x^j / j! 
       \quad and \quad
  \sum_{ 0 \leq j \leq m, \;j \, {\small odd}}
   {}^b \widetilde{h}_i^{(j)} x^j / j!
 \end{eqnarray}
 respectively,
 where ${}^b\tilde{h}_i^{(j)} = {}^bh_i^{(j)}$ $( i > 0)$ and
 ${}^b \tilde{h}_0^{(j)} = \hf \; {}^b h_0^{(j)}$
 (see Section \ref{subsec_binf} for  notations).
  \label{prop_dtab}
\end{proposition}
\begin{demo}{Proof}
 Again it suffices to calculate $\Delta^+_{m, s, \Lambda} (x)$. 
 \begin{eqnarray*}
   2 \Delta^+_{m, s, \Lambda} (x) 
   & \stackrel{(1)}{=} & - \sum_{j=0}^m \sum_{i \in \Z}
          ( \Lambda( u^j E_{ii}) / j! )
            \frac{ (\eta_j (x; -i + 1/2) - \eta_j (x; -i - 1/2) ) }
                 {\sinh (x/2)}   \\
    & \quad & + c_0 \frac{\cosh (x/2) -1}{\sinh (x/2)} 
    + \sum_{j =1}^m \frac{ c_j \eta_j (x; -1/2 ) }{\sinh (x/2)}
                                     \\
    & \stackrel{(2)}{=} &    \frac{ 
                      \sum_{j=0}^m \sum_{i \in \Z}
                         \Lambda ( u^j E_{ii} - u^j E_{i+1, i+1})
                           \eta_j (x; -i -1/2) 
                  }     {  \sinh (x/2)}  \\
    & \quad & + c_0 \frac{\cosh (x/2) -1}{\sinh (x/2)} 
    + \sum_{j =1}^m \frac{ c_j \eta_j (x; -1/2 ) }{\sinh (x/2)}
                                     \\
    & \stackrel{(3)}{=} & 
                     \frac{ 
                      \sum_{j=0}^m \sum_{i \in \Z_{+} } 
                       {}^b\tilde{h}_i^{(j)} \eta_j (x;  -i -1/2)
                           }   
            { \sinh (x/2)}  
   -  \frac{c_0 }{\sinh (x/2)}.
 \end{eqnarray*}
 The identity (1) is obtained from a similar identity in the proof of
 Proposition \ref{prop_deltagl} by putting $s = -\hf$.
 The identity (2) is obtained by shifting the index $i$ to $i+1$
 in the first half of the first summation of the left hand side of (2).
 The identity (3) is obtained by splitting the summation into two,
 $ \sum_{i \in \Z} = \sum_{i \in \Z_{+}} + \sum_{i \in - \Bbb N}$,
 changing the index $i$ to $-i -1$ in the second summation,
 and then using formulas(\ref{eq_bweight}) and (\ref{eq_property}).
\end{demo}

Take an irreducible quasifinite highest weight $\hB$-module $V$
with central charge $c$ and 
$$
   \Delta(x)=\frac{ F (x) }{2 \sinh (x/2 )}
$$
where $F(x)$ is an
even quasipolynomial such that $F(0) = 0$
(see Theorem \ref{th_delta}). We may write
\begin{equation}
  F(x) + c \cosh (x/2 )
 = \sum_{s \in \Bbb C} \sum_{j=0}^{m_s} a_{s,j} \eta_j (x; s -1/2)
       \label{eq_quasipoly}
\end{equation}
where $a_{s,j} \in \Bbb C$ and $a_{s,j} \neq 0$ for 
only finitely many $s \in \C$. We can fix the ambiguities in 
expressing $F(x)$ in the form (\ref{eq_quasipoly}) caused by
the symmetries (\ref{eq_property}) 
by the following rules in choosing the parameter $s$:
when $s \in \Z$, we require $s \leq 0;$ 
when $s \in \hz$, we require $s \leq \hf;$ 
when $s \notin \hf \Z$, we require that $\mbox{Im } s > 0$ 
if $\mbox{Im } s \neq 0$, or $s - [s] < \hf$ if $ s \in \Bbb R$, where
$\mbox{Im } s$ denotes the imaginary part of $s$
and $[s]$ the closest integer to $s$ which is not larger than $s$.

Decompose the set 
$\{ s \in \Bbb C \mid a_{s,j}\neq 0 \mbox{ for some $j$} \}$ 
into a disjoint union of equivalence classes 
under the equivalence condition: $s \sim s' $ if and only if
$s = \pm s' \; (\bmod{\Bbb Z})$. 
Pick a representative $s$ in an equivalence class $S$ such that
$ s =0$ if the equivalence class lies in $\Z$ and 
$ s = \hf$ if the equivalence class lies in $\hz$. 
Let $S=\{ s, s -k_1, s -k_2,\ldots \}$ 
be such an equivalence class and let $m = \max \limits_{s \in S} m_s$. 
Put $k_0 =0$.
It is easy to see by the rules in picking the parameter $s$ that
if $s = 0$ or $\hf$, then $k_1, k_2, \cdots \in \Bbb N.$

We associate to $S$ the $\gm$-module $ L_s^{[m]} (\lambda_s)$ 
in one of the following three ways:

First, if $s \notin \Bbb Z/2$, let 
$ {}^a h_{k_r}^{(j)} = a_{s + k_r, j} \; 
  (j =0, \ldots, m_s, r =0, 1, 2, \ldots)$. 
We associate to $S$ a $\hglm$-module $L_s^{[m]}(\lambda_S)$ 
with central charges and labels
$$ 
 c_j = \sum_{k_r} {}^a h_{k_r}^{(j)},
  \quad
 {}^a \lambda_i^{(j)}=\sum _{k_r \geq i} {}^a \tilde{h}_{k_r}^{(j)}
$$
where ${}^a \tilde{h}_k^{(j)} = {}^a h_k^{(j)} - c_j \delta_{k,0}$.

Secondly, if $s =0$, let 
$ {}^b \tilde{h}_{ k_r}^{(j)} = a_{k_r, j} (j =0, \ldots, m_0, r \in \Z_{+})$. 
We associate to $S$ a $\binf^{[m_0]}$-module $L_0^{[m_0]}(\lambda_S)$ 
with central charges and labels
\begin{eqnarray*}
 c_j =  \sum_{k_r} {}^b \tilde{ h}_{k_r}^{(j)} \; (j \mbox{ even}), &&
   c_j = 0 \; (j \mbox{ odd})   \\
 {}^b \lambda_0^{(j)}=\sum_{k_r \geq 0}h_{k_r}^{(j)}  \; (j \mbox{ odd}), &&
 {}^b \lambda_i^{(j)}=\sum_{k_r \geq i}h_{k_r}^{(j)}
\end{eqnarray*}
where $  i \in \Bbb N, j = 0, \ldots, m_0.$

Thirdly, if $s = \frac 12$, let 
$ {}^c h_{k_r}^{(j)} = a_{1/2 +k_r, j}
  (j = 0, \ldots, m_{1/2}, r  \in \Z_{+})$. 
We associate to $S$ the $\cinf^{[m_{1/2}]}$-module 
$L_{\frac 12}^{[m_{1/2}]}(\lambda_S)$ 
with central charges and labels
$$ 
 c_j =  \sum_{k_r} {}^c h_{k_r}^{(j)} \; (j \mbox{ even}),
   \quad c_j = 0 \; (j \mbox{ odd}),  \quad
 {}^c \lambda_i^{(j)}=\sum_{k_r \geq i} {}^c h_{k_r}^{(j)}
$$
where $ i \in \Bbb N, j = 0, \ldots, m_{1/2}$.

Denote by $\{ s_1, s_2, \ldots, s_N \}$ a set of 
representatives of equivalence classes in the set
$\{ s \in \Bbb C \mid a_{s,j}\neq 0 \mbox{ for some $j$} \}$.
By Theorem \ref{th_pullback}, the 
$\hB$-module $L_{\vec s}^{[\vec m], -}(\vec \lambda)$ 
is irreducible for $\vec s = (s_1, s_2, \ldots, s_N)$
satisfying the $(\star -)$ condition. Then  we have
$$ \Delta^-_{\vec m, \vec s, \vec \lambda} (x)
   = \sum_i \Delta^-_{m_i, s_i, \lambda(i)} (x),
      \quad c=\sum_i c_0(i).
$$

Summarizing the above, together with Theorem \ref{th_pullback},
Propositions \ref{prop_deltagl}, \ref{prop_deltab}
and \ref{prop_deltac}, we have proved the following theorem for the
$-$ case. Similarly one can prove the $+$ case with the help
of Theorem \ref{th_pullback},
Propositions \ref{prop_dtagl}, \ref{prop_deltad} and \ref{prop_dtab}.
\begin{theorem}
   Let $V$ be an irreducible quasifinite highest weight $\hBpm$-module
 with central charge $c$ and 
 $$
   \Delta^{\pm} (x) = \frac{ F(x)}{2 \sinh (x/2 )}
 $$
 for some even quasipolynomial $F(x)$ which is written in
 the form (\ref{eq_quasipoly}). Then $V$ is isomorphic to
 the tensor product of the modules $L_s^{[m], \pm}(\Lambda_S)$
 for different equivalence classes $S$.
    \label{th_realization}
\end{theorem}
\begin{remark} \rm
 A different choice of representative $s \notin \hf \Z$ in an equivalence 
 class $S$ has the effect of shifting $\hglm$ via the automorphism
 $\nu^i$ for some $i$. It is not difficult to see that any irreducible 
 quasifinite highest weight module $L(\hBpm; \xi)$ can be obtained  
 as above in an essentially unique way up to the shift.
 Note that we have always put 
 $c_j = 0 \; (j \mbox{ odd})$ when $s =0$ or $\hf$ in the $-$ case
 (resp. $- \hf$ in the $+$ case).
 We could have defined $\binfm$, $\binfmtwo$, $\cinfm$ and $\dinfm$ to be the
 central extensions of $\Binfm$, $\Binfmtwo$, $\Cinfm$ and $\Dinfm$ by
 $R_m ' = \bigoplus_{ j =1}^{[m/2]} \C u^{2j} \subset R_m$.
 We have made our choice for the convenience of notations.
\end{remark}
\section{Unitary QHWM's of $\hBpm$}
  \label{sec_unit}
  Consider the following
anti-linear anti-involution $\omega$ of the Lie algebra $\hD$:
$$ \omega ( t^k f(D) ) = t^{-k}\overline{f}(D -k), 
  \quad \omega( C) = C,
$$
where $\overline{f}$ denotes the complex conjugate.
Note that $\hBpm$ is $\omega$-invariant:
\begin{eqnarray*}
 \omega ( t^k f( D +k/2)) & = & t^{-k} \overline{f}(D- k/2)   \\
 \omega ( t^k f( D +(k + 1)/2)) & = & t^{-k} \overline{f}(D + (- k +1)/2).
\end{eqnarray*}
The goal of this section is 
to classify all unitary quasifinite highest weight modules 
over $\hBpm$ with respect to the anti-involution $\omega$.
\begin{remark} \rm
 \begin{enumerate}
 \item[1)] The anti-involution $\omega$ on $\D$ is
 compatible with the standard anti-involution $\omega '$ of $\gl$ 
 (given by the 
 complex conjugate transpose of a matrix)
 under the homomorphism $ \phi_s = \phi_s^{[0]}: \BO \longrightarrow \gl $:
 $$
  \omega '(\phi_s ( t^k f(D) ) ) = \phi_{\bar{s} }(\omega ( t^k f(D))),
 $$
 where $s \in \C$ and $\bar{s}$ denotes its complex conjugate.

 \item[2)] With respect to the anti-involution $\omega '$ of $\hgl$ defined by 
 $\omega '(A)={}^t\overline{A}$ and $\omega '(c) = \overline{c}$,
 a highest weight $\hgl$-module 
 with highest weight $\Lambda$ and central charge $c$ is unitary 
 if and only if all ${}^a h_i$ are non-negative integers
 and $c= \sum_i {}^ah_i$. This follows from the unitarity of 
 finite dimensional modules over ${\frak gl}_n$.

 Similarly, a highest weight module of $\binf$ or
 $\binftwo$ with highest weight $\Lambda$ 
 and central charge $c$ is unitary with respect to $\omega '$
 if and only if the numbers ${}^b h_i \; (i \in \Bbb Z_+)$
 are non-negative integers and $c = \hf {}^b h_0 + \sum_{i \geq 1} {}^b h_i.$

 A highest weight $\cinf$-module with highest weight $\Lambda$
 and central charge $c$ is unitary if and only if the 
 numbers ${}^c h_i \; (i \in \Bbb Z_+)$ are non-negative integers 
 and $c=\sum_{i \geq 0} {}^c h_i.$

 A highest weight $\dinf$-module with highest weight $\Lambda$
 and central charge $c$ is unitary if and only if the 
 numbers ${}^d h_i \; (i \in \Bbb Z_+)$ are non-negative integers 
 and $c = ({}^d h_0 + {}^d h_1)/2 + \sum_{i \geq 2} {}^d h_i.$
  \end{enumerate}
    \label{rem_compat}
\end{remark}
\begin{remark}
   A highest weight $\dinf$-module with highest weight $\Lambda$
 and central charge $c$ is unitary if and only if the 
 numbers ${}^d h_i \in \Z_{+} \; (i \in \Bbb N)$, $c \in \hf \Z_{+}$ 
 and $c \geq {}^d h_1 /2 + \sum_{i \geq 2} {}^d h_i.$
  \label{rem_subtle}
\end{remark}
\begin{lemma}
 \begin{enumerate}
  \item[1)] Let $V$ be a unitary quasifinite highest weight module
 over $\hB$, and let $b(w)$ be its characteristic polynomial.
 Then $b(w)$ has only real roots and any non-zero root
 of $ b (w)$ is simple; 0 is a double root of $ b (w)$ if it is a root.
 \item[2)] Let $V$ be a unitary quasifinite highest weight module
 over $\hBp$, and let $b(w)$ be its characteristic polynomial.
 Then $b(w)$ has only real simple roots. 
 \end{enumerate}
   \label{lem_unitary}
\end{lemma}
\begin{demo}{Proof}
  We will prove part 1). The proof of part 2) is similar.

   Let $v_\Lambda$ be a highest weight vector of $V$. Since $b( w)$ is
 an even polynomial, we may assume that
 $\deg b(w) = 2n$, and $ b( w) = g( w^2)$ for some polynomial $ g(w)$.
 Then the first graded subspace $V_1$ of $V$ has a basis
 $$
  \{ t^{-1}( D - 1/2)^{2l} v_\Lambda,
  \quad 0\leq l< n \}.
 $$
 Consider the action of
 $S = - \frac{1}{12} ( 4 D^3 + 4 \Delta_3 +1) \in \hB$ 
 on $\End V_1$. It is straightforward to check that
 $$
  S^j \left( t^{-1} v_\Lambda
      \right) = \left( t^{-1} ( D -1/2)^{2j}
                           \right)  v_\Lambda.
 $$
 It follows that $ g(S) ( t^{-1} v_\Lambda) = 0$ and
 $ \{ S^j ( t^{-1} v_\Lambda ), \quad 0\leq j <n \} $
 is a basis of $V_1$.
 We conclude from the above  that $g(w)$ is the characteristic 
 polynomial of the operator $S$ on $V_1$. Since the operator $S$
 is self-adjoint, all the  roots of $g(w)$ are real.

 Similarly, let $b_2(w)$ be the polynomial of minimal degree such that 
 $$ t^{-2} b_2 ( D -1) v_\Lambda = 0. $$
 Since $b_2 ( w)$
 is an odd polynomial, we may assume $b_2 ( w) = w g_2 ( w^2)$
 for some polynomial $g_2 ( w)$ of degree $k$. Consider
 the action of $T = - \frac{1}{6} (D^3 + \Delta_3 +2) \in \hB$
 on the second graded subspace $V_2$ of $V$. One can check that
 $$
  T^j \left( \left( t^{-2} ( D -1)
             \right)  v_\Lambda
      \right) = \left( t^{-2} ( D - 1)^{2j +1}
                           \right)  v_\Lambda.
 $$
 It follows that $ g_2 (T) ( (t^{-2} ( D -1)) v_\Lambda) = 0$ and
 $ \{ T^j  (t^{-2}( D -1)) v_\Lambda ), \quad 0\leq j <n \} $
 is a basis of $V_2$.

 We then conclude that  $g_2 ( w)$ is the characteristic 
 polynomial of the operator $T$ on $V_2$. Since the operator $T$
 is self-adjoint, all the  roots of $g_2 ( w)$ are real.
 By Lemma \ref{lem_charact}, $b(w)$ divides $ w b_2 ( w - 1/2)$.
 Given a root $a$ of $b(w)$, it follows that
 $$ a (a -1/2) g_2 \left( (a - 1/2)^2
                   \right) = 0. $$
 So either $ a$ is equal to $ 0, \hf,$ or $(a - 1/2)^2$ is real.
 But we know already that $a^2$ is real as well
 since $a^2$ is a root of $g( w)$.
 Thus $a $ is real.

 Now let $r$ be a root of $g( w)$ of multiplicity $m$,
 then we may write $ g( w) = h( w) (w -r)^m$ for some polynomial $ h(w)$.
 Denote $ v = (S -r)^{m -1} h( S)(t^{ -1}v_\lambda ) $.
 It is a non-zero vector in $V_1$. However, on the other hand,
 $$ 
  ( v, v) = ( c( S) (t^{ -1} v_\lambda ),
     (S-r)^{2m -2}c( S)(t^{ -1} v_\lambda) ) =0 \quad \mbox{ if } m\geq 2.
 $$
 Hence the unitarity condition forces $m =1$. Consequently, any non-zero
 root of $b(w) = g(w^2)$ has multiplicity 1 and $0$ 
 is a double root if it is a root.
\end{demo}
\begin{lemma}
   If $L_{\vec s}^{[\vec m], \pm}(\vec \Lambda)$ is a unitary quasifinite
 irreducible highest weight module over $\hBpm$, then $m =0$.
    \label{lem_onlyif}
\end{lemma}
\begin{demo}{Proof}
   By Theorem \ref{th_realization},
 Corollary \ref{cor_diff} and Lemma \ref{lem_unitary}, we see that 
 $ a_{s,j} = 0 $ $(j \geq 1 )$ in (\ref{eq_quasipoly})
 in the $-$ case. The argument for the $+$ case is the same.
\end{demo}
\begin{theorem}
  \begin{enumerate}
 \item[1)] An irreducible quasifinite highest weight $\hB$-module
 with central charge $c$
 is unitary if and only if it is positive primitive with real exponents $e_i$,
 or equivalently if and only if
 \begin{equation}
    \Delta(x) =
      \frac { \sum_i n_i ( \cosh (e_i x) - \cosh (x/2) ) }{2 \sinh (x/2)}
  \label{eq_simple}
 \end{equation}
 where $e_i$ are distinct real numbers different from $\pm \hf$,
 $n_i \in \Bbb N$, $c \in \hf \Z_{+}$ and $c \geq \sum_i n_i.$
 \item[2)] A quasifinite irreducible highest weight $\hBp$-module 
 with central charge $c$
 is unitary if and only if it is positive primitive with real exponents $e_i$
 or equivalently if and only if
 \begin{eqnarray*}
    \Delta(x) =   \frac {\sum_i n_i (\cosh (e_i x) - 1 )}{2 \sinh (x/2)}
 \end{eqnarray*}
 where $e_i$ are distinct non-zero real numbers, $n_i \in \Bbb N$,
 $c \in \hf \Z_{+}$ and $c \geq \sum_{i \neq i_0} n_i + n_{i_0}/2$. 
 Here ${i_0}$ is the index such that $e_{i_0} = 1$.
 \item[3)] Any unitary quasifinite highest weight $\hBpm$-module
 can be obtained
 by taking a tensor product of $N$ unitary irreducible highest
 weight modules over $\frak g^{[\vec 0]}$ and restricting to $\hBpm$ 
 via $\widehat{\phi}_{\vec s}^{[\vec 0]}$, where
 $\vec s=(s_1, s_2, \ldots, s_N) \; ( s_i \in \Bbb R )$
 satisfies the $(\star \pm)$ condition. 
   \end{enumerate}
   \label{th_unit}
\end{theorem}
\begin{demo}{Proof}
   The ``only if'' part of 1) was implied by the Lemma \ref{lem_onlyif}.
 The part 3) follows from Remark \ref{rem_compat}. Now the ``if'' part 
 of 1) follows from part 3) since by Theorem \ref{th_realization}
 we have realized all irreducible unitary
 quasifinite highest weight $\hB$-modules with $ \Delta(x)$
 of the form (\ref{eq_simple}). The proof of part 2) is similar
 to part 1) (cf. Remark \ref{rem_subtle}).
\end{demo}
     
By Theorem \ref{th_realization}, a positive primitive $\hBpm$-module is
the pullback of a tensor product of $N$ unitary irreducible highest
 weight modules $V_i$ over $\frak g^{[0]}$ 
 via a homomorphism $\widehat{\phi}_{\vec s}^{[0]}$, where
 $\vec s=(s_1, s_2, \ldots, s_N) \; ( s_i \in \C )$
 satisfies the $(\star \pm)$ condition. Since the $\Z$-gradations
 of $\hBpm$ and $\frak g^{[0]}$ are compatible under the
 homomorphism $\widehat{\phi}_{\vec s}^{[0]}$, 
 the $q$-character formula of a positive primitive module 
 of $\hBpm$ is given by the product of the 
 $q$-character formulas of $V_i$, which are in turn given by
 formulas (\ref{eq_Aqchar}), (\ref{eq_Bqchar}), and (\ref{eq_Cqchar})
 in $-$ case (resp. (\ref{eq_Dqchar}) in $+$ case).
\begin{remark} \rm
    One can introduce a tensor category ${\cal O}_{b,\pm}$
 consisting of all positive
 primitive modules of $\hBpm$ with non-negative half-integral 
 central charge as was done for $\W$ in \cite{W1} (also cf. \cite{W2})
 equipped with the usual tensor product. 
 One can show that ${\cal O}_{b,\pm}$ is a semisimple tensor category.
 A tensor  product of two irreducible modules in such a category
 is decomposed into an (infinite) sum of irreducibles with
 finite multiplicities.  Certain reciprocity laws 
 can be established between these multiplicities 
 and some branching coefficients in finite dimensional Lie groups
 as a formal consequence of various duality results
 we will establish in the next sections.
\end{remark}
\section{FFR's of QHWM's over $\hB$ with 
$c \in - \Bbb N /2$}
   \label{sec_boson}
\subsection{Dual pair $\Wpairdc$}
 In this subsection we study the free field realizations (FFR's) of certain
primitive $\hB$-modules with negative 
integral central charges in some bosonic Fock spaces and establish
certain duality results. These duality results are intimately 
related to duality results obtained in \cite{W2}. We refer the reader
to \cite{W2} for more detail. In the
case of $\hD$ similar results were obtained in \cite{KR2}.

Let us take a pair of bosonic ghost fields
$$
 \gamma^+ (z) = \sum_{ n \in \hz }
                         \gamma^+_n z^{-n - \hf}, 
  \quad \gamma^- (z) = \sum_{ n \in \hz }
                         \gamma^-_n z^{-n - \hf}
$$
with 
             %
%
the following commutation relations
$$ [ \gamma^+_m, \gamma^-_n ]
= \delta_{m+n, 0}.$$
We define the Fock space ${\cal F}^{\bigotimes -1}$ 
of the fields $\gamma^+ (z)$ and $\gamma^-(z)$, 
generated by the vacuum $\vac$, satisfying
\begin{eqnarray*}
  \gamma^{+}_n \vac = \gamma^{-}_n \vac = 0 ,  \quad  n \in \hn.
\end{eqnarray*}
Now we take $l$ pairs of bosonic ghost fields 
$\gamma^{+,p} (z), \gamma^{-,p} (z) \;(p = 1, \dots, l)$
and consider the
corresponding Fock space $\Fminusl$.

It is well known and straightforward to verify that
\begin{eqnarray}
  E(z,w) & \equiv & \sum_{i,j \in \Z} E_{ij} z^{i-1}w^{-j } 
  = - \sum_{p =1}^l :\gamma^{+,p}  (z) \gamma^{-,p} (w):
\end{eqnarray}
defines an action of $\hgl$ on $\Fminusl$ with central charge $-l$.

It is known \cite{FF} that
the Fourier components of the following generating functions
\begin{eqnarray}
  e^{pq}_{**} (z) & \equiv & \sum_{i,j \in \Z} e^{pq}_{**}(n)  
                   z^{ -n -1}
  = : \gamma^{+,p} (z) \gamma^{+,q} ( -z): 
    \quad (p \neq q )          \label{eq_plusgamma}     \\ 
  e^{pq} (z) & \equiv & \sum_{i,j \in \Z} e^{pq}(n)
                   z^{ -n -1}
  = : \gamma^{-,p}(z) \gamma^{-,q} ( -z):
    \quad (p \neq q )      \label{eq_minusgamma}       \\
  e_*^{pq} (z) & \equiv & \sum_{i,j \in \Z} e_*^{pq} (n)
                   z^{ -n -1}
  = : \gamma^{+,p} (z) \gamma^{-,q} (z):
   \quad (p,q = 1, \dots, l )  \label{eq_mixgamma}   
\end{eqnarray}
span an affine algebra ${\frak {gl}}^{(2)} (2l)$
of type $A^{(2)}_{2l -1}$ of central charge $ -1$. 
The horizontal subalgebra of the affine algebra 
${\frak {gl}}^{(2)} (2l)$ spanned by
$e^{pq}_{**} (0), e_*^{pq} (0)$, $e^{pq} (0)$ $(p, q = 1, \ldots, l)$ 
is isomorphic to the Lie algebra ${ \frak {so}} (2l)$.
The Borel subalgebra of ${ \frak {so}} (2l)$ is taken to
be the one spanned by 
$e_{*}^{pq} ( p \leq q), e_{**}^{pq}, p, q = 1, \ldots, l$
and the Cartan subalgebra is spanned by
$e_{*}^{pp},  p  = 1, \ldots, l.$
The action of ${ \frak {so}} (2l)$ can be lifted to
an action of $SO(2l)$ on $\Fminusl$ and then extends to $O(2l)$ naturally. 
For example the operator which commutes with 
$\psi^{\pm,k} (z)$ $(k =1, \ldots, l-1)$ and
sends $\psi^{\pm,l} (z)$ to $\psi^{\mp, l} (z)$
lies in $ \in O(2l) - SO(2l)$.
Put 
\begin{eqnarray}
 \begin{array}{rcl}
  \lefteqn{
    \sum_{i,j \in \Z} (E_{ij} - ( -1)^{i+j} E_{1-j, 1-i})z^{i-1}w^{-j}
          }   \\
   & = & \sum_{ k=1}^l \left(\gamma^{+,k} (z) \gamma^{-,k} (w)
                        - \gamma^{+,k} ( -w) \gamma^{-,k} ( -z)
                  \right).
 \end{array}
   \label{eq_Fminuslcinf}
\end{eqnarray}
The operators $E_{ij} - ( -1)^{i+j} E_{1-j, 1-i} \;( i,j \in \Z )$
span $\cinf$ with central charge $-l$. 
It is known \cite{W2} that the actions of $O(2l)$ and $\cinf$
commute with each other and form a dual pair on $\Fminusl$
in the sense of Howe \cite{H1, H2}. 

$O(2l)$ is a semi-direct product of $SO(2l)$ by $\Z_2$.
 If $\lambda$ is a representation of $SO(2l)$ of
highest weight $(m_1, m_2, \ldots, m_l )$ ($m_l \neq 0$),
then the induced representation of $(m_1, m_2, \ldots, m_l )$ 
to $O(2l)$ is irreducible and its restriction to
$SO(2l)$ is a sum of $(m_1, m_2, \ldots, m_l )$ and
$(m_1, m_2, \ldots, - m_l )$.
We denote this irreducible representation $\lambda$ of $O(2l)$ by
$(m_1, m_2, \ldots, \overline{m}_l )$, where $m_l$ is chosen to be 
greater than 0. If $ m_l = 0$, the representation
$\lambda = (m_1, m_2, \ldots, m_{l-1}, 0 )$ extends to two
different representations of $O(2l)$, denoted by
$\lambda$ and $\lambda \bigotimes det$,
where $det$ is the $1$-dimensional non-trivial representation of
$O(2l)$. We denote 
\begin{eqnarray*}
  \Sigma(D) =
   & & \left\{ (m_1, m_2, \ldots, \overline{m}_l ) \mid
       m_1 \geq m_2 \geq \ldots \geq m_l > 0, m_i \in \Z;  \right.
                               \\
   & &   (m_1, m_2, \ldots, m_{l-1}, 0 ) \bigotimes det, \\
   & &   \left.   (m_1, m_2, \ldots, m_{l-1}, 0 ) \mid
       m_1 \geq m_2 \geq \ldots \geq  m_{l-1} \geq 0, m_i \in \Z \right\}. 
\end{eqnarray*}
 
Define a map $\Lambda^{ \frak{dc}}$
from $\Sigma(D)$ to ${\cinf}_0^*$
by sending $ \lambda = (m_1, \cdots, \overline{ m}_l) \; (m_l > 0 )$ 
to 
$$
 \Lambda^{ \frak{dc}} (\lambda) =
    (-l -m_1) {}^c \hL_0 + \sum_{k =1}^l (m_k -m_{k+1}) {}^c \hL_k
                      \quad (m_{l+1} =0),
$$
sending $(m_1, \cdots, m_j, 0, \cdots, 0 ) \; (j < l )$ to 
$$ \Lambda^{ \frak{dc}} (\lambda) =
    (-l -m_1) {}^c \hL_0 + \sum_{k =1}^j (m_k -m_{k+1}) {}^c \hL_k, $$
and sending $ (m_1, \cdots, m_j, 0, \cdots, 0 ) \bigotimes {det} $ to 
 $$  \Lambda^{ \frak{dc}} (\lambda) =
   (-l -m_1) {}^c \hL_0 + \sum_{k =1}^{j-1} (m_k -m_{k+1}) {}^c \hL_k
     + (m_j -1) {}^c \hL_j + {}^c \hL_{2l-j}, $$
if $m_1 \geq \ldots \geq m_{j} > m_{j+1} = \ldots = m_l = 0.$

The following theorem was proved in \cite{W2}.
\begin{theorem}
 We have the following $\Mpairdc$-module decomposition:
 \begin{equation}
   \Fminusl = \bigoplus_{\lambda \in \Sigma (D)} 
          V( O(2l); \lambda) \otimes L \left(\cinf;
                                 \Lambda^{ \frak{dc}} (\lambda), -l
                               \right)
 \end{equation}
 where $V( O(2l); \lambda)$ is the irreducible $O (2l)$-module
 parametrized by $\lambda \in \Sigma (D)$ and 
 $L \left(\cinf; \Lambda^{ \frak{dc}} (\lambda), -l \right)$
 is the irreducible $\cinf$-module of highest weight
 $\Lambda^{ \frak{dc}} (\lambda) $ and central charge $ -l$.
  \label{th_Mpairdc}
\end{theorem}

Note that $\hB$ acts on $\Fminusl$ via the composition of
the homomorphism 
$\widehat{\phi}_{1/2}$ and the action of 
$\cinf$ given by (\ref{eq_Fminuslcinf}).
Introduce the following generating function
$$ T^n_{1/2} (z) = \sum_{k \in \Z} T_k^{n, 1/2} z^{ -k-n-1}
  \quad ( n \in 2 \Z_{+} +1) 
$$ 
for the basis $T^{n, 1/2}_k $ defined in (\ref{eq_basis}) of $\hB$.
It follows from Proposition \ref{prop_ext} that 
\begin{eqnarray}
  c & = & -l   \\
\widehat{\phi}_{1/2}
   \left( e^{xD} - e^{-xD} \right)
      & = & \phi_{1/2} \left( e^{xD} - e^{-xD} \right) -l \tanh(x/4).
  \label{eq_diff}
\end{eqnarray}
It follows that 
\begin{equation}
 \widehat{\phi}_{1/2}(T_0^{n, 1/2}) = \phi_{1/2}(T_0^{n, 1/2}) - \alpha_n l
  \label{eq_half}
\end{equation}
for some constant $\alpha_n$ determined by (\ref{eq_diff})
depending on $n$ only.

\begin{lemma}
 One can recast the action of $\hB$ in terms of $ T_{ 1/2}^n (z)$ as
 \begin{eqnarray*}
  T_{ 1/2}^n (z) = - \sum_{k =1}^l
           \left( : \gamma^{+,k} (z) \partial_z^n \gamma^{-,k} (z):
             +  : \partial_z^n \gamma^{+,k} ( -z) \gamma^{-,k} ( -z):
           \right) - \alpha_n l z^{ -n -1}.
 \end{eqnarray*}
  \label{lem_dc}
\end{lemma}
\begin{demo}{Proof}
 We have the following calculations.
 \begin{eqnarray*}
  T_{1/2}^n (z) & = & \sum_{k \in \Z} T_k^{n, 1/2} z^{ -k-n-1} \\
   & \stackrel{(1)}{=} & - \sum_{k \in \Z}
     ( [ -j]_n + ( -1)^{k +1} [ j - k -1]_n ) E_{j -k, j} z^{ -k-n-1}   
      - \alpha_n l z^{ -n -1}    \\
   & \stackrel{(2)}{=} &  - \sum_{k, j \in \Z}
    [ -j]_n  \left( E_{j -k, j} + ( -1)^{k +1}E_{ 1-j, k-j +1} \right)
      z^{ -k-n-1} - \alpha_n l z^{ -n -1}    \\
   & \stackrel{(3)}{=} &   - \sum_{k =1}^l
           \left( : \gamma^{+,k} (z) \partial_z^n \gamma^{-,k} (z):
             +  : \partial_z^n \gamma^{+,k} ( -z) \gamma^{-,k} ( -z):
           \right) - \alpha_n l z^{ -n -1},
 \end{eqnarray*}
 where equation (1) is given by the homomorphism 
 $\widehat{\phi}_{1/2}$ (cf. (\ref{eq_emb}) and (\ref{eq_half})),
 (2) is obtained by change of variables in the second summation
 (replacing $j - k -1$ with $j$), and (3) follows by taking the $n$-th
 derivative of (\ref{eq_Fminuslcinf}) with respect to $w$ and then
 putting $ w = z$.
\end{demo}

For $ \lambda = (m_1, \cdots, \overline{ m}_l)\in \Sigma(D) \; (m_l > 0 ) $ 
we let $e (\lambda)$ be the set of the exponents $k$
with multiplicity $m_k - m_{k +1}$ $(k =1, \ldots, l)$
and the exponent $0$ with multiplicity $ -l -m_1$, where $m_{l +1} = 0$;
for $\lambda = (m_1, \cdots, m_j, 0, \cdots, 0 )\in \Sigma(D)$ 
$ ( j < l) $, we let $e (\lambda)$ be the set of the exponents $k$
with multiplicity $m_k - m_{k +1}$ $(k =1, \ldots, j)$
and the exponent $0$ with multiplicity $ -l -m_1$, where $m_{j +1} = 0$;
for $ (m_1, \cdots, m_j, 0, \cdots, 0 ) \bigotimes {det} \in \Sigma(D)$ 
$ ( j < l) $
we let $e (\lambda)$ be the set of the exponents $k$
with multiplicity $m_k - m_{k +1}$ $(k =1, \ldots, j -1)$,
the exponent $j$ with multiplicity $m_j -1$,
the exponent $2l -j$ with multiplicity $1$
and the exponent $0$ with multiplicity $ -l -m_1$.

Now we can state a duality theorem between $O(2l)$ and $\hB$
(recall Convention~\ref{conven_1}).
\begin{theorem}
 $O(2l)$ and $\hB$ form a dual pair on $\Fminusl$. More
 explicitly we have the following $\Wpairdc$-module decomposition:
 \begin{equation}
   \Fminusl = \bigoplus_{\lambda \in \Sigma (D)} I_{\lambda}
       \equiv \bigoplus_{\lambda \in \Sigma (D)} 
          V( O(2l); \lambda) \otimes L \left(\hB;
                                 e (\lambda), -l
                               \right).
 \end{equation}
   \label{th_Wpairdc}
\end{theorem}
\begin{demo}{Proof}
  Since the actions of $O(2l)$ commutes with $\cinf$, the 
 action of $O(2l)$ commutes with the action of $\hB$ given in
 Lemma \ref{lem_dc} by Proposition \ref{prop_phihalf}. 
 So the decomposition of
 the Fock space into isotypic subspaces with respect to
 the dual pair $\Fpaircc$ can be regarded as decomposition with respect to
 the dual pair $\Wpaircc$ as well. 
 By Theorems~\ref{th_Mpairdc} and \ref{th_pullback} 
 each isotypic subspace is
 irreducible under the joint action of $Sp(2l)$ and $\hB$.
 By Proposition~\ref{prop_deltac} and Theorem~\ref{th_Mpairdc} the highest
 weight of the representation $I_{\lambda}$ with respect to $\hB$
 is given as follow.
 \begin{eqnarray*}
        \Delta_{m, s, \Lambda} (x)
                & = &  \frac{ (-l -m_1) +   \sum_{k =1}^j
                              ( m_k -m_{k+1}) \cosh (k x)
                          }   
                           { 2 \sinh (x/2)}    
              -  \frac{ -l \cosh (x/2)}{2 \sinh (x/2)}    
 \end{eqnarray*}
 for $ \lambda = (m_1, \cdots, \overline{ m}_l)  \in \Sigma(D)\; (m_l > 0 )$.

 The cases $\lambda = (m_1, \cdots, m_j, 0, \cdots, 0 )$ 
 $ ( j < l) $ and $ (m_1, \cdots, m_j, 0, \cdots, 0 ) \bigotimes {det}$
 can be treated similarly.
\end{demo}
\subsection{Dual pair $\Wpairbc$}
 In this subsection we will realize certain
primitive $\hB$-modules with negative 
half-odd-integral central charges in some Fock spaces
and establish some duality theorems. We need a bosonic field
 $ \chi (z) = \sum_{ n \in \hz} \chi_n z^{ -n - \hf}$ which 
satisfies the following commutation relations:
$$[ \chi_m , \chi_n ] = ( -1)^{m +\hf} \delta_{m, -n}, \quad m,n \in \hz .$$
Let $\Fminuslhalf$ be the tensor product of the Fock space 
of $l$ pairs of bosonic ghost fields
$\gamma^{\pm, k} (z) \; ( k = 1, \ldots, l)$
and the Fock space $\Fminushalf$ of $\chi (z)$.

It is known \cite{FF} that
the Fourier components of the generating functions 
(\ref{eq_plusgamma}), (\ref{eq_minusgamma}), (\ref{eq_mixgamma})
and the following generating functions
\begin{eqnarray*}
 \zeta (z) \equiv \sum_{n \in \Z}
            \zeta(n) \NN & = & : \chi (z) \chi ( -z):  \nonumber     \\
 e^{p} (z) \equiv \sum_{n \in \Z}
            e^{p}(n) \NN & = & :
              \gamma^{-,p}(z) \chi ( -z):  \nonumber     \\
 e_*^{p} (z) \equiv \sum_{n \in \Z}
           e_*^{p}(n) \NN & = & : \gamma^{+,p}(z) \chi (z):  \nonumber    \\
\end{eqnarray*}
span an affine algebra $A^{(2)}_{2l}$ of central charge $-1$
when acting on $\Fminushalf$. 
Its horizontal subalgebra is isomorphic to ${\frak {so}}(2l +1)$.
We take the Borel subalgebra 
$\frak{b}({\frak {so}}(2l +1))$ to be the one spanned by
$e_{*}^{pq}(0) (p \leq q), e_*^{p}(0), e_{**}^{pq}(0), p, q = 1, \ldots, l$.
The Cartan subalgebra $\frak{h}({\frak {so}}(2l +1))$
is spanned by
$e_{*}^{pp}(0) , p = 1, \ldots, l$.

The action of ${\frak {so}}(2l +1)$ can be lifted to $ SO( 2l +1)$
on $\Fminuslhalf$ and then extends naturally to $ O( 2l +1)$. 
For example the operator which commutes with 
$\Psi^{\pm,k} (z), k = 1, \ldots, l$ and sends $\chi (z)$ to $- \chi (z)$
lies in $O(2l +1) - SO( 2l +1)$. Let 
\begin{eqnarray}
 \begin{array}{rcl}
  \lefteqn{
   \sum_{i,j \in \Z} (E_{ij} - ( -1)^{i +j} E_{1-j, 1-i})z^{ i-1} w^{-j}
          }              \\
  & = & \sum_{ k=1}^l \left( :\gamma^{+,k} (z) \gamma^{-,k} (w):
                             + :\gamma^{+,k} ( -w) \gamma^{-,k} ( -z):
                      \right)
              + : \chi (z) \chi (-w):.
  \end{array}
 \label{eq_Fminusodd}
\end{eqnarray}
The operators $ E_{ij} - ( -1)^{i +j}E_{1-j, 1-i} \; ( i,j \in \Z )$
span $\cinf$ with central charge $ -l -1/2$. It is shown \cite{W2}
that $ O( 2l +1)$ and $\cinf$ form a dual pair on $\Fminuslhalf$.

An irreducible modules of $SO(2l+1)$ is
parametrized by its highest weight
$\lambda = (m_1, \ldots, m_l ), m_1 \geq \ldots \geq m_l \geq 0, m_i \in \Z.$
It is well known that
$ O (2l+1)$ is isomorphic to the direct product $ SO (2l+1) \times \Z_2$
by sending the minus identity matrix to $-1 \in \Z_2 = \{ \pm 1 \}.$
Denote by $det$ the non-trivial one-dimensional representation of $O(2l+1)$.
An irreducible representation $\lambda$ of $SO(2l+1)$ extends
to two different irreducible representations of $O(2l+1)$
by tensoring with the trivial representation and non-trivial
representation of $\Z_2$, 
denoted by $\lambda$ and $\lambda \bigotimes det$.
All irreducible representations of
$O(2l+1)$ can be obtained this way.
Then we can parametrize irreducible 
representations of $O(2l+1)$ by $ \Sigma(B)$ consisting of 
highest weights
$ (m_1, \ldots, m_l)$ and $ (m_1, \ldots, m_l) \bigotimes det $.

Define a map $ \Lambda^{ \frak{bc}}$ from $\Sigma(B)$ to ${\cinf}_0^*$ 
by sending $\lambda = (m_1, m_2, \ldots, m_l)$ to
$$ \Lambda^{ \frak{bc}} (\lambda) =
     (-l -m_1 -1/2) {}^c \hL_0 + \sum_{k =1}^j (m_k -m_{k+1}) {}^c \hL_k $$
and sending $\lambda = (m_1, m_2, \ldots, m_l) \bigotimes det $ to
$$ \Lambda^{ \frak{bc}} (\lambda) =
     (-l -m_1 -1/2) {}^c \hL_0 + \sum_{k =1}^{j-1} (m_k -m_{k+1}) {}^c \hL_k
    + (m_j -1) {}^c \hL_j + {}^c \hL_{2l -j +1}, $$
where $ m_1 \geq \ldots \geq m_j > m_{j+1} = \ldots = m_l = 0 .$
The following theorem is proved in \cite{W2}.
\begin{theorem}
  We have the following $\Mpairbc$-module decomposition:
 \begin{eqnarray*}
   \Fminuslhalf = \bigoplus_{\lambda \in \Sigma (B)} 
          V( O(2l +1); \lambda) \otimes L \left(\cinf;
                                 \Lambda^{ \frak{bc}} (\lambda), -l -1/2
                               \right)
 \end{eqnarray*}
 where $V( O(2l +1); \lambda)$ is the irreducible $O(2l +1)$-module
 parametrized by $\lambda \in \Sigma (B)$ and 
 $L \left(\cinf; \Lambda^{ \frak{bc}} (\lambda), -l -1/2 \right)$
 is the irreducible highest weight $\cinf$-module of highest weight
 $\Lambda^{ \frak{bc}} (\lambda) $
 and central charge $-l -1/2$.
  \label{th_Mpairbc}
\end{theorem}

Note that $\hB$ now acts on $\Fminuslhalf$ via the composition
of the homomorphism 
$\widehat{\phi}_{1/2}$ and the action of $\cinf$ 
given by formula (\ref{eq_Fminusodd}). The following lemma
is analogous to Lemma \ref{lem_dc}.
\begin{lemma}
We can recast the action of $\hB$ in terms of generating function
$ T_{1/2}^n (z) = \sum_{k \in \Z} T_k^{n, 1/2} z^{ -k-n-1} $ as 
\begin{eqnarray*}
 T_{1/2}^n (z) & = & - \sum_{k =1}^l
           \left( : \gamma^{+,k} (z) \partial_z^n \gamma^{-,k} (z):
             +  : \partial_z^n \gamma^{+,k} ( -z) \gamma^{-,k} ( -z):
           \right)   \\
          & &  -  :\chi (z) \partial_z^n \chi (-z): 
               - \alpha_n (l + 1/2) z^{-k-n-1},  \quad n \in 2 Z_{+} +1.
\end{eqnarray*}
\end{lemma}

For $\lambda = (m_1, \cdots, m_j, 0, \cdots, 0 ) \in \Sigma(B)$, 
we let $e (\lambda)$ be the set of the exponents $k$
with multiplicity $m_k - m_{k +1}$ $(k =1, \ldots, j)$, where $m_{j +1} = 0$,
and the exponent $0$ with multiplicity $ -l -m_1 -1/2$; 
for $ \lambda = (m_1, \cdots, m_j, 0, \cdots, 0 ) 
\bigotimes {det} \in \Sigma(B)$ 
we let $e (\lambda)$ be the set of the exponents $k$
with multiplicity $m_k - m_{k +1}$ $(k =1, \ldots, j -1)$,
the exponent $j$ with multiplicity $m_j -1$,
the exponent $2l -j +1$ with multiplicity $1$,
and the exponent $0$ with multiplicity $ -l -m_1 -1/2$.

Now we have the following duality theorem on the joint 
action on $\Fminuslhalf$ of the dual pair $\Wpairbc$. The proof is
similar to that of Theorem \ref{th_Wpairdc} which is now based
on Theorem \ref{th_Mpairbc} and Proposition \ref{prop_deltac}.
\begin{theorem}
  We have the following $\Wpairbc$-module decomposition:
 \begin{eqnarray*}
   \Fminuslhalf = \bigoplus_{\lambda \in \Sigma (B)} 
          V( O(2l +1); \lambda) \otimes L \left(\hB;
                                 e (\lambda),-l -1/2
                               \right).
 \end{eqnarray*}
\end{theorem}
\section{FFR's of QHWM's over $\hB$ 
with $c \in \Bbb N$}
   \label{sec_fermion}
 In this section we will realize certain primitive $\hB$-modules with 
positive integral central charges in some fermionic Fock spaces.
Similar results for $\hD$ were obtained in \cite{FKRW} (see also \cite{KR2}).
Let us take a pair of fermionic fields
$$ 
 \psi^{+}(z) = \sum_{n \in \underline{\Z} } \psi^{+}_n z^{ -n -\hf + \epsilon},
  \quad \psi^{-}(z) = \sum_{n \in \underline{\Z} }
                           \psi^{-}_n z^{ -n -\hf + \epsilon}, 
  \quad \underline{\Z} = \hf + \Z \mbox{ or } \Z
$$
with the following anti-commutation relations
\begin{eqnarray*}
 [\psi^{+}_m, \psi^{-}_n ]_{+} & = & \delta_{m+n, 0} \\
 {[} {\psi}^{\pm}_m, {\psi}^{\pm}_n {]}_{+} & = & 0.
\end{eqnarray*}
We take the convention here and below
that $\epsilon = 0$ if $\underline{\Z} = \hz;$ and 
$\epsilon = \hf$ if $ \underline{\Z} = \Z.$
Denote by $\cal F$ the Fock space of the fields 
$\psi^{-}(z)$ and $\psi^{+}(z)$, 
generated by the vacuum $\vac$, satisfying
\begin{eqnarray*}
  \psi^{+}_n \vac = \psi^{-}_n \vac = 0 \quad ( n \in \hn), &&\mbox{ when }
   \underline{\Z} = \hz,                       \\
 \psi^{+}_n \vac = \psi^{-}_{n+1} \vac = 0 \quad (n \in \Z_{+} ),
   && \mbox{ when } \underline{\Z} = \Z. 
\end{eqnarray*}

Now we take $l$ pairs of fermionic fields, 
$\psi^{\pm,p} (z) \;( p = 1, \dots, l)$
and consider the corresponding Fock space $\Fl$.

Introduce the following ``twisted'' generating functions
\begin{eqnarray}
  E(z,w) & \equiv & \sum_{i,j \in \Z} E_{ij} z^{i-1 + 2\epsilon }
     w^{-j  } 
   = \sum_{k =1}^l :\psi^{+,k} (z) \psi^{-,k} (w):,  \\
  e^{pq} (z) & \equiv & \sum_{n \in \Z} e^{pq} (n) z^{ -n -1 +2 \epsilon}
    =  : \psi^{-,p}(z) \psi^{-,q} ( -z): 
    ,  \label{eq_minuspsi}   \\ 
  e_{**}^{pq} (z) & \equiv &
                   \sum_{n \in \Z} e_{**}^{pq}(n) z^{ -n -1 +2 \epsilon}
    =  : \psi^{+,p}(z) \psi^{+,q} ( -z):, \label{eq_pluspsi}  \\
  e_*^{pq} (z) & \equiv & \sum_{n \in \Z} e_*^{pq} (n) z^{ -n -1 +2 \epsilon}
    =  : \psi^{+,p}(z) \psi^{-,q} (z): +  \delta_{p,q} \epsilon z^{-1}, 
                            \label{eq_mixpsi} 
\end{eqnarray}
where $ p,q = 1, \ldots, l$. It is known \cite{FF} that
the Fourier components of these generating functions
span a representation of the twisted affine algebra
$ {\frak gl}^{(2)} ( 2l)$ 
of type $ A^{(2)}_{2l-1}$ with central charge $ 1$. 

\vspace{0.15in}
\noindent{\bf I. Case ${\bf \underline{\Z} =\hz}$: dual pair ${\bf \Wpaircc}$}
\vspace{0.1in}

The horizontal subalgebra of $ {\frak gl}^{(2)} ( 2l)$ 
spanned by the operators 
$e^{pq}(0), e_*^{pq}(0),$ $ e_{**}^{pq}(0)$, $ (p, q = 1, \cdots, l) $
is isomorphic to Lie algebra $\frak {sp}(2l)$.
We identify the Borel subalgebra ${\frak b}( {\frak {sp}}(2l) )$
with the one generated by $e_*^{pq}(0)\; (p \leq q), e_{**}^{pq}(0)$
$ (p, q = 1, \ldots, l )$ and the Cartan subalgebra with
the one generated by $e_*^{pp}(0)\;  (p= 1, \ldots, l )$. Let
 \begin{eqnarray}
   \begin{array}{rcl}
    \lefteqn{
      \sum_{i,j \in \Z} ( E_{i,j} -  (-1)^{i +j} E_{1-j, 1-i}) z^{i-1} w^{-j} 
            }     \\
     & = & \sum_{k =1}^l \left( : \psi^{+,k} (z) \psi^{-,k} (w):
                          + : \psi^{+,k} ( -w) \psi^{-,k} ( -z):
                  \right).
   \end{array}
  \label{eq_Flcinf}
 \end{eqnarray}
The operators $E_{i,j} -  (-1)^{i +j} E_{1-j, 1-i} \; (i, j \in \Z)$
span $\cinf$ with central charge $l$. 
The action of $\frak {sp}(2l)$ on $\Fl$ can
be integrated to $Sp(2l)$. It is known \cite{W2}
that the actions of $Sp(2l)$ and $\cinf$ commute with
each other on $\Fl$ and they form a dual pair.

Finite dimensional irreducible modules of $Sp (2l)$ are parametrized
by the highest weights in 
$$
 \Sigma (C) = \{ \lambda
 = (m_1, \ldots, m_l), m_1 \geq \ldots \geq m_l, m_i \in \Z_{+} \}.
$$
We define a map $\Lambda^{\frak{cc} }$ from $\Sigma(C)$ to ${\cinf}_0^*$ by
sending $(m_1, \ldots, m_l)$ to
$$ \Lambda^{\frak{cc} } (\lambda) =
       (l -j ) {}^c \hL_0 + \sum_{k=1}^j{}^c \hL_{m_k},$$
where $ m_1 \geq \ldots \geq m_j > m_{j +1} = \ldots = m_l = 0.$
We quote the following theorem from \cite{W2}.
\begin{theorem} 
  We have the following $\Fpaircc$-module decomposition:
 \begin{eqnarray}
  \Fl =  \bigoplus_{\lambda \in \Sigma (C) } 
          V(Sp( 2l); \lambda) \otimes L \left(\cinf;
                                 \Lambda^{ \frak{cc} }(\lambda), l
                               \right)
   \label{decom_cc}
 \end{eqnarray}
 where $V(Sp(2l); \lambda)$ is the irreducible $Sp(2l)$-module
 parametrized by $\lambda \in \Sigma (C)$, and
 $L \left(\cinf; \Lambda^{\frak{cc} } (\lambda), l \right)$ 
 is the irreducible highest weight $\cinf$-module of highest weight
 $ \Lambda^{\frak{cc}} (\lambda)$ and central charge $l$.
  \label{th_Fpaircc}
\end{theorem}

Note that $\hB$ acts on $\Fl$ via the composition
of the homomorphism $\widehat{\phi}_{1/2}$
and the action of $\cinf$ given by formula (\ref{eq_Flcinf}).
Similarly as in Lemma \ref{lem_dc} we can rewrite the action of $\hB$ 
in terms of generating function $T^n_{1/2} (z)$ as
\begin{eqnarray}
      T^n_{1/2} (z)
    & = & - \sum_{k =1}^l
           \left( : \psi^{+,k} (z) \partial_z^n \psi^{-,k} (z):
                    + : \partial_z^n \psi^{+,k} ( -z) \psi^{-,k} ( -z):
           \right)  \nonumber   \\
    &   & + \alpha_n l z^{ -k -n -1}.
  \label{vo_Fl}
\end{eqnarray}

Given $\lambda = (m_1, \ldots, m_l) \in \Sigma(C)$,
where $ m_1 \geq \ldots \geq m_j > m_{j +1} = \ldots = m_l = 0$,
let $ e (\lambda)$ be the set of exponents $ m_k$ $(k =1, \ldots, j) $
with multiplicity 1 and the exponent $0$ with multiplicity $l -j$.
\begin{theorem}
  $Sp(2l)$ and $\hB$ form a dual pair on $\Fl$. Furthermore
  we have the following $\Wpaircc$-module decomposition:
 \begin{eqnarray*}
  \Fl =  \bigoplus_{\lambda \in \Sigma (C) } I_{\lambda}
      \equiv \bigoplus_{\lambda \in \Sigma (C) } 
          V(Sp( 2l); \lambda) \otimes L \left(\hB;
                                 e (\lambda), l
                               \right).
 \end{eqnarray*}
  \label{th_Wpaircc}
\end{theorem}
\begin{demo}{Proof}
  Since the actions of $Sp(2l)$ commutes with $\cinf$, the 
 action of $Sp(2l)$ commutes with the action of $\hB$ given by
 formula (\ref{vo_Fl}) by Proposition \ref{prop_phihalf}. 
 So the decomposition of
 the Fock space into isotypic subspaces with respect to
 the dual pair $\Fpaircc$ can be regarded as decomposition with respect to
 the dual pair $\Wpaircc$ as well. By Theorem~\ref{th_pullback} 
 each isotypic subspace is
 irreducible under the joint action of $Sp(2l)$ and $\hB$. 
 By  Proposition~\ref{prop_deltac} and Theorem \ref{th_Fpaircc} the highest
 weight of the representation $I_{\lambda}$ with respect to $\hB$
 is given as follows.
 \begin{eqnarray*}
    \Delta_{m, s, \Lambda} (x)
      & = & \frac{ ( l-j) + \sum_{k=1}^j  \cosh (m_k x) } { 2 \sinh (x/2)}    
        - \frac{l \cosh (x/2)}{2 \sinh (x/2)}.
 \end{eqnarray*}
\end{demo}
\vspace{0.15in}
\noindent{\bf II. Case ${\bf \underline{\Z} =\Z}$: dual pair ${\bf \Wpairdb}$}
\vspace{0.1in}

    In this case 
the horizontal subalgebra of $ {\frak gl}^{(2)} ( 2l)$ spanned by
the operators $e^{pq}, e_*^{pq}, e_{**}^{pq}, \; (p, q = 1, \cdots, l)$
is isomorphic to Lie algebra $\frak {so}(2l)$.
In particular, the operators $e_*^{pq}\; ( p, q = 1, \cdots l)$ form
a subalgebra ${\frak {gl}}(l)$ in the horizontal ${\frak {so}}(2l)$.
We identify the Borel subalgebra ${\frak b}( {\frak {so}}(2l) )$
with the one spanned by $e_*^{pq}\; (p \leq q), e_{**}^{pq}, \;
p, q = 1, \cdots, l.$ 
The action of $\frak {so}(2l)$ can be lifted to $Pin(2l)$ on $\Fl$.
Recall that $Pin(2l)$ (resp. $Spin(2l)$)
is the double covering group of $O(2l)$ (resp. $SO(2l)$).

Denote $ {\bf 1}_l = (1, 1, \ldots, 1) \in \Z^l $ and 
${\bar{\bf 1} }_l  = (1, 1, \ldots, 1, -1) \in \Z^l $.
An irreducible representation of $Spin (2l)$ 
which does not factor to $SO (2l)$ 
 is an irreducible representation of ${\frak {so}} (2l)$
given by its highest weight 
\begin{eqnarray}
  \lambda = \hf {\bf 1}_l + (m_1, m_2, \ldots, m_l)    \label{eq_wtplus}
\end{eqnarray}
or 
\begin{eqnarray}
\lambda = \hf {\bar{\bf 1} }_l + (m_1, m_2, \ldots, -m_l)  \label{eq_wtminus}
\end{eqnarray}
where $m_1 \geq \ldots \geq m_l \geq 0, m_i \in \Z.$
A representation of $Pin (2l)$ induced from
$\lambda $ of $Spin(2l)$ of the form (\ref{eq_wtplus})
is decomposed into a sum of two irreducible $Spin(2l)$-modules
of highest weights (\ref{eq_wtplus}) and (\ref{eq_wtminus}). We use 
$\lambda = \hf |{\bf 1}_l | + (m_1, m_2, \ldots, \overline{m}_l), m_l \geq 0$
to denote this irreducible representation of $Pin (2l)$. Denote 
\begin{eqnarray}
 \Sigma(Pin) = \{ \hf |{\bf 1}_l | + (m_1, m_2, \ldots, \overline{m}_l), 
   m_1 \geq \ldots \geq m_l \geq 0, m_i \in \Z \}.  
  \label{eq_wtpin}
\end{eqnarray}

Let
 \begin{eqnarray}
   \begin{array}{rcl}
    \lefteqn{
      \sum_{i,j \in \Z} ( E_{i,j} -  (-1)^{i +j} E_{-j, -i}) z^i w^{-j} 
            }     \\
     & = & \sum_{k =1}^l \left( : \psi^{+,k} (z) \psi^{-,k} (w):
                          - : \psi^{+,k} ( -w) \psi^{-,k} ( -z):
                  \right).
   \end{array}
  \label{eq_Flbinf}
 \end{eqnarray}
It is known \cite{W2} that the operators 
$E_{i,j} -  (-1)^{i +j} E_{-j, -i} (i, j \in \Z)$
span $\binf$ with central charge $l$ and 
$Pin(2l)$ and $\binf$ form a dual pair on $\Fl$.
We define a map $\Lambda^{\frak{db} }$
from $\Sigma(Pin) $ to ${\binf}_0^*$ by sending 
$ \lambda = (m_1, \ldots, \overline{ m}_l) $ to
 $$ \Lambda^{\frak{db}}(\lambda) =
     (2l - 2j) {}^b\hL_0 + \sum_{k =1}^j {}^b\hL_{m_k}  $$
if $m_1 \geq \ldots \geq m_j > m_{j +1} = \ldots = m_l = 0.$
           
We need to quote the following theorem from \cite{W2}.
\begin{theorem} 
  We have the following $\Fpairdb$-module decomposition:
 \begin{eqnarray*}
  \Fl =  \bigoplus_{\lambda \in \Sigma (Pin) } 
          V(Pin(2l); \lambda) \otimes L \left(\binf;
                                 \Lambda^{\frak{db}}(\lambda), l
                               \right)
 \end{eqnarray*}
 where $V(Pin(2l); \lambda)$ is the irreducible $Pin(2l)$-module
 parametrized by $\lambda \in \Sigma (Pin)$, and
 $L \left(\binf; \Lambda^{\frak{db}} (\lambda), l \right)$ 
 is the irreducible highest weight $\binf$-module of highest weight
 $ \Lambda^{\frak{db}} (\lambda)$ and central charge $l$.
  \label{th_Fpairdb}
\end{theorem}

Note that $\hB$ acts on $\Fl$ via the composition
of the homomorphism $\widehat{\phi}_0$
and the action of $\binf$ given by (\ref{eq_Flbinf}).
We can rewrite the action of $\hB$ in terms of generating function
$T^n (z)$ $(n \in 2 \Z_{+} +1)$ as
$$
 T^n (z) = - \sum_{k =1}^l
           \left( : \psi^{+,k} (z) \partial_z^n \psi^{-,k} (z):
                    - : \partial_z^n \psi^{+,k} ( -z) \psi^{-,k} ( -z):
           \right).
$$

Given $\lambda = (m_1, \ldots, \overline{ m}_l) \in \Sigma(Pin)$, where
$$
 m_1 \geq \ldots \geq m_j > m_{j +1} = \ldots = m_l = 0,
$$
let $e (\lambda)$ be the set of exponents $m_k + 1/2$
$(k = 1, \ldots, j)$ with multiplicity $1$.
The following theorem can be proved in the same way
as Theorem \ref{th_Wpaircc}, based on Theorem \ref{th_Fpairdb} 
and Proposition \ref{prop_deltab}.
\begin{theorem}
  $Pin(2l)$ and $\hB$ form a dual pair on $\Fl$.
 More explicitly we have the following $\Wpairdb$-module decomposition:
 \begin{eqnarray*}
  \Fl = \bigoplus_{\lambda \in \Sigma (Pin) } 
          V(Pin(2l); \lambda) \otimes L \left(\hB;
                                 e (\lambda), l
                               \right).
 \end{eqnarray*}
   \label{th_Wpairpin}
\end{theorem}
\section{FFR's of QHWM's over $\hB$ with $c \in \Bbb N - 1/2$}
   \label{sec_fermionhalf}
 In this section we will realize certain
primitive $\hB$-modules with positive 
half-integral central charges in some Fock spaces.

\vspace{0.15in}
\noindent{\bf I. Case ${\bf  \underline{\Z} = \hz}$:
 dual pair ${\bf \Wpairospc}$}
\vspace{0.1in}
                
    We need a bosonic field
 $ \chi (z) = \sum_{ n \in \hz} \chi_n z^{ -n - \hf}$ which 
satisfies the following commutation relations:
$$[ \chi_m , \chi_n ] = ( -1)^{m + \hf}\delta_{m, -n}, \quad m,n \in \hz .$$
Denote by $\Fminushalf$ the Fock space of $\chi (z)$ generated by
a vacuum vector which is annihilated by $\chi_n, n\in \hz_{+}$.
Let $\Flminushalf$ be the tensor product of the Fock space 
of $l$ pairs of fermionic fields $\psi^{\pm, k} (z) \; (k = 1, \ldots, l )$
and the Fock space $\Fminushalf$ of $\chi (z)$.

It is known \cite{FF} that
the Fourier components of the following ``twisted'' generating functions
together with (\ref{eq_minuspsi}), (\ref{eq_pluspsi}) and (\ref{eq_mixpsi})
\begin{eqnarray*}
 \zeta (z) & \equiv & \sum_{n \in \Z}
            \zeta(n) \NN  = : \chi (z) \chi ( -z):,   \nonumber     \\
 e^{p} (z) & \equiv &\sum_{n \in \Z}
            e^{p}(n) \NN  =  : \psi^{-,p}(z) \chi ( -z):,   \nonumber     \\
 e_*^{p} (z) & \equiv & \sum_{n \in \Z}
           e_*^{p}(n) \NN  = : \psi^{+,p}(z) \chi (z):,  \nonumber    \\
\end{eqnarray*}
generate a representation of the affine 
superalgebra ${\frak {gl}}^{(2)}(1, 2l )$
of type $ A^{(2)}(0, 2l-1)$ \cite{K1} with central charge $ 1$. 
Denote
$$   e^{p} \equiv e^{p}(0), e_*^{p} \equiv e_*^{p}(0), 
e^{pq} \equiv e^{pq}(0), e_*^{pq} \equiv e_*^{pq}(0), 
e_{**}^{pq} \equiv e_{**}^{pq}(0), \;
p, q = 1, \cdots, l. $$ 
The horizontal subalgebra in ${\frak {gl}}^{(2)}(1, 2l )$
spanned by the operators 
$e^{p}, e_*^{p}, e^{pq}, e_*^{pq}, e_{**}^{pq},$ $(p, q = 1, \cdots, l)$
is isomorphic to Lie superalgebra $\frak {osp}(1, 2l)$.
We identify the Borel subalgebra ${\frak b}( {\frak {osp}}(1, 2l) )$
with the one generated by $e_{**}^{pq}, e_*^{pq}\; (p \leq q),\;
p, q = 1, \cdots, l.$ Let

 \begin{eqnarray}
   \begin{array}{rcl}
    \lefteqn{
      \sum_{i,j \in \Z} ( E_{i,j} - ( -1)^{i+j} E_{1-j,1-i}) z^{i-1} w^{-j}
            }     \\
   & = & \sum_{k =1}^l \left( : \psi^{+,k} (z) \psi^{-,k} (w):
                          + : \psi^{+,k} ( -w) \psi^{-,k} ( -z):
                  \right)
                  + :\chi (z) \chi (-w):
  \end{array}
  \label{eq_Flminushalfcinf}
 \end{eqnarray}

It is shown in \cite{W2} that 
$E_{i,j} - ( -1)^{i+j} E_{1-j,1-i} \; (i, j \in \Z)$ span 
$\cinf$ with central charge $ l - \hf$ and $\Fpairospc$ 
form a dual pair on $\Flminushalf$.

Finite dimensional irreducible representations of ${\frak {osp}}(1, 2l)$
are parametrized by the highest weights \cite{K1} 
$$ \Sigma (Osp) =
 \{ (m_1, m_2, \ldots, m_l ) \mid 
   \quad m_1 \geq m_2 \geq \ldots \geq m_l \geq 0 , m_i \in \Z \}. $$
Define a map $\Lambda^{ \frak{ospc}}$ from
$\Sigma (Osp)$ to $ {\cinf}_0^*$ by sending
$\lambda = (m_1, \ldots, m_l )$ to
$$\Lambda^{ \frak{ospc}} (\lambda) =
   (l -j -1/2) {}^c \hL_0 + \sum_{k =1}^j {}^c \hL_{m_k},$$
if $m_1 \geq \ldots \geq m_j > m_{j+1} = \ldots = m_l =0.$
We quote the following theorem from \cite{W2}.
\begin{theorem} 
 We have the following $\Fpairospc$-module decomposition:
 \begin{eqnarray*}
  \Flminushalf =  \bigoplus_{\lambda \in \Sigma (Osp) } 
          V({ \frak{osp}}(1, 2l); \lambda) \otimes L \left(\cinf;
                                 \Lambda^{ \frak{ospc}}(\lambda), l -1/2
                               \right)
 \end{eqnarray*}
 where $V({ \frak{osp}}(1, 2l); \lambda)$ is the irreducible 
 module of ${ \frak{osp}}(1, 2l)$
 parametrized by $\lambda \in \Sigma (Osp)$, and
 $L \left(\cinf; \Lambda^{ \frak{ospc}} (\lambda), l -1/2 \right)$ 
 is the irreducible highest weight $\cinf$-module of highest weight
 $ \Lambda^{ \frak{ospc}} (\lambda)$ and central charge $l -1/2$.
  \label{th_Fpairospc}
\end{theorem}

Note that $\hB$ acts on $\Fl$ via the the composition
of homomorphism $\widehat{\phi}_{1/2}$
and the action of $\cinf$ given by (\ref{eq_Flminushalfcinf}).
We can rewrite the action of $\hB$ in terms of generating function
$ T^n_{1/2} (z)$ $(n \in 2 \Z_{+} +1)$ as
\begin{eqnarray*}
     T^n_{1/2} (z) & =  &
      - \sum_{k =1}^l
           \left( : \psi^{+,k} (z) \partial_z^n \psi^{-,k} (z): 
                  - : \partial_z^n \psi^{+,k} ( -z) \psi^{-,k} ( -z):
           \right)   \\
          & &  - :\chi (z) \partial_z^n \chi (-z):.
\end{eqnarray*}

Given $\lambda  = (m_1, m_2, \ldots, m_l ) \in \Sigma(Osp)$, where
$m_1 \geq \ldots \geq m_j > m_{j +1} = \ldots = m_l =0,$
let $ e (\lambda)$ be the set of exponents
$m_k$ $(k =1, \ldots,j) $ with multiplicity 1 and
the exponent $0$ with multiplicity $l -j - \hf$.
The following duality theorem can be proved in a similar way as
Theorem \ref{th_Wpaircc}, based on Theorem \ref{th_Fpairospc} 
and Proposition \ref{prop_deltac}.
\begin{theorem}
 ${ \frak{osp}}(1, 2l)$ and $\hB$ form a dual pair on $\Flminushalf$. 
  Furthermore we have the following $\Wpairospc$-module decomposition:
 \begin{eqnarray*}
  \Flminushalf = \bigoplus_{\lambda \in \Sigma (Osp) } 
          V({ \frak{osp}}(1, 2l); \lambda) \otimes L \left(\hB;
                                 e (\lambda), l -1/2
                               \right).
 \end{eqnarray*}
\end{theorem}
\vspace{0.15in}
\noindent{\bf II. Case ${\bf \underline{\Z} = \Z}$: 
dual pair ${\bf \Wpairbb }$ }
\vspace{0.1in}

Recall that a fermionic field
 $\varphi (z) = \sum_{n \in \Z} \varphi_n z^{-n}$ which 
satisfy the following commutation relations:
\begin{eqnarray*}
   [ \varphi_m , \varphi_n ]_{+} = ( -1)^m \delta_{m, -n},
     \quad m,n \in \Z. 
\end{eqnarray*}
In this case the Fock space $\Flhalf$ is the tensor product
of the Fock space of $l$ pairs of fermionic fields 
$\psi^{\pm, k} (z), k = 1, \ldots, l$
and the Fock space $\Fhalf$ of $\varphi (z)$
generated by a vacuum vector which is annihilated by
$ \varphi_m, m \in \Bbb N$.

The Fourier components of the following generating functions
together with (\ref{eq_minuspsi}), (\ref{eq_pluspsi}) and (\ref{eq_mixpsi})
\begin{eqnarray*}
 \eta (z) & \equiv & \sum_{n \in \Z}
            \eta(n) z^{-n} =  : \varphi (z) \varphi ( -z):, \nonumber   \\
 e^{p} (z) & \equiv & \sum_{n \in \Z}
            e^{p}(n) z^{-n} = : \psi^{-,p}(z) \varphi ( -z):, \nonumber \\
 e_*^{p} (z) & \equiv & \sum_{n \in \Z}
           e_*^{p}(n) z^{-n} = : \psi^{+,p}(z) \varphi (z):  \nonumber  \\
\end{eqnarray*}
span an affine algebra of type $ A^{(2)}_{ 2l}$
on $\Flhalf$. The horizontal subalgebra of 
$ A^{(2)}_{ 2l}$ is isomorphic to ${\frak {so}}( 2l +1 )$. 
The action of ${\frak {so}}( 2l +1 )$ can be lifted to
an action of $Spin(2l +1)$ on $\Flhalf$. 
It is well known that an irreducible representation of $Spin (2l +1)$ 
which does not factor to $SO (2l +1)$ 
is an irreducible representation of ${\frak {so}} (2l +1)$
parametrized by its highest weight 
\begin{eqnarray}
\lambda = \hf {\bf 1}_l + (m_1, m_2, \ldots, m_l),  
 \quad m_1 \geq \ldots \geq m_l \geq 0.
  \label{eq_spinodd}
\end{eqnarray}
 
Denote 
$$ \Sigma(PB) = \left\{ \hf {\bf 1}_l  + (m_1, m_2, \ldots, m_l ) \mid
   m_1 \geq \ldots \geq m_l \geq 0, m_i \in \Z \right \}.   $$

Let

 \begin{eqnarray}
   \begin{array}{rcl}
    \lefteqn{
     \sum_{i,j \in \Z} ( E_{i,j} - ( -1)^{i+j} E_{-j,-i}) z^i w^{-j}
            }     \\
   & = & \sum_{k =1}^l \left( : \psi^{+,k} (z) \psi^{-,k} (w):
                          + : \psi^{+,k} ( -w) \psi^{-,k} ( -z):
                  \right)
                  + :\varphi (z) \varphi (-w):
   \end{array}
 \label{eq_Flhalfbinf}
 \end{eqnarray}

The Fock space $\Flhalf$ splits into a sum of
two subspaces $\Flhalf_e$ and $\Flhalf_o$, where
$\Flhalf_e$ consists all even vectors while $\Flhalf_o$
consists all odd vectors according to the $ \Z_2$ gradation
on the vector superspace $\Flhalf$. 
The action of ${\frak {so}}( 2l +1 )$
can be lifted to $Spin( 2l)$ on $\Flhalf_e$ and $\Flhalf_o$.
$Spin( 2l)$ and $\binf$ form a dual pair on $\Flhalf_e$ and $\Flhalf_o$ 
\cite{W2}.

Note that $\hB$ acts on $\Flhalf$ via the composition
of the homomorphism $\widehat{\phi}_0$
and the action of $\binf$ given by (\ref{eq_Flhalfbinf}).
We can rewrite the action of $\hB$ in terms of generating function
$T^n (z)$ $(n \in 2 \Z_{+} +1)$ as
\begin{eqnarray*}
     T^n (z) & =  &
      - \sum_{k =1}^l
           \left( : \psi^{+,k} (z) \partial_z^n \psi^{-,k} (z): 
                  + ( -1)^n : \partial_z^n \psi^{+,k} ( -z) \psi^{-,k} ( -z):
           \right)   \\
          & &  - (-1)^n :\varphi (z) \partial_z^n \varphi (-z):.
\end{eqnarray*}
                 
Given $\lambda = \hf {\bf 1}_l  + (m_1, m_2, \ldots, m_l ) \in \Sigma(PB)$,
where 
$$
 m_1 \geq \ldots \geq m_j > m_{j +1} = \ldots = m_l =0,
$$
let $V(Spin(2l+1); \lambda)$ be the irreducible $Spin(2l+1)$-module
 parametrized by $\lambda$ and
let $e (\lambda)$ be the set of exponents $m_k + 1/2$ $(k =1, \ldots, j)$
with multiplicity $1$ 
(as in Section \ref{sec_fermion}).
The following duality theorem on the commuting actions
of $Pin(2l)$ and $\hB$ follows from a corresponding duality theorem 
of a dual pair $\Fpairbb$ in \cite{W2} and similar argument
as in Theorem~\ref{th_Wpaircc}. 
\begin{theorem}
   We have the following $\Wpairbb$-module decomposition:
 \begin{eqnarray*}
   \Flhalf_e = \bigoplus_{\lambda \in \Sigma(PB) } 
          V(Spin(2l+1); \lambda) \otimes L \left(\hB;  
                                 e (\lambda), l+1/2
                               \right)         \\
   \Flhalf_o = \bigoplus_{\lambda \in \Sigma(PB) } 
          V(Spin(2l+1); \lambda) \otimes L \left(\hB;  
                                 e  (\lambda), l+1/2
                               \right) .
 \end{eqnarray*}
       \label{th_Wpairoddpin}
\end{theorem}
\section{FFR's of QHWM's over $\hBp$ with $ c \in \Bbb N$}
   \label{sec_n}
 
Let us take a pair of fermionic fields
$$ 
 \psi^{\pm}(z) = \sum_{n \in \underline{\Z} }
  \psi^{\pm}_n z^{ -n -\hf + \epsilon},
    \quad \underline{\Z} = \hf + \Z \mbox{ or } \Z.
$$
In the case $\underline{\Z} = \hz$ the anti-commutation
relations among $\psi^{\pm}_n$ is equivalent to the following
operator product expansions (OPE)
$$ \psi^{+}(z) \psi^{-}(w) \sim \frac{1}{z-w}, \quad
   \psi^{+}(z) \psi^{+}(w) \sim 0, \quad
   \psi^{-}(z) \psi^{-}(w) \sim 0.
$$
                 
Take $l$ pairs of fermionic fields, 
$\psi^{\pm,p} (z) \;( p = 1, \dots, l)$
and consider the corresponding Fock space $\Fl$.
Introduce the following generating functions
\begin{eqnarray}
  E(z,w) & \equiv & \sum_{i,j \in \Z} E_{ij} z^{i-1 + 2\epsilon }
     w^{-j  } 
   = \sum_{p =1}^l :\psi^{+,p} (z) \psi^{-,p} (w):  \label{eq_genelin} \\
  e^{pq} (z) & \equiv & \sum_{n \in \Z} e^{pq} (n) z^{ -n-1 + 2\epsilon }
    =  : \psi^{-,p}(z) \psi^{-,q} (z): 
    \quad (p \neq q )   \label{eq_orthaffine}    \\ 
  e_{**}^{pq} (z) & \equiv & \sum_{n \in \Z}
                            e_{**}^{pq}(n)z^{ -n-1 + 2\epsilon }
    =  : \psi^{+,p}(z) \psi^{+,q} (z):
    \quad (p \neq q )  \label{eq_orth}    \\
  e_*^{pq} (z) & \equiv & \sum_{n \in \Z} e_*^{pq} (n) z^{ -n-1 + 2\epsilon }
    =  : \psi^{+,p}(z) \psi^{-,q} (z): +  \delta_{p,q} \epsilon 
         \label{eq_orthaffinestar}  
\end{eqnarray}
where $p,q = 1, \dots, l$, and the normal ordering $::$ means that 
the operators annihilating $\vac$ are moved to the right and 
multiplied by $-1$.

It is well known \cite{F1, KP} that the operators $e^{pq}(n), e_*^{pq}(n), 
e_{**}^{pq}(n)$, $p,q = 1, \dots, l$, $n \in \Z$ 
form a representation of the affine algebra
$\widehat{so}(2l)$ of level $ 1$. 
The operators $e^{pq}(0), e_*^{pq}(0), e_{**}^{pq}(0)
 \;(p, q = 1, \cdots, l)$ form
the horizonal subalgebra $\frak {so}(2l)$ in $\widehat{ \frak {so}}(2l)$.
We identify the Borel subalgebra ${\frak b}( {\frak {so}}(2l) )$
with the one generated by $e_{**}^{pq} \;(p \neq q), e_*^{pq}\; (p \leq q),
p, q = 1, \cdots, l.$

From now on we need to treat the two cases 
$\underline{\Z} = \hf + \Z \mbox{ or } \Z $ separately.
First consider the case $\underline{\Z} = \hf + \Z.$
It follows from (\ref{eq_genelin}) that
 \begin{eqnarray}
   \begin{array}{rcl}
     \lefteqn{
      \sum_{i,j \in \Z} ( E_{ij} - E_{1-j,1-i} ) z^{i- 1} w^{-j} =
             }          \\
    & & \sum^l_{k=1} (: \psi^{+,k} (z) \psi^{-,k} (w) : 
    - : \psi^{+,k} (w) \psi^{-,k} (z) :).
   \label{eq_gener}
   \end{array}
 \end{eqnarray}
One can show that $ E_{ij} - E_{1-j, 1-i} \;( i,j \in \Z)$
span $\dinf$. The action of the horizontal subalgebra ${\frak {so}}(2l)$
can be integrated into an action of $SO(2l)$ and extended to an
action of $O(2l)$ naturally. It is known \cite{W2} that
the action of $\dinf$ commutes with the action of $ O(2l)$
on $\Fl$ and moreover $\dinf$ and $ O(2l)$ form a dual 
pair\footnote{This fact was also known to V. Kac and C. Yan.}
by the same argument as in finite dimensional dual pairs case \cite{H1, H2}.

We define a map 
$\Lambda^{ \frak{dd}}: \Sigma (D) \longrightarrow {\dinf}_0^* $
by sending $ \lambda = (m_1, \cdots, \overline{ m}_l) \; (m_l > 0 )$ to 
$$  \Lambda^{ \frak{dd}} (\lambda) =
   (l -i ) {}^d \hL_0 + ( l -i ) {}^d \hL_1 + \sum_{k =1}^i {}^d \hL_{m_k}, 
$$
sending $(m_1, \cdots, m_j, 0, \ldots, 0 )\; (j <l )$ to 
$$  \Lambda^{ \frak{dd}} (\lambda) =
   (2l -i -j) {}^d \hL_0 + (j-i ) {}^d \hL_1 + \sum_{k =1}^i {}^d \hL_{m_k}, 
$$
and sending $ (m_1, \ldots, m_j, 0, \ldots, 0 ) \bigotimes {det} 
   \; ( j < l)$ to 
 $$  \Lambda^{\frak{dd}} (\lambda) =
  (j-i) {}^d \hL_0 + (2l -i -j) {}^d \hL_1 + \sum_{k =1}^i {}^d \hL_{m_k}, $$
if $m_1 \geq \ldots m_{i} > m_{i+1} = \ldots = m_{j} =1
   > m_{j+1} = \ldots = m_l = 0.$

The following theorem was proved in \cite{W2}.
\begin{theorem}
  We have the following $(O(2l), \dinf )$-module decomposition:
 \begin{eqnarray*}
   \Fl = \bigoplus_{\lambda \in \Sigma (D) } 
          V(O(2l); \lambda) \otimes L \left(\dinf;
                                 \Lambda^{\frak{dd}} (\lambda), l
                               \right)
 \end{eqnarray*}
 where $V(O(2l); \lambda)$ is the irreducible $O(2l)$-module
 parametrized by $\lambda \in \Sigma (D)$ and 
 $L \left(\dinf; \Lambda^{\frak{dd}} (\lambda), l  \right)$
 is the irreducible 
 highest weight $\dinf$-module of highest weight
 $ \Lambda^{\frak{dd}} (\lambda)$ and central charge $l$.
  \label{th_Fpair}
\end{theorem}
We can obtain the action of $\hBp$ on $\Fl$ by composing
the action of $\dinf$ and the homomorphism $\hphio$
given by the formula (\ref{eq_embd}). Introduce
the generating function
\begin{eqnarray}
  W^n (z) = \sum_{k \in \Z} W^n_k z^{ -k -n -1} \quad (n \in 2 \Z_{+} +1).
 \label{eq_field}
\end{eqnarray}
\begin{lemma}
   On $\Fl$ we have
 \begin{eqnarray} 
 W^n (z) = \hf \sum^l_{k=1} (: \partial_z^n \psi^{-,k} (z) \psi^{+,k} (z) : 
    + : \partial_z^n \psi^{+,k} (z) \psi^{-,k} (z) :).
   \label{eq_ffr}
 \end{eqnarray}
  \label{lem_ffr}
\end{lemma}
\begin{demo}{Proof}
 We calculate $W^n (z)$ as follows.
 \begin{eqnarray}
   W^n (z) & = & - \hf \sum_{k, j \in \Z}
           ( [ -j]_n - [ j -k -1]_n ) E_{j -k, j} z^{ -k -n -1} 
                \label{eq_free1} \\
           & = &  - \hf \sum_{k, j \in \Z}
            [ -j]_n (E_{j -k, j} - E_{1 -j, k -j +1})  z^{ -k -n -1} 
                \label{eq_free2}  \\
           & = & \hf \sum^l_{k=1}
             (- : \psi^{+,k} (z) \partial_z^n \psi^{-,k} (z) : 
    + : \partial_z^n \psi^{+,k} (z) \psi^{-,k} (z) :).
                \label{eq_free3} 
 \end{eqnarray}
 Here (\ref{eq_free1}) is given by (\ref{eq_field}) and
 (\ref{eq_embd}), (\ref{eq_free2}) is
 obtained by shifting the indices from $j -k -1$
 to $j$ in the second part of
 the right hand side of (\ref{eq_free1}), and  
 (\ref{eq_free3}) is obtained by taking $n$-th derivatives 
 of (\ref{eq_gener}) with respect to $w$. It is clear
 that (\ref{eq_free3}) is the same as (\ref{eq_ffr}).
\end{demo}

For $ \lambda = (m_1, \cdots, \overline{ m}_l)\in \Sigma(D) $ where
$m_1 \geq \ldots \geq m_i > m_{i +1} = \ldots = m_l = 1$,
we let $e (\lambda)$ be the set of exponents $m_k$ $(k =1, \ldots, i)$
with multiplicity $1$ and the exponent $1$ with multiplicity $l -i$;
for $\lambda = (m_1, \cdots, m_l)\in \Sigma(D)$ where 
$$
  m_1 \geq \ldots m_{i} > m_{i+1} = \ldots = m_{j} =1
   > m_{j+1} = \ldots = m_l = 0 \;( j < l),
$$ 
we let $e (\lambda)$ be the set of 
exponents $m_k$ $(k =1, \ldots, i)$ of multiplicity $1$,
exponent $1$ of multiplicity $j -i$;
for $ (m_1, \cdots, m_l) \bigotimes {det} \in \Sigma(D)$ where 
$$
 m_1 \geq \ldots m_{i} > m_{i+1} = \ldots = m_{j} =1
   > m_{j+1} = \ldots = m_l = 0\; ( j < l),
$$
we let $e (\lambda)$ be the set of 
exponents $m_k$ $(k =1, \ldots, i)$ of multiplicity $1$,
exponent $1$ of multiplicity $2l -i -j$.
We will simply write $(0, \cdots, 0 ) $
and $(0, \cdots, 0 ) \bigotimes {det}$ as $0$ and $det$ respectively.
\begin{theorem}
   $O(2l)$ and $\hBp$ form a dual pair on $\Fl$. 
Moreover we have the following $\Wpairdd$-module decomposition:
 \begin{equation}
   \Fl = \bigoplus_{\lambda \in \Sigma (D)} 
          V( O(2l); \lambda) \otimes L \left(\hBp;
                                 e (\lambda), l
                               \right).
 \end{equation}
  \label{th_Wpairdd}
\end{theorem}
\begin{demo}{Proof}
   By Theorem \ref{th_pullback} the $\dinf$-module
 $L \left(\dinf; \Lambda^{\frak{dd}} (\lambda), l  \right)$
 regarded as a $\hBp$-module via the pullback by $\hphio$
 remains irreducible. Then this theorem follows from 
 Theorem \ref{th_Fpair} once we determine the corresponding
 $\Delta (x)$ for this $\hBp$-module. 
 It follows from the definition of $ \Lambda^{\frak{dd}} (\lambda)$
 by using Proposition~\ref{prop_deltad} that 
 for $ \lambda = (m_1, \cdots, \overline{ m}_l)\in \Sigma(D) $ where
 $m_1 \geq \ldots \geq m_i > m_{i +1} = \ldots = m_l = 1$,
 \begin{eqnarray*}
 \Delta (x) & = & \sum_{k =1}^i
        \frac {\cosh (m_k x) + (l -i) \cosh x +
          \hf ( (l -i) - (l -i)) }{2 \sinh (x/2)}
        - \frac {l}{2 \sinh (x/2)}.
 \end{eqnarray*}
 The computations of $\Delta (x)$ for the remaining
 $\lambda \in \Sigma(D)$ are similar.
\end{demo}

We have an immediate corollary.
\begin{corollary}
   The space of 
 invariants of $O(2l)$ in the Fock space $\Fl$ is 
 naturally isomorphic to the irreducible 
 module $ L (\dinf; 2l\; {}^d \hL_0 )$ of central charge $l$
 with highest weight vector $\vac$,
 or equivalently to the irreducible module 
 $L \left(\hBp; e(0), l \right)$.
   \label{cor_dl}
\end{corollary}
\begin{remark} \rm
  The Dynkin diagram of $\dinf$ admits an automorphism of order 2 denoted
 by $\sigma$. $\sigma$ induces naturally an automorphism of order 2
 of $\dinf$, which is denoted again by $\sigma$ by abuse of notation.
 $\sigma$ acts on the set of highest weights of $\dinf$
 by mapping $ \lambda = \; {}^d h_0 \; {}^d  \hL_0 + \; {}^d h_1 \;{}^d \hL_1 
  + \sum_{i \geq 2} \; {}^d h_i  \; {}^d \hL_i$ to
 $\sigma(\lambda ) = \; {}^d h_1  \; {}^d \hL_0 + \; {}^d h_0  \; {}^d \hL_1 
  + \sum_{i \geq 2} \; {}^d h_i  \; {}^d \hL_i$. 
 In this way one can obtain an irreducible module of the semi-product
 $\sigma \ltimes \dinf$
 on $ L(\dinf; \lambda) \bigoplus L(\dinf; \sigma(\lambda ) )$
 if $\sigma(\lambda ) \neq \lambda $ and on $ L(\dinf; \lambda) $
 if $\sigma(\lambda ) = \lambda $.

 It was noted in \cite{W2} that 
 $( SO(2l), \sigma \ltimes \dinf)$ form a dual pair on $\Fl$. In particular
 the space of invariants of $\Fl$ under the action fo $SO(2l)$ is 
 isomorphic to the $\dinf$-module
 $L(\dinf; 2l \; {}^d \hL_0 ) \bigoplus L(\dinf; 2l \; {}^d \hL_1 )$
 or equivalently the $\hBp$-module $L \left(\hBp; e(0), l \right)
 \bigoplus L \left(\hBp; e (det), l \right)$. 
   \label{rem_semi}
\end{remark}
          
Now we consider the case $\underline{\Z} = \Z.$
It follows from (\ref{eq_genelin}) that
 \begin{eqnarray}
   \begin{array}{rcl}
     \lefteqn{
      \sum_{i,j \in \Z} ( E_{ij} - E_{-j,-i} ) z^{i} w^{-j} =
             }          \\
    & & \sum^l_{k=1} (: \psi^{+,k} (z) \psi^{-,k} (w) : 
    - : \psi^{+,k} (w) \psi^{-,k} (z) :).
   \end{array}
 \end{eqnarray}
One can show that $ E_{ij} - E_{1-j, 1-i} \;( i,j \in \Z)$
span $\binftwo$. The action of the horizontal subalgebra ${\frak {so}}(2l)$
can be integrated into an action of $Spin(2l)$ (cf. e.g. \cite{BtD}
for more on spin groups)
and then extended naturally
to $Pin (2l)$. $Pin (2l)$ and $\binftwo$ form a dual pair
on $\Fl$ \cite{W2}.

We define a map $\Lambda^{\frak{db} }$
from $\Sigma(Pin) $ (see (\ref{eq_wtpin}) for notation)
to ${\binf}_0^*$ by sending 
 $$ \lambda = (m_1, \ldots, \overline{ m}_l) , \quad
 m_1 \geq m_2 \geq \ldots \geq m_l \geq 0$$
to
 $$ \Lambda^{\frak{db}}(\lambda) =
     (2l - 2j) {}^b \hL_0 + \sum_{k =1}^j {}^b \hL_{m_k} $$
if $m_1 \geq \ldots m_j > m_{j +1} = \ldots = m_l = 0.$

The following theorem was proved in \cite{W2}.
\begin{theorem}
 We have the following $\Fpairdbtwo$-module decomposition:
 \begin{eqnarray*}
  \Fl =  \bigoplus_{\lambda \in \Sigma (Pin) } 
          V(Pin(2l); \lambda) \otimes L \left(\binftwo;
                                 \Lambda^{\frak{db}}(\lambda), l
                               \right)
 \end{eqnarray*}
 where $V(Pin(2l); \lambda)$ is the irreducible $Pin(2l)$-module
 parametrized by $\lambda \in \Sigma (Pin)$, and
 $L \left(\binftwo; \Lambda^{\frak{db}} (\lambda), l \right)$ 
 is the irreducible highest weight $\binftwo$-module of highest weight
 $ \Lambda^{\frak{db}} (\lambda)$ and central charge $l$.
   \label{th_pinpair}
\end{theorem}

We define the action of $\hBp$ on $\Fl$ by the composition
of the action of $\binftwo$ and the homomorphism $\hphic$
given by (\ref{eq_embd}). It follows that
the action of $Pin (2l)$ commutes with that of $\hBp$.
Introduce the following generating function
\begin{equation}
  W^n_{1/2} (z) = \sum_{k \in \Z} W^{n, 1/2}_k z^{ -k -n }  \label{eq_tw}
\end{equation}
It follows from Proposition \ref{prop_ext} that
the representation of $\hBp$ on $\Fl$ has central charge $l$ and
\begin{eqnarray}
 \widehat{\phi}_{- 1/2}
   \left( e^{xD} - e^{-xD} \right)
       = \phi_{- 1/2} \left( e^{xD} - e^{-xD} \right) -l \tanh(x/4).
  \label{eq_diffd}
\end{eqnarray}
Therefore
\begin{eqnarray*}
 \widehat{\phi}_{ -1/2}(W_0^{n, 1/2}) = \phi_{ -1/2}(W_0^{n, 1/2}) - \alpha_n l
\end{eqnarray*}
for some constant $\alpha_n$ determined by (\ref{eq_diffd})
depending on $n$ only. Similarly as in Lemma \ref{lem_ffr}
one can recast the action of $\hBp$ in terms of $ W_{ 1/2}^n (z)$ as
 \begin{eqnarray}
  \begin{array}{rcl}
   \lefteqn{
  W_{ 1/2}^n (z)} \\
   & = &
     \hf \sum_{k =1}^l
           \left( : \partial_z^n \psi^{-,k} (z) \psi^{+,k} (z):
             +  : \partial_z^n \psi^{+,k} (z) \psi^{-,k} (z):
           \right) - \alpha_n l z^{ -n }.
  \end{array}
   \label{eq_twfield}
 \end{eqnarray}

Given $ \lambda = (m_1, \ldots, \overline{ m}_l) \in \Sigma (Pin)$ where
$ m_1 \geq m_2 \geq \ldots \geq m_l \geq 0$, we let 
$e (\lambda)$ be the set of 
exponents $m_k + \hf$ $(k =1, \ldots, j)$ of multiplicity $1$ and 
exponent $ \hf$ of multiplicity $l -j$.

The proof of the following theorem is obtained
in an analogous way as for Theorem \ref{th_Wpairdd} by using
Theorems \ref{th_pullback}, \ref{th_pinpair}
and Proposition~\ref{prop_dtab}. 
\begin{theorem}
    We have the following $\Wpairdbtwo$-module decomposition:
 \begin{eqnarray*}
  \Fl =   \bigoplus_{\lambda \in \Sigma (Pin) } 
          V(Pin(2l); \lambda) \otimes L \left(\hBp;
                                 e (\lambda), l
                               \right).
 \end{eqnarray*}
\end{theorem}
\section{FFR's of QHWM's over $\hBp$ with $ c \in \hf + \Bbb Z_{+}$}
   \label{sec_nhalf}
   Introduce a neutral fermionic field 
$\phi (z) = \sum_{n \in \underline{\Z} } \phi_n z^{-n -\hf + \epsilon}$ 
which satisfies the following commutation relations:
$$[ \phi_m , \phi_n ]_{+} = \delta_{m, -n},
     \quad m,n \in \underline{\Z}.$$
Denote by $\Fhalf$ the Fock space of $\phi (z)$
generated by a vacuum vector $\vac$, which is annihilated by 
$\phi_n, n \in \underline{\Z}_{+}$.
Denote by $\Flhalf$ the $\Z_2$-graded tensor product of $\Fhalf$ 
and the Fock space $\Fl$ of $l$ pairs of fermionic fields 
$\psi^{\pm, k} (z) (k =1, \dots, l)$.

Denote by 
\begin{eqnarray}
 e^{p} (z) \equiv \sum_{n \in \Z}
 e^{p}(n) z^{-n -1 + 2 \epsilon} & = & : \psi^{-,p}(z) \phi (z):
             \nonumber       \\
 e_*^{p} (z) \equiv \sum_{n \in \Z}
 e_*^{p}(n) z^{-n -1 + 2 \epsilon} & = & : \psi^{+,p}(z) \phi (z):
    \quad (p = 1, \dots, l).         \label{eq_oddorth}
\end{eqnarray}
Then the Fourier components of $e^p (z), e^p_* (z)$, together with
generating functions (\ref{eq_orthaffine}),
(\ref{eq_orth}) and (\ref{eq_orthaffinestar})
$$
 e_*^{pq}(n),\; e^{pq}(n) (p \neq q),\; e_{**}^{pq}(n) (p \neq q),\;
 e^{p} (n), e_*^{p}(n) \quad ( n \in \Z, \; p, q = 1, \dots, l) 
$$
generate an affine algebra $\widehat{{\frak {so}}}(2l+1)$ of level $1$ 
\cite{F1, KP}.  
$ e^{pq}(0)\; (p \neq q), e_{**}^{pq}(0)\; (p \neq q),$
$e_*^{pq}(0), e^{p}(0), e_*^{p}(0)$ $( p, q = 1, \dots, l) $ 
generate the horizontal subalgebra 
${\frak {so}}(2l+1)$ of $\widehat{{\frak {so}} }(2l+1)$.
In particular, we identify the Borel subalgebra
${\frak b}( {\frak {so}}(2l +1) )$ with the one generated by 
$e_{**}^{pq}(0) (p \neq q), e_*^{pq}(0) (p \leq q), \;
e_*^p(0), p, q = 1, \cdots, l.$

From now on we need to consider the two cases 
$\underline{\Z} = \hf + \Z \mbox{ and } \Z $ separately.

First consider the case $\underline{\Z} = \hf + \Z.$
Introduce a generating function
 \begin{eqnarray}
   \begin{array}{rcl}
     \lefteqn{
   \sum_{i,j \in \Z} ( E_{ij}
    - E_{1-j,1-i} ) z^{i- 1} w^{-j}
             }             \nonumber \\
   & & = \sum^l_{k=1} (: \psi^{+,k} (z) \psi^{-,k} (w) :
    - : \psi^{+,k} (w) \psi^{-,k} (z) :) 
    + :\phi (z) \phi (w):.
   \end{array}
 \end{eqnarray}
 One can show that this
defines an action of $\dinf$ on $\Flhalf$ with central charge $l + \hf$.
The action of the horizontal subalgebra 
${\frak {so}}(2l +1)$ can be lifted to an action of $O(2l +1)$. 
The action of $O(2l +1)$ commutes with that of $\dinf$ generated by 
$ E_{ij} - E_{1-j, 1-i} \; (i,j \in \Z )$ on $\Flhalf$.

Define a map $\Lambda^{ \frak{bd} } $ from $\Sigma (B)$ 
to ${\dinf}_0^*$ by sending
$$\lambda = (m_1, m_2, \ldots, m_l)$$
 to
$$  \Lambda^{ \frak{bd} } =
   (2l+1 - i -j) {}^d \hL_0 + (j-i){}^d \hL_1 + \sum_{k=1}^i {}^d \hL_{m_k}  
$$
and sending
$$\lambda = (m_1, m_2, \ldots, m_l) \bigotimes det $$
 to
$$   \Lambda^{ \frak{bd} } =
    (j - i) {}^d \hL_0 + (2l +1-i -j){}^d \hL_1
    + \sum_{k =1}^i {}^d \hL_{m_k} $$
assuming that 
$$ m_1 \geq \ldots \geq m_i > m_{i +1} = \ldots = m_j = 1 >
m_{j+1} = \ldots = m_l = 0 .$$
 
The following theorem is quoted from \cite{W2}.
\begin{theorem}
 We have the following $\Fpairbd$-module decomposition:
 \begin{eqnarray*}
   \Flhalf = \bigoplus_{\lambda \in \Sigma(B) } 
          V(O(2l+1); \lambda) \otimes L \left(\dinf;  
                                 \Lambda^\frak{bd}_{+}(\lambda), l +1/2
                               \right)
 \end{eqnarray*}
 where $V(O(2l+1); \lambda)$ is the irreducible $O(2l+1)$-module
 parametrized by $\lambda  \in \Sigma(B) $ and 
 $L \left(\dinf; \Lambda^\frak{bd} (\lambda), l+1/2 \right)$
 is the irreducible highest weight $\dinf$-module with 
 highest weight $ \Lambda^\frak{bd} (\lambda)$ and central charge $l +1/2$.
  \label{th_Fpairbd}
 \end{theorem}
 
The action of $\hBp$ on $\Flhalf$ is given by the composition
of the action of $\dinf$ on $\Flhalf$ and the homomorphism $\hphio$.
Similarly as in Lemma \ref{lem_ffr} we can show that
\begin{eqnarray}
 \begin{array}{rcl}
  \lefteqn{
  W^n (z) = }    \\
  && \hf \sum^l_{k=1} \left(
  : \partial_z^n \psi^{-,k} (z) \psi^{+,k} (z) : 
    + : \partial_z^n \psi^{+,k} (z) \psi^{-,k} (z) : \right)
    + \hf : \partial_z^n \phi(z) \phi (z) : .
 \end{array}
   \label{eq_bl}
\end{eqnarray}
For $\lambda = (m_1, m_2, \ldots, m_l) \in \Sigma(B)$
we let $e (\lambda) $ be the set of
exponents $m_k$ $(k =1, \ldots, i)$ of multiplicity $1$, the
exponent $1$ of multiplicity $j -i$ 
(the exponent $0$ has multiplicity $l -j + \hf$);
and for $\lambda = (m_1, m_2, \ldots, m_l) \bigotimes det $
we let $e (\lambda) $ be the set of 
exponents $m_k$ $(k =1, \ldots, i)$ of multiplicity $1$, the
exponent $1$ of multiplicity $2l +1 -i -j$ 
(the exponent $0$ has multiplicity $-l +j - \hf$), where
$$ m_1 \geq \ldots \geq m_i > m_{i +1} = \ldots = m_j = 1 >
m_{j+1} = \ldots = m_l = 0 .$$

The following theorem can be proved in an analogous way
as Theorem~\ref{th_Wpairdd} by using 
Theorems \ref{th_pullback}, \ref{th_Fpairbd}
and Proposition~\ref{prop_deltad}. 

\begin{theorem}
 We have the following $\Wpairbd$-module decomposition:
 \begin{eqnarray*}
   \Flhalf = \bigoplus_{\lambda \in \Sigma(B) } 
          V(O(2l+1); \lambda) \otimes L \left(\hBp;
                                 e (\lambda), l +1/2
                               \right).
 \end{eqnarray*}
  \label{th_Wpairbd}
\end{theorem}
 
The following corollary is immediate.
\begin{corollary}
    The space of invariants of $O(2l +1)$ in $\Flhalf$ is naturally
 isomorphic to
 the irreducible $\dinf$-module $L ( \dinf; (2l +1) \; {}^d \hL_0 )$
 or equivalently to the irreducible $\hBp$-module 
 $L \left(\hBp; e (0), l+1/2 \right)$.
   \label{cor_bl}
\end{corollary}
\begin{remark} \rm
 $( SO(2l +1), \sigma \ltimes \dinf)$ form a dual pair on $\Flhalf$.
 In particular the space of invariants of $\Flhalf$ with respect to
 $SO(2l +1)$ is isomorphic to the $\dinf$-module
 $L(\dinf; (2l+1) \; {}^d \hL_0 ) \bigoplus L(\dinf;  (2l+1) \; {}^d \hL_1 )$
 or equivalently the $\hBp$-module 
 $L (\hBp; e (0), l+1/2 )
  \bigoplus L (\hBp; e (det), l+1/2 )$.
\end{remark}    

Now we consider the case $\underline{\Z} = \Z.$
Introduce the following generating function
 \begin{eqnarray*}
   \begin{array}{rcl}
    \lefteqn{
      \sum_{i,j \in \Z} ( E_{i,j} -  E_{-j,-i}) z^i w^{-j} 
            }     \\
     & = & \sum_{k =1}^l \left( : \psi^{+,k} (z) \psi^{-,k} (w):
                          - : \psi^{+,k} ( w) \psi^{-,k} ( z):
                  \right)
                  + :\phi (z) \phi (w): .
   \end{array}
 \end{eqnarray*}
This defines a representation of $\binftwo$ on $\Flhalf$ of central
 charge $l +1/2$.

The Fock space $\Flhalf$ splits into a sum of
two subspaces $\Flhalf_e$ and $\Flhalf_o$, where
$\Flhalf_e$ consists all even vectors while $\Flhalf_o$
consists all odd vectors according to the $ \Z_2$ gradation
on the vector superspace $\Flhalf$. 
The action of ${\frak {so}}( 2l +1 )$
can be lifted to $Spin( 2l)$ on $\Flhalf_e$ and $\Flhalf_o$ respectively.
$Spin( 2l)$ and $\binf$ form a dual pair on $\Flhalf_e$ and $\Flhalf_o$ 
\cite{W2}.

Define a map $\Lambda^{\frak{bb} }$ from $\Sigma (PB)$ to 
${\binf}_0^*$ by sending
$\lambda = \hf {\bf 1}_l + (m_1, m_2, \ldots, m_l)$ to
$$  \Lambda^{\frak{bb} } (\lambda)
     = (2l +1 - 2j ) {}^b \hL_0 + \sum_{k=1}^j {}^b \hL_{m_k} $$
if $m_1 \geq \ldots \geq m_j > m_{j+1} = \ldots = m_l =0.$
                 
The following theorem was proved in \cite{W2}.
\begin{theorem}
  We have the following $\Fpairbbtwo$-module decomposition:
 \begin{eqnarray*}
   \Flhalf_e & = & \bigoplus_{\lambda \in \Sigma(PB) } 
          V(Spin(2l+1); \lambda) \otimes L \left(\binftwo;  
                                 \Lambda^{\frak{bb}} (\lambda), l + \hf
                               \right)                 \\
   \Flhalf_o & = & \bigoplus_{\lambda \in \Sigma(PB) } 
          V(Spin(2l+1); \lambda) \otimes L \left(\binftwo;  
                                 \Lambda^{\frak{bb}} (\lambda), l + \hf
                               \right)  
 \end{eqnarray*}
 where $V(O(2l+1); \lambda)$ is the irreducible $Spin(2l+1)$-module
 parametrized by 
 $\lambda = \hf {\bf 1}_l + (m_1, m_2, \ldots, m_l)$ and 
 $ L \left(\binftwo; \Lambda^{\frak{bb}} (\lambda), l + \hf \right)$ 
 is the irreducible 
 highest weight $\binftwo$-module with central charge $l + \hf$.
   \label{th_Fpairbbtwo}
 \end{theorem}

Now the action of $\hBp$ on $\Flhalf_e$ and $\Flhalf_o$
can be obtained by the composition of
the action of $\binftwo$ and the homomorphism $\hphic$
given by (\ref{eq_embd}). In terms of 
$W^n_{1/2} (z) = \sum_{k \in \Z} W^{n, 1/2}_k z^{ -k -n }$
we obtain in an anagolous way as in (\ref{eq_twfield}) that 
 \begin{eqnarray*}
  W_{ 1/2}^n (z) &=&
    \hf \sum_{k =1}^l
           \left( : \partial_z^n \psi^{-,k} (z) \psi^{+,k} (z):
             +  : \partial_z^n \psi^{+,k} (z) \psi^{-,k} (z):  \right)  \\
     & &       + \hf : \partial_z^n \phi(z) \phi (z) :
     - \alpha_n (l +1/2) z^{ -n } .
 \end{eqnarray*}

For $\lambda = \hf {\bf 1}_l + (m_1, m_2, \ldots, m_l) \in \Sigma (PB)$ where 
$$m_1 \geq \ldots \geq m_j > m_{j +1} = \ldots = m_l = 0, $$
we let $e (\lambda)$ be the set of
exponents $m_k + \hf$ $(k =1, \ldots, j)$ of multiplicity $1$ and the
exponent $ \hf$ of multiplicity $l + 1/2 -j$.
The following theorem can be proved in an analogous way
as Theorem~\ref{th_Wpairdd} by using instead
Theorems \ref{th_pullback}, \ref{th_Fpairbbtwo}
and Proposition~\ref{prop_dtab}. 
\begin{theorem}
    We have the following $\Wpairbbtwo$-module decomposition:
 \begin{eqnarray*}
   \Flhalf_e & = & \bigoplus_{\lambda \in \Sigma(PB) } 
          V(Spin(2l+1); \lambda) \otimes L \left(\hBp;  
                                e (\lambda), l + \hf
                               \right)                 \\
   \Flhalf_o & = & \bigoplus_{\lambda \in \Sigma(PB) } 
          V(Spin(2l+1); \lambda) \otimes L \left(\hBp;  
                                 e (\lambda), l + \hf
                               \right)  .
 \end{eqnarray*}
\end{theorem}
\section{FFR's of QHWM's over $\hBp$ with $ c \in - \hf \Bbb N$}
   \label{sec_minusn}
\subsection{Case $ c \in - \Bbb N$}
       Let us take a pair of bosonic ghost fields
$$ \gamma^{\pm} (z) = \sum_{ n \in \hz }
                         \gamma^{\pm}_n z^{-n - \hf} .$$
Equivalently we have the following operator product expansions
$$
   \gamma^+ (z)\gamma^- (w) \sim \frac{1}{z-w}, \quad
   \gamma^+ (z)\gamma^+ (w) \sim 0, \quad
   \gamma^- (z)\gamma^- (w) \sim 0.
$$

We take $l$ pairs of bosonic ghost fields 
$\gamma^{\pm ,p} (z) \;(p = 1, \dots, l)$
and consider the corresponding Fock space $\Fminusl$.

Introduce the following generating functions
\begin{eqnarray}
  E(z,w) & \equiv & \sum_{i,j \in \Z} E_{ij} z^{i -1} w^{-j } 
  = - \sum_{p =1}^l :\gamma^{+,p}  (z) \gamma^{-,p} (w):
                 \label{eq_hglboson} \\
  e^{pq}_{**} (z) & \equiv &
              \sum_{i,j \in \Z} e^{pq}_{**}(n) z^{ -n -1 } 
  = : \gamma^{+,p} (z) \gamma^{+,q} (z): 
    \quad (p \neq q )           \nonumber       \\ 
  e^{pq} (z) & \equiv &
               \sum_{i,j \in \Z} e^{pq}(n)z^{-n -1}
  = : \gamma^{-,p}(z) \gamma^{-,q} (z):
    \quad (p \neq q )      \label{eq_genfunc}     \\
  e_*^{pq} (z) & \equiv & \sum_{i,j \in \Z} e_*^{pq} (n)z^{-n -1}
  = : \gamma^{+,p} (z) \gamma^{-,q} (z):
   \quad (p,q = 1, \dots, l )  \nonumber    
\end{eqnarray}
where the normal ordering $::$ means that the operators
annihilating $\vac$ are moved to the right.

It is well known
that the operators $E_{ij}$ $ (i,j \in \Z )$
form a representation in $\Fminusl$ of
the Lie algebra $\hgl$ with central charge $ -l$;
the operators 
 $$ e^{pq}(n), e_*^{pq}(n), 
  e_{**}^{pq}(n) \quad ( p,q = 1, \dots, l, n \in \Z )$$
form an affine algebra
$\widehat{{\frak {sp}} }(2l)$ with central charge $ -1$ \cite{FF}. 
The operators $e^{pq}(0),$ $ e_*^{pq}(0),$ $ e_{**}^{pq}(0)$
$ (p, q = 1, \cdots l )$ span the horizontal subalgebra 
${\frak {sp}}(2l)$ in $\widehat{{\frak {sp}} }(2l)$.
In particular, operators $e_*^{pq}$ $( p, q = 1, \cdots l )$ form
a Lie subalgebra ${\frak {gl}}(l)$ 
in the horizontal subalgebra ${\frak {sp}}(2l)$.
We identify the Borel subalgebra ${\frak b} ({\frak {sp} }(2l) )$
with the one generated by 
$e_{**}^{pq}(0), e_*^{pq}(0) (p \leq q), \; p, q = 1, \cdots ,l$.
The action of the horizontal subalgebra ${\frak {sp}}(2l)$
can be lifted to an action of Lie group $Sp(2l)$ on $\Fminusl$.

It follows from (\ref{eq_hglboson}) that
\begin{eqnarray*}
   \sum_{i,j \in \Z} (E_{ij} -E_{1-j, 1-i}) z^{i-1} w^{-j } 
  =  \sum_{p =1}^l (:\gamma^{+,p} (w) \gamma^{-,p} (z):
     - :\gamma^{+,p}  (z) \gamma^{-,p} (w):).
\end{eqnarray*}
The operators $E_{ij} -E_{1-j, 1-i}$ $(i,j \in \Z)$
span the Lie algebra $\dinf$ with central charge $-l$.
The actions of $Sp(2l)$ and of $\dinf$ on $\Fminusl$
commutes with each other and form a dual pair.
We now define a map $\Lambda^{\frak{cd} }$ from $\Sigma (C)$ to ${\dinf}_0^*$ 
which maps $\lambda = (m_1, \ldots, m_l) $ to
$$ \Lambda^{\frak{cd} } (\lambda) =
    (-2l -m_1 -m_2) {}^d \hL_0 + \sum_{k =1}^l (m_k -m_{k +1}) {}^d \hL_k $$
with the convention $ m_{l+1} = 0$ here and below.
The following theorem was proved in \cite{W2}.
\begin{theorem}
  We have the following $\Mpaircd$-module decomposition:
 \begin{eqnarray*}
   \Fminusl = \bigoplus_{\lambda \in \Sigma (C)} 
          V( Sp(2l); \lambda) \otimes L \left(\dinf;
                                 \Lambda^{ \frak{cd}} (\lambda), -l
                               \right)
 \end{eqnarray*}
 where $V( Sp(2l); \lambda)$ is the irreducible highest
 weight $Sp (2l)$-module
 of highest weight $\lambda$ and 
 $L \left(\dinf; \Lambda^{ \frak{cd}} (\lambda), -l  \right)$ 
 the irreducible highest weight $\dinf$-module of highest weight
 $\Lambda^{ \frak{cd}} (\lambda) $
 and  central charge $-l$.
  \label{th_Bpair}
\end{theorem}

The action of $\hBp$ on $\Fminusl$ is given by the composition
of the action of $\dinf$ and the homomorphism $\hphic$. 
A similar argument as in Lemma \ref{lem_ffr} shows that 
\begin{eqnarray}
 W^n (z) = \hf \sum_{p =1}^l \left(
     :\partial_z^n\gamma^{+,p} (z) \gamma^{-,p} (z):
     - :\gamma^{+,p}  (z) \partial_z^n \gamma^{-,p} (z):
     \right).
  \label{eq_cl}
\end{eqnarray}
Given $\lambda = (m_1, \ldots, m_l) \in \Sigma (C)$, 
we let $e (\lambda)$ be the set of 
exponents $k$ with multiplicity $m_k - m_{k +1}$ $(k =1, \ldots, l)$
where $m_{l +1} = 0$
(the exponent $0$ has multiplicity $-l -m_1$).
\begin{theorem}
   We have the following $\Wpaircd$-module decomposition:
 \begin{eqnarray*}
   \Fminusl = \bigoplus_{\lambda \in \Sigma (C)} 
          V( Sp(2l); \lambda) \otimes L \left(\hBp;
                                 e (\lambda), -l
                               \right).
 \end{eqnarray*}
   \label{th_Wpaircd}
\end{theorem}
\begin{demo}{Proof}
    By Theorem \ref{th_pullback} the $\dinf$-module
 $L \left(\dinf; \Lambda^{\frak{cd}} (\lambda), -l  \right)$
 regarded as a $\hBp$-module via the pullback by $\hphio$
 remains irreducible. To prove the thereom it suffices
 to determine the corresponding $\Delta (x)$ for this $\hBp$-module. 
 It follows from the definition of $ \Lambda^{\frak{cd}} (\lambda)$
 by using Proposition~\ref{prop_deltad} that 
 for $\lambda = (m_1, \ldots, m_l) \in \Sigma(C) $
 \begin{eqnarray*}
 \Delta (x) & = & 
      \frac{ \sum_{k =1}^l (m_k - m_{k +1})\cosh (kx)
            + \hf ( (-2l -m_1 -m_2) - (m_1 -m_2) ) - ( -l) }
                {2 \sinh (x/2)}  \\
            & = &  \frac{ \sum_{k =1}^l (m_k - m_{k +1})\cosh (kx)
            + (-l -m_1)  }
                {2 \sinh (x/2)} - \frac{( -l)} {2 \sinh (x/2)}.
 \end{eqnarray*}
\end{demo}

The following corollary is immediate.
\begin{corollary} 
    The space of invariants of $Sp (2l)$ in the Fock space
  $\Fminusl$ is naturally isomorphic to the irreducible
  module $L(\dinf; -2l {}^d \hL_0 )$ or equivalently 
   to the irreducible $\hBp$-module 
 $L \left(\hBp;e( 0), -l  \right)$.
   \label{cor_cl}
\end{corollary}
\subsection{Case $ c \in - \Bbb N + 1/2$}
      We denote by $\Fhalfminusl$ the tensor product
of the Fock space $\Fminusl$ of $l$ pairs of bosonic ghost fields
and the Fock space $\Fhalf$ of a neutral fermionic field $\phi (z)$.
It is known \cite{FF} that the Fourier components of 
the generating functions $e^{pq}(z), e^{pq}_* (z),$ 
$ e^{pq}_{**}(z)$
in (\ref{eq_genfunc}) and the following generating functions
\begin{eqnarray}
 e^{p} (z) \equiv \sum_{n \in \Z}
 e^{p}(n) \NN & = & : \gamma^{-,p}(z) \phi (z):
         \nonumber       \\
 \tilde{e}_*^{p} (z) \equiv \sum_{n \in \Z}
 \tilde{e}_*^{p}(n) \NN & = & : \gamma^{+,p}(z) \phi (z):
    \quad (p = 1, \dots, l )   \nonumber 
  \label{eq_oddsuper}
\end{eqnarray}
span the affine superalgebra 
$\widehat{ {\frak {osp}} }(1, 2l)$ of level $-1$. 
The operators $ e^{pq}(0),$ $e_*^{pq}(0),$ $e_{**}^{pq}(0),$
$\tilde{e}^{p}(0), \tilde{e}_*^{p}(0)$ $ (p, q = 1, \dots, l) $ 
generate the horizontal subalgebra
${\frak {osp}}(1,2l)$ of the affine superalgebra
$\widehat{{\frak {osp}} }(1, 2l) $.
We identify the Borel subalgebra ${\frak b} ({\frak {osp} }(1, 2l) )$
with the one generated by 
$  e_*^{pq}(0) (p \leq q), e_{**}^{pq}(0), \tilde{e}_*^{p}(0),$ 
$ p, q = 1, \cdots, l. $
Introduce the generating function
 \begin{equation}
      \begin{array}{rcl}
      \lefteqn{
  \sum_{i,j \in \Z}  \left( E_{ij}  - E_{1-j,1-i}
                     \right) z^{i- 1} w^{-j}
              }                                \\
       & &
  = \sum^l_{k=1} \left(: \gamma^{+,k} (z) \gamma^{-,k} (w) :
                     - : \gamma^{+,k} (w) \gamma^{-,k} (z) :
                 \right) 
                     + :\phi (z) \phi (w):.
      \end{array}
 \end{equation}
This defines a representation of $\dinf$ of central
  charge $-l + \hf$ on $\Fhalfminusl$. It is known \cite{W2} the action of
the horizontal subalgebra ${\frak {osp}}(1, 2l)$ commute with 
that of Lie algebra $\dinf$ on $\Fhalfminusl$.
$({\frak {osp}}(1, 2l), \dinf)$ form a dual pair on $\Fhalfminusl$.

We define a map $\Lambda^{ \frak{ospd}}$
from $\Sigma (osp) $ to ${\dinf}^*_0$ by sending
$ \lambda = (m_1, \ldots, m_l )$ to
$$ \Lambda^{ \frak{ospd}} (\lambda )=
   (-2l +1 -m_1 -m_2) {}^d \hL_0 + \sum_{k =1}^l (m_k - m_{k +1} ){}^d \hL_k.$$
The following duality theorem was quoted from \cite{W2}.
\begin{theorem}
  We have the following $\Mpairospd$-module decomposition
   \begin{eqnarray*}
    \Fhalfminusl = \bigoplus_{\lambda \in \Sigma (Osp)} 
          V( {\frak {osp}}(1, 2l); \lambda) \otimes L \left(\dinf;
                             \Lambda^{ \frak{ospd}} (\lambda), -l +1/2 
                               \right)
   \end{eqnarray*}
 where $V({\frak {osp}}(1, 2l); \lambda)$ is the irreducible
 module of ${\frak {osp}}(1, 2l) (1, 2l)$ of highest weight $\lambda$, and 
 $L \left(\dinf; \Lambda^{ \frak{ospd}} (\lambda), -l +1/2 \right)$
 is the irreducible highest weight $\dinf$-module of highest weight
 $\Lambda^ \frak{ospd} (\lambda) $ and central charge $-l +1/2$.
   \label{th_pairospd}
\end{theorem}

We define the action of $\hBp$ on $\Fhalfminusl$ to be the
composition of the action of $\dinf$ and the homomorphism
$\hphio$ given by (\ref{eq_embd}). In a similar way as obtaining
(\ref{eq_cl}) we have
\begin{eqnarray}
 W^n (z) & = & \hf \sum^l_{k=1} 
    \left(: \gamma^{+,k} (z) \partial_z^n\gamma^{-,k} (z) :
                     - : \gamma^{-,k} (z) \partial_z^n \gamma^{+,k} (z) :
                 \right)    \nonumber  \\
         & &      + \hf :\phi (z)  \partial_z^n \phi (z):.
   \label{eq_osp}
\end{eqnarray}
Then of course the action of the horizontal subalgebra
${\frak {osp}}(1, 2l)$ commute with 
that of Lie algebra $\hBp$ on $\Fhalfminusl$.

Given $\lambda = (m_1, \ldots, m_l) \in \Sigma (Osp)$, 
we let $e (\lambda) $ be the set of
exponents $k$ of multiplicity $m_k - m_{k +1}$ $(k =1, \ldots, l)$
where $m_{l +1} = 0$ (the exponent $0$ has multiplicity $-l -m_1 +\hf$).
We obtain the following theorem in a similar way as Theorem \ref{th_Wpaircd}
by using Theorems~\ref{th_pullback}, \ref{th_pairospd}, and
Proposition~\ref{prop_deltad}.

\begin{theorem}
  We have the following $\Wpairospd$-module decomposition
   \begin{eqnarray*}
    \Fhalfminusl = \bigoplus_{\lambda \in \Sigma (Osp)} 
          V( {\frak {osp}}(1, 2l); \lambda) \otimes L \left(\hBp;
                             e (\lambda), -l +1/2 
                               \right).
   \end{eqnarray*}
   \label{th_Wpairospd}
\end{theorem}

  The following corollary is immediate.
\begin{corollary}
    The space of invariants of ${\frak {osp}}(1, 2l)$ in the Fock space
  $\Fhalfminusl$ is naturally isomorphic to the irreducible
  module $L(\dinf; (-2l +1) {}^d \hL_0 )$ or equivalently 
  to the irreducible $\hBp$-module 
 $L \left(\hBp; 0, -l+ 1/2  \right)$.
   \label{cor_osp}
\end{corollary}
\section{Vertex algebra associated to $\hBp$}
  \label{sec_voa}
\subsection{Vertex algebra structure on the vacuum module of $\hBp$}
  In Example \ref{ex_para2} we have constructed the vacuum $\hBp$-modules
$M_c$ and $V_c$ $ ( c \in \C)$. We want to show that $M_c$ and $V_c$
carry a natural structure of a vertex algebra, cf.  
\cite{B, FLM, DL, K2}. 

 Denote by $\vac$ the highest weight vector of $M_c$ and $V_c$.
$M_c$ (resp. $V_c$) is spanned by the vectors
\begin{equation}
  W^{n_j}_{-k_j -n_j -1} \cdots W^{n_1}_{-k_1 -n_1 -1} \vac,\;
  n_i \in 2 \Z_{+} + 1, \; k_i \in \Z .
   \label{eq_vector}
\end{equation}

By the definition of $M_c$ the vector $\vac$ is annihilated by
$\cal P$, namely $W^n_k \vac = 0$ for $k +n \geq 0$.
Operators from $\cal P$ is often referred to as annihilation operators 
and the operators $W^n_k $ from $\hBp - \cal P$ as creation operators.
We have known that $W_k^1 = - t^k (D + (k +1)/2 )\; (k \in \Z)$ 
form the Virasoro algebra
with central charge $c$ when acting on $M_c$ or $V_c$.
Also it follows by a direct computation using (\ref{eq_12}) that
$$
  \left[ W^1_{-1}, W^n_k \right] = -(n +k) W^n_{k -1}
$$
or equivalently
\begin{eqnarray}
   \left[ W^1_{-1}, W^n (z) \right] = \partial_z W^n (z).
 \label{eq_transl}
\end{eqnarray}
Denote $T = W^1_{-1}$. Define a linear map 
$Y(\cdot, z) : M_c \longrightarrow (\End M_c) [[z, z^{-1} ]]$ 
by 
$$ Y(W^n_{-n -1} \vac, z) = W^n (z) \quad (n \in 2 \Z_{+} + 1)$$
and in general we associate to a vector $a$ of the form (\ref{eq_vector})
the following generating function
$$
 Y(a, z)
  = :\partial_z^{(k_j)} W^{n_j}(z) \cdots \partial_z^{(k_1)} W^{n_1}(z):
$$
where $\partial_z^{(k)}$ denotes $\partial_z^k / k!$ and
$::$ denotes the standard normal ordering from the right to left which
moves the annihilators to the right.
Clearly $ Y(a, z) \vac \mid_{z =0} = a.$
\begin{proposition}
 $ (V_l, \vac, T, Y(\cdot, z) )$ is a vertex 
 algebra isomorphic to the space of invariants 
  of $\Fl$ for $\underline{\Z} = \hz$
 with respect to $O(2l)$. The generating
 fields $W^n (z)$ are given by formula (\ref{eq_field}).
    \label{prop_voa}
\end{proposition}
    
 By Corollary \ref{cor_dl}, $V_l$ and the space of invariants 
of $\Fl$ with respect to $O(2l)$ are naturally isomorphic 
as $\hBp$-modules. Since $\Fl$
has a natural vertex algebra structure, the vertex algebra
structure on $V_l$ is ensured by
the well-known fact that the space of invariants of an
automorphism group of a vertex algebra is a vertex algebra.
We will give a direct proof in Lemma \ref{lem_closeope} below
that the generating fields $W^n (z)$ $( n \in 2\Z_{+} +1)$
are closed under operator product expansions. 

Consider the case of $\Fl$ with $\underline{\Z} = \hz$. Denote 
$$ \Psi^{m,n} (z) = \sum_{k=1}^l
  \left(  : \partial^m \psi^{+,k} (z) \partial^n \psi^{-,k} (z):
        + : \partial^m \psi^{-,k} (z) \partial^n \psi^{+,k} (z):
  \right),
 \quad m,n \in \Z_+.
$$
Note the obvious symmetry $\Psi^{m,n} (z) = - \Psi^{n, m} (z)$.
\begin{lemma}
   $ \Psi^{m,n} (z)$ is a linear combination of
 $ \partial^{i} W^{m +n -i} (z), 0 \leq i \leq  m +n$
 and $i \equiv m +n -1 \bmod{2}$.
   \label{lem_comb}
\end{lemma}
\begin{demo}{Proof}
 We prove by induction on $m +n$.
 When $m+n = 0$ or $ 1,$ the statement is obvious.
 
 Assume that for $m +n = 2k -1 $ $( k \in \Bbb N)$ the statement is true.
 Then by this induction assumption,
 $\partial \Psi^{2k -1 -m, m} (z)$
 $(m = 0, \ldots, k-1 )$ is a linear combination of
 $ \partial^{i} W^{2k +2 -i} (z)$ ($ 0 \leq i \leq 2k -1$ and $ i$ odd).
 Since 
 $$ \partial \Psi^{2k -1 -m, m} (z)
   = \Psi^{2k -m, m} (z) + \Psi^{2k -1 -m, m +1} (z)
 $$
  it follows by a little algebra that the linear span of
 $\Psi^{2k -m, m} (z), m =0, \ldots, k -1$
 is equal to the linear span of 
 $\partial \Psi^{2k -1 -m, m} (z), m =0, \ldots, k -1.$
 This  proves that the statement is true for $m +n = 2k $.

 Therefore it follows that $\partial \Psi^{2k -m, m} (z),
 m = 0, \ldots, k-1$ is a linear combination of
 $ \partial^{i} W^{2k +1 -i} (z)$ ($ 0 \leq i \leq 2k +1$ and $i$ even).
 Since 
 $$ \partial  \Psi^{2k -m, m} (z)
   = \Psi^{2k +1 -m, m } (z) +  \Psi^{2k +1 -m, m} (z)
 $$
 by a little linear algebra it is easy to show that the linear span of
 $ \Psi^{2k +1 -m, m} (z)$ $( m = 0, \ldots, k)$
 is equal to the linear span of 
 $\partial \Psi^{2k -m, m} (z)$ $( m = 0, \ldots, k -1)$ and $ W^{2k +1} (z)$.
 This proves that the statement is true for $m +n = 2k +1$.
\end{demo}
\begin{lemma}
  $ \partial^i W^j (z), i \in \Z_{+}, j \in 2 \Z_{+} + 1$
 are closed under the operator product expansions. 
   \label{lem_closeope}
\end{lemma}
\begin{demo}{Proof}
 For the simplicity of notations, we first assume $l =1$ and
 write $\psi^{\pm} (z)$ for $\psi^{\pm, 1} (z)$.  
   By Wick's theorem, we have
 \begin{eqnarray*}
  W^m (z) W^n (z)
 & \sim & \frac14 \left(
             :\partial_z^m \psi^- (z) \psi^+ (z): + 
             :\partial_z^m \psi^+ (z) \psi^- (z):
          \right) \cdot      \\
  & &     \left(
             :\partial_w^n \psi^- (w) \psi^+ (w): + 
             :\partial_w^n \psi^+ (w) \psi^- (w):
          \right)       \\
 & \sim & \frac14 \left( \partial_z^m ((z-w)^{ -1}) 
               \right)   
               \left( :\psi^+ (z) \partial_w^n \psi^- (w): +
                      :\psi^- (z) \partial_w^n \psi^+ (w): 
               \right)     \\
  & & + \frac14 \left( \partial_w^n ( (z-w)^{ -1}) 
               \right) 
               \left( :\partial_z^m \psi^- (z) \psi^+ (w): +
                      :\partial_z^m \psi^+ (z) \psi^- (w):
               \right)  \\
  & & + \frac12   \partial_z^m ((z-w)^{ -1}) \partial_w^n ( (z-w)^{ -1})  \\
  & \sim &  \frac{ ( -1)^m m!}{ 4 (z -w)^{m +1} }               %
               \left( :\psi^+ (z) \partial_w^n \psi^- (w): +
                      :\psi^- (z) \partial_w^n \psi^+ (w): 
               \right)     \\
  & & +  \frac{ n!}{ 4 (z -w)^{n +1} }       
               \left( :\partial_z^m \psi^- (z) \psi^+ (w): +
                      :\partial_z^m \psi^+ (z) \psi^- (w):
               \right)  \\
  & & + \frac{( -1)^m (m+n)!}{2 (z -w)^{m +n +2}}.
\end{eqnarray*}
 
 By taking the Taylor expansions of $\psi^{\pm} (z)$ and
 $\partial_z^m\psi^{\pm} (z)$ at $ z =w$, we can easily see
 that all the fields appearing on the right hand side of the 
 above OPE are linear combinations of $ \Psi^{m,n} (z)$.
 Thus our lemma follows from Lemma \ref{lem_comb}.
 For the general case $\Fl$, it is clear that
 the only modification in the final OPE formula above is 
 that the central term should be mutiplied by $l$.
\end{demo}
     
The following corollary holds in the cases $c =l$ $( l \in \Bbb N )$
by the computation of OPE in Lemma \ref{lem_closeope} since
the maximal order of poles at $z = w$ appearing there is $m +n +2$.
Since the central term depends on the central charge $c$ linearly,
the corollary remains true for an arbitrary central charge.
\begin{corollary}
  For arbitrary central charge $c$ we have  as a formal power series 
\begin{eqnarray*}
  (z -w)^{m +n +2} \left[ W^m (z), W^n (w) \right] = 0.
\end{eqnarray*}
  \label{cor_locality}
\end{corollary}
       
Now the following theorem follows by a general argument
in the theory of vertex algebras, cf. \cite{L}, 
Proposition 3.1 in \cite{FKRW},
or Theorem 4.5 in \cite{K2}, since all the requirements there are satisfied
(also cf. \cite{G}).
\begin{theorem}
 $  (M_c, \vac, T, Y(\cdot, z) )$ and $(V_c, \vac, T, Y(\cdot, z) )$
 are vertex algebras.
  \label{th_vacuumvoa}
\end{theorem}
 
The following proposition follows from Corollaries \ref{cor_bl},
\ref{cor_cl}, \ref{cor_osp} and Theorem~\ref{th_vacuumvoa}.
\begin{proposition}
   \begin{enumerate}
   \item[1)] $ (V_{l +1/2}, \vac, \omega, Y(\cdot, z) )$ is a vertex 
 algebra isomorphic to the space of invariants 
  of $\Flhalf$ for $\underline{\Z} = \hz$
 with respect to $O(2l +1)$. The generating
 fields $W^n (z)$ are given by formula (\ref{eq_bl}).

   \item[2)] $ (V_{-l }, \vac, \omega, Y(\cdot, z) )$ is a vertex 
 algebra isomorphic to the space of invariants 
  of $\Fminusl$ with respect to $Sp(2l)$. The generating
 fields $W^n (z)$ are given by formula (\ref{eq_cl}).

   \item[3)] $ (V_{-l +1/2}, \vac, \omega, Y(\cdot, z) )$ is a vertex 
 algebra isomorphic to the space of invariants 
  of $\Fhalfminusl$ with respect to ${\frak {osp}}(1, 2l)$. The generating
 fields $W^n (z)$ are given by formula (\ref{eq_osp}).
  \end{enumerate}
   \label{prop_vary}
\end{proposition}
\begin{corollary}
   The irreducible representations of $\hBp$ 
 appearing in the decompositions of Fock spaces
 $\Fl, \Flhalf, \Fminusl$ and $\Fhalfminusl$ are representations of
 the vertex algebra associated to $\hBp$ with central charges
 $l, l + 1/2, -l$ and $ -l + 1/2$ respectively.
\end{corollary}
\begin{remark} \rm
 It is interesting
 to compare some 
 theorems in \cite{DLM} by Dong-Li-Mason 
 with our Theorems \ref{th_Wpairdd}, 
 \ref{th_Wpairbd}, \ref{th_Wpaircd} and \ref{th_Wpairospd} in the
 framework of vertex algebras. Their results have the virtue of
 being general while the statements of our theorems 
 are more powerful and precise thanks to being concrete.
\end{remark}
\subsection{Vertex algebra $V_c$ for $c \in \hf \Z$ 
and $\cal W$ algebras}
   A particular important class of vertex algebras is the so-called
$\cal W$-algebras. One can associate a $\cal W$-algebra ${\cal W} {\frak g}$ 
to an arbitrary complex simple Lie algebra $\frak g$ (see \cite{BS, FeF}
and references therein). In general such a $\cal W$-algebra can
be defined in terms of intersections of screening operators,
cf. e.g. \cite{FeF}. Denote by ${\cal W}(\widehat{ \frak g} /\frak g, k)$
the $\cal W$-algebra arising as the space of invariants of 
$\frak g$ in the vacuum module of the affine algebra $\widehat{ \frak g}$
of level $k$, where $\frak g$ is the horizontal 
subalgebra of $\widehat{ \frak g}$.
In the case when $\frak g$ is simply-laced of
rank $l$, ${\cal W} {\frak g}$ with central charge $l$
can be shown to be isomorphic to ${\cal W}(\widehat{ \frak g} /\frak g, 1)$
cf. e.g. \cite{BS, F2, BBSS, FKRW}. 

It is known \cite{F1, KP} that 
the basic representation of $\widehat{ {\frak {so}} }( 2l)$
is isomorphic to the even part of the vector superspace $\Fl$.
From Remark~\ref{rem_semi} we see that the 
space of invariants of $\Fl$ with respect to $SO(2l)$ is isomorphic to 
the $\hBp$-module $L \left(\hBp; 0, l \right)
 \bigoplus L \left(\hBp; e (det), l \right)$ (see Section \ref{sec_n}).
The highest weight vector of $L \left(\hBp; e (det), l \right)$
is given by 
$\prod_{k =1}^l (\psi^{+, k}_{ - 1/2}\psi^{-, k}_{ - 1/2}) \vac$ \cite{W2} 
which lies in the even part of $\Fl$. Thus 
${\cal W}(D_l^{(1)} /D_l, 1)$
is isomorphic to the space of invariants of $\Fl$ with respect to $SO(2l)$.
In this way we have reached the following conclusion.
\begin{theorem}
 The $\cal W$-algebra ${\cal WD}_l $ with central
 charge $l$ is isomorphic to
 a sum of the vertex algebra $V_l$ and the irreducible
 module $L \left(\hBp; e (det), l \right)$ of $V_l$.
   \label{th_wdl}
\end{theorem}

 As in \cite{FKRW}, Theorem \ref{th_wdl} implies the following
corollary.
\begin{corollary}
   All positive primitive modules over $\hBp$ with a positive integral
 central charge $l$ are irreducible modules over the vertex algebra $V_l$. 
\end{corollary}
\begin{remark} \rm
  It is an important question to determine
 whether these positive primitive modules
 are all irreducible modules over the vertex algebra $V_l$.
\end{remark}
\begin{remark} \rm
 Theorem \ref{th_wdl} provides a new way to compute the $q$-character
 of ${\cal WD}_l $ with central charge $l$:
 \begin{eqnarray}
  \ch_q {\cal WD}_l
   & = &   \ch_q L \left(\hBp; 0, l \right)
     + q^l \ch_q L \left(\hBp; e (det), l \right)    \nonumber \\
   & = &    \ch_q L \left( \dinf, 2l \, {}^d \hL_0 \right)
      + q^l \ch_q L \left( \dinf, 2l \, {}^d \hL_1 \right).
   \label{eq_newch}
 \end{eqnarray}
 Here $q^l$ accounts for the weight $l$ of the highest weight
 vector $\prod_{k =1}^l (\psi^{+, k}_{ - 1/2}\psi^{-, k}_{ - 1/2}) \vac$
 of $L \left(\hBp; e (det), l \right)$ in $\Fl$. The $q$-character
 formulas of $\dinf$-modules $L \left( \dinf, 2l \, {}^d \hL_0 \right)$
 and $L \left( \dinf, 2l \, {}^d \hL_1 \right)$ can be read from the Appendix.
 A straightforward computation by using (\ref{eq_newch}) yields
 \begin{eqnarray*}
   \ch_q {\cal WD}_l
   & = & \prod_{1 \leq i < j \leq l} (1 -q^{j -i}) (1 - q^{i +j -2}) /
     \varphi (q)^l    \\
   & = & \frac{ \prod_{i =1}^l  ( \prod_{n =1}^{e_i} (1 -q^n) )}
    {\varphi (q)^l}
 \end{eqnarray*}
 where $\varphi (q) = \prod_{i \geq 1} (1 -q^i)$
 and $e_i = 2i -1$ $(i =1, \cdots, l -1)$, $e_l = l -1$ 
 are the exponents of the simple Lie algebra
 $\frak {so}( 2l)$. The same formula was earlier
 deduced in \cite{K0} and \cite{BS}. The character formula implies
 that ${\cal WD}_l$ with central charge $l$ is freely generated by fields
 $W^i (z) (i =1, 3, \ldots, 2l -3)$ of conformal weights $i +1$ and
 the field $\prod_{i =1}^l: \psi^{+, k} (z) \psi^{-, k} (z):$
 of conformal weight $l$ corresponding to the vector 
 $\prod_{k =1}^l (\psi^{+, k}_{ - 1/2}\psi^{-, k}_{ - 1/2}) \vac$.
\end{remark}
\begin{remark} \rm
  Combining Corollary \ref{cor_dl} with Theorem \ref{th_wdl},
 we obtain another dual pair $(SO(2l), {\cal WD}_l)$ on the
 Fock space $\Fl$. When restricting to
 the even part of $\Fl$ which is isomorphic
 to the basic representation of $\widehat{ {\frak {so}}}(2l)$,
 we recover a duality theorem of I. Frenkel \cite{F2}.
 Since the language of vertex algebras was not
 available at the time of \cite{F2}, $ {\cal WD}_l$ was replaced
 by the Lie algebra of Fourier components of fields in $ {\cal WD}_l$.
\end{remark}
The even part of the fermionic Fock space
$\Flhalf$ for $\underline{\Z} = \hz$ is isomorphic to
the basic representation of $\widehat{ {\frak {so}} }( 2l +1)$ \cite{F1, KP}.
According to Corollary \ref{cor_bl}, the space of invariants
of $\Flhalf$ with respect to $SO ( 2l +1)$ is isomorphic to
the $\hBp$-module $L (\hBp; 0, l+1/2 )
  \bigoplus L (\hBp; e (det), l+1/2 )$. 
The highest weight vector of 
$L (\hBp; e (det), l+1/2 )$ in $\Flhalf$ is 
$\prod_{k =1}^l (\psi^{+, k}_{ - 1/2}\psi^{-, k}_{ - 1/2}) \phi_{ -1/2} \vac$,
cf. \cite{W2}. Observe that it is an odd vector in $\Flhalf$.
Thus  ${\cal W}(B^{(1)}_l /B_l, 1)$
is isomorphic to $L (\hBp; 0, l+1/2 )$. 
Combining with Proposition \ref{prop_vary} we have proved
the following theorem.
\begin{theorem}
 The vertex algebra $V_{l +1/2}$ associated to $\hBp$ with 
 $c = l +1/2$ is isomorphic to the $\cal W$-algebra
  ${\cal W}(B^{(1)}_l /B_l, 1)$.
   \label{th_wbl}
\end{theorem}
\begin{corollary}
   All positive primitive modules over $\hBp$ with a 
 central charge $l + 1/2 \in 1/2 + \Z_+$ are 
 irreducible modules over the vertex algebra $V_{ l + 1/2}$. 
\end{corollary}
\begin{remark} \rm
  It remains to determine
 whether these positive primitive modules
 are all irreducible modules over the vertex algebra $V_{ l + 1/2}$.
\end{remark}
   
  It follows that the character formula of 
the $\cal W$-algebra ${\cal W}(B^{(1)}_l /B_l, 1)$
is the same as the $q$-character formula
of $L \left(\dinf; (2l +1) {}^d \hL_0 \right)$ (cf. Appendix):
$$\frac{ \prod_{i =1}^l  ( \prod_{n =1}^{e_i} (1 -q^n) )}
    {\varphi (q)^l} \cdot
   \hf \left( {\prod_{n \geq 1} ( 1 +q^{ n+ l - 1/2})}
            + {\prod_{n \geq 1} ( 1 -q^{ n+ l - 1/2})}
       \right).
$$
where $ e_i = 2i -1$ $(i =1, \cdots, l)$ are the exponents of the
Lie algebra $B_l$.
A different method was used in \cite{BS}
to obtain the same formula. In particular we have seen that
${\cal W}(B^{(1)}_l /B_l, 1)$ lies inside a vertex superalgebra,
denoted by ${\cal WB} (0, l)$, of central charge $l -1/2$
which is isomorphic to 
the space of invariants of $\Flhalf$ with respect to $SO ( 2l +1)$.
The $q$-character of the latter can be calculated by using the
$q$-character formulas of $\dinf$-modules given in the Appendix as
\begin{eqnarray}
  \ch_q L \left(\dinf; (2l +1) {}^d \hL_0 \right)
      + q^{l +1/2} \ch_q L \left(\dinf; (2l +1) {}^d \hL_1 \right) \nonumber \\
     = \frac{ \prod_{i =1}^l  ( \prod_{n =1}^{e_i} (1 -q^n) )}
    {\varphi (q)^l} {\prod_{n \geq 1} ( 1 +q^{ n+ l - 1/2})}.
   \label{eq_wb}
\end{eqnarray}
Here $q^{l +1/2}$ accounts for the grading of the highest weight
vector of the $\dinf$-module
$L \left(\dinf; (2l +1) {}^d \hL_1 \right)$ in $\Flhalf$.

Formula (\ref{eq_wb}) indicates that
vertex superalgebra ${\cal WB} (0, l)$ 
is freely generated by the bosonic fields
$W^i (z) (i =1, 3, \ldots, 2l-1)$ of conformal weights $i +1$ and
the fermionic field
$\prod_{i =1}^l: \psi^{+, k} (z) \psi^{-, k} (z) \varphi (z):$
of conformal weight $l +1/2$ corresponding to the vector 
$\prod_{k =1}^l (\psi^{+, k}_{ - 1/2}\psi^{-, k}_{ - 1/2})
   \varphi_{ - 1/2}) \vac$. 
The phenomenon that ${\cal W}( B^{(1)}_l /B_l, 1)$
lies in such a vertex superalgebra was conjectured in \cite{BS}.
Here we find the explicit model for ${\cal WB} (0, l)$ with 
central charge $l -1/2$. 
${\cal WB} (0, l)$ were also obtained by the quantized Drinfeld-Sokolov
reduction in \cite{Ito}. We summarize these in the following theorem
which resembles remarkably Theorem~\ref{th_wdl}.
\begin{theorem}
  The vertex superalgebra ${\cal WB} (0, l)$ 
 is the sum of the vertex algebra $V_{l +1/2}$ 
 and the irreducible module
 $L \left(\hBp;  e (det), l+1/2\right)$ of $V_{l +1/2}$.
   \label{th_wsuperbl}
\end{theorem}

The following theorem 
easily follows from Proposition \ref{prop_vary} and similar argument
which leads to Theorems~\ref{th_wbl}.
\begin{theorem}
 The vertex algebra $V_{ -l}$ associated to $\hBp$ with central
 charge $ -l$ is isomorphic to the $\cal W$-algebra
 ${\cal W}(C_l^{(1)} /C_l, -1)$.
 The vertex algebra $V_{ -l + 1/2}$ associated to $\hBp$ with central
 charge $ -l +1/2$ is isomorphic to 
 ${\cal W}(B^{(1)}(0, l) /B(0, l), -1)$. 
\end{theorem}
\begin{remark} \rm
 Theorems \ref{th_wdl} and \ref{th_wsuperbl} tell us
 what is the minimal set of generating fields 
 for $V_c$ with $c \in  \hf \Bbb N$ (or rather for the vertex (super)algebra
 which are the $\Z_2$ extension of $V_c$).
 The question remains open for  $V_c$ for a general negative half-ingeral
 central charge. A similar
 question was addressed in the case of $\W$ in \cite{W} (also
 see \cite{EFH}). The negative central charge case
 turned out to be more subtle and difficult.
\end{remark}
\section{Appendix}
   The following proposition gives the $q$-character formula
(compatible with the $\Z$-gradation of $\dinf$ induced from that of $\hgl$)
of a highest weight representation of $\dinf$ with highest
weight $\Lambda = {}^d \hL_{n_1} + {}^d \hL_{n_2} + \ldots +
 {}^d \hL_{n_k} + {}^d h_0 {}^d \hL_0 $, where 
$n_1 \geq \ldots n_k \geq 1$, ${}^d h_0 \in \Z_{+}$. 
Denote ${}^d h_1 $ to be the number of $n_i$'s which are equal to $1$.
Then the central charge $c = \hf ( {}^d h_0 - {}^d h_1 ) + k$. 
We retain the notation in Section \ref{sec_classical}.
\begin{proposition}
  The $q$-character formula of $L( \dinf; \Lambda)$ corresponding
 to the principal gradation of $\dinf$ is
\begin{eqnarray}
   \ch_q L(\dinf; \Lambda)
  & = & \frac{ \prod\limits_{1 \leq i < j \leq k} (1 - q^{n_i - n_j +j -i})}
             { \prod\limits_{1 \leq i \leq k} \varphi_{n_i + k -i}(q)} \cdot
 \frac{\varphi(q^2)^{\overline{2c}}
                         \prod\limits_{j >0} \varphi_{2c - 2j}(q)} 
              {\varphi(q)^{ [c] - \overline{2c} } } \cdot     \nonumber   \\
   & \quad & \prod_{i=0}^{n_1-1}
               \frac{ \varphi_{2c + i + n_1}(q)}
                    { \varphi_{2c + n_1 + i - {}^d \lambda_{i +1}} (q) } \cdot
          \prod_{0 \leq i <j \leq n_1} 
           \frac{ 1 - q^{2c +j +i - {}^d \lambda_{i +1} - {}^d \lambda_{j+1}}}
                {1-q^{2c +j +i}} \cdot   \nonumber   \\
   & &\cdot \hf
        \left( \prod_{j \in \Bbb N} 
                 \left( 1 + q^{- ( \Lambda + \rho, \epsilon_j)}
                 \right)
                   + \prod_{j \in \Bbb N} 
                 \left( 1 - q^{- (\Lambda + \rho, \epsilon_j)}
                 \right)
       \right) .     
  \label{eq_Dqchar}
\end{eqnarray}
  \label{prop_dform}
\end{proposition}
\begin{demo}{Proof}
 Note that the $q$-characters of $\dinf$-modules are specialized
 characters of type $(2, 1, 1, \ldots)$, in contrast to the
 $q$-characters of modules of $\binf, \cinf$ which are of
 type $(1, 1, \ldots)$. So the way of computing such a $q$-character
 of a $\dinf$-module  will be different. We will proceed as follows.
 By the Weyl-Kac character formula
 \begin{eqnarray*}
  e(-\Lambda ) \ch L( \Lambda)  
    \equiv \frac{\cal N}{\cal D} = \frac
         {\sum_{w \in W} \epsilon (w) e( w( \Lambda + \rho)-(\Lambda +\rho))}
         {\prod _{\alpha \in \Delta_+}(1-e(-\alpha))}.
 \end{eqnarray*}
 Here $\cal N$ (resp. $\cal D$) denotes the numerator (resp. denominator).
 Given $ \vec s = (2, 1, 1,\ldots)$, 
 define a homomorphism $F_{\vec s}: \Bbb C [[ e(-\alpha_i)]] \longrightarrow
 \Bbb C [[q]]$ by $F_{\vec s} (e(-\alpha_i) = q^{s_i}$.
 Put  $ \vec t =(t_0, t_1, t_2, \ldots)$
 where $t_i = \langle \Lambda + \rho, \alpha_i^{\vee} \rangle$. Then

 \begin{eqnarray}
  F_{\vec s}  (\cal N)
   &  = & \sum_{ w \in W} \epsilon (w)
                q^{\langle \Lambda + \rho - w( \Lambda + \rho),
                      \hL_0^{\vee} + \rho^{\vee}  \rangle}  \nonumber  \\
   &  = & q^{\langle \Lambda + \rho ,
               \hL_0^{\vee} + \rho^{\vee} \rangle}
                 \sum_{w \in W} \epsilon (w) q^{\langle \Lambda + \rho,
                       -w( \hL_0^{\vee} + \rho^{\vee} )\rangle}  \nonumber  \\
   & = &  q^{ \langle \Lambda + \rho, \hL^{\vee}_0 + \rho^{\vee} \rangle} 
      F_{\vec t} \left( \sum_{w \in W} \epsilon (w)
                         e( w( \hL_0 +\rho )) \right) 
              \nonumber \\
   & = &  q^{\langle \Lambda + \rho, \hL^{\vee}_0 + \rho^{\vee} \rangle}
      F_{\vec t} \left( \ch L( \hL_0) \cdot \sum \epsilon (w)
               e( w( \rho) - \rho) \cdot e(\rho) \right)  \nonumber \\
   & = & F_{\vec t} (e( -\rho -\hL_0)) F_{\vec t} \left( \ch L( \hL_0) \right)
       F_{\vec t} \left( \prod_{\alpha \in \Delta_+} (1 - e( -\alpha) \right)
       F_{\vec t} (e(\rho))     \nonumber    \\
   & = &  F_{\vec t} \left( e(- \hL_0)\ch L( \hL_0) \right) \cdot
       F_{\vec t} \left( \prod_{\alpha \in \Delta_+}(1 -e( -\alpha)) \right).
    \label{eq_iden}
 \end{eqnarray}
By Theorem \ref{th_Fpairbbtwo} (putting $l = 0$ there), the even part
of ${\cal F}^{\bigotimes \hf}$ is isomorphic to the irreducible $\dinf$-module
 with highest weight $\hL_0$ and central charge $\hf$, so
 \begin{eqnarray}
  F_{\vec s} ( e( -\hL_0) \ch L( \hL_0))
   & = & \hf \left( \prod_{j \in \Bbb N} (1 + q^{j -\hf})
                   + \prod_{j \in \Bbb N} (1 - q^{ j - \hf})
             \right) ,   \label{eq_qbasic} \\
   e( - \hL_0) \ch L(\hL_0)
   & = & \hf \left( \prod_{j \in \Bbb N} (1 + e( \epsilon_j))
     + \prod_{j \in \Bbb N} (1 - e( \epsilon_j))
             \right).   \label{eq_basic}
 \end{eqnarray}
 By combining equations (\ref{eq_iden}) and (\ref{eq_basic}) we obtain
 \begin{eqnarray}
   \ch_q L(\dinf; \Lambda)
  &= & \prod_{\alpha \in \Delta_+} \frac{1 - q^{( \Lambda +\rho, \alpha)}}
  {1 - q^{(\hL_0 + \rho, \alpha)}} \cdot   \nonumber \\        
  &&\cdot \hf \left( \prod_{j \in \Bbb N} 
                 \left( 1 + q^{- ( \Lambda + \rho, \epsilon_j)}
                 \right)
                   + \prod_{j \in \Bbb N} 
                 \left( 1 - q^{- (\Lambda + \rho, \epsilon_j)}
                 \right)
            \right)    \nonumber   \\
   & = & \prod_{1 \leq i < j}
      \frac{1 - q^{{}^d \lambda_i - {}^d \lambda_j + j -i} }
           {1 - q^{j -i}} \cdot  
   \prod_{0\leq i <j}
        \frac{1 - q^{2c - {}^d \lambda_{i +1} - {}^d \lambda_{j +1} +j +i } }
             {1 - q^{j +i +1}} \cdot  \nonumber  \\
  & &\cdot \hf
        \left( \prod_{j \in \Bbb N} 
                 \left( 1 + q^{- ( \Lambda + \rho, \epsilon_j)}
                 \right)
                   + \prod_{j \in \Bbb N} 
                 \left( 1 - q^{- (\Lambda + \rho, \epsilon_j)}
                 \right)
       \right)         \label{eq_terms}
 \end{eqnarray}

 A little manipulation of algebra shows that the first product on the 
 right hand side of (\ref{eq_terms}) is the same as the first term
 in (\ref{eq_Dqchar}). We rewrite the second term in (\ref{eq_terms})
 into the product of three terms as follows:
 \begin{eqnarray*}
   \prod_{0\leq i <j}
      \frac{1 - q^{2c +j +i } }{1 - q^{j +i +1}}  \cdot
   \prod_{0\leq i < n_1 <j}
      \frac{1 - q^{2c - {}^d \lambda_{i +1} - {}^d \lambda_{j +1} +j +i } }
           {1 - q^{2c +j +i } }\cdot \\
   \cdot \prod_{0\leq i <j \leq n_1}
      \frac{1 - q^{2c - {}^d \lambda_{i +1} - {}^d \lambda_{j +1} +j +i } }
           {1 - q^{2c +j +i } }.
 \end{eqnarray*}
  A little further manipulation of algebra will show the first,
 second and third terms in the above formula are equal to
 the second, third and fourth terms of (\ref{eq_Dqchar}).
\end{demo}

\end{document}